\documentclass[11pt,a4paper]{article}
\usepackage[T1]{fontenc}
\usepackage{tgpagella}

\usepackage{amsmath} 
\usepackage{amsthm} 
\usepackage{amssymb}	
\usepackage{mathrsfs}
\usepackage{graphicx}
\usepackage{eucal}
\usepackage{esint} 
\usepackage{stmaryrd}
\usepackage[shortlabels]{enumitem}
\usepackage{pdfpages}
\usepackage[textwidth=16cm,centering,headheight=24pt]{geometry} 
\usepackage{comment}
\usepackage{scalerel,stackengine}
\stackMath
\newcommand\reallywidehat[1]{%
\savestack{\tmpbox}{\stretchto{%
  \scaleto{%
    \scalerel*[\widthof{\ensuremath{#1}}]{\kern-.6pt\bigwedge\kern-.6pt}%
    {\rule[-\textheight/2]{1ex}{\textheight}}
  }{\textheight}%
}{0.5ex}}%
\stackon[1pt]{#1}{\tmpbox}%
}
\def\Xint#1{\mathchoice
{\XXint\displaystyle\textstyle{#1}}%
{\XXint\textstyle\scriptstyle{#1}}%
{\XXint\scriptstyle\scriptscriptstyle{#1}}%
{\XXint\scriptscriptstyle\scriptscriptstyle{#1}}%
\!\int}
\def\XXint#1#2#3{{\setbox0=\hbox{$#1{#2#3}{\int}$ }
\vcenter{\hbox{$#2#3$ }}\kern-.6\wd0}}

\def\dashint{\Xint-}

\newcommand{\abs}[1]{\left| #1 \right|} 
\newcommand{\norm}[1]{\left\| #1 \right\|}
 
\newcommand{\R}{\mathbb R}

\newcommand{\N}{\mathbb N}

\renewcommand{\H}{\mathcal{H}}

\renewcommand{\S}{\mathbb{S}}

\newcommand{\eps}{\varepsilon}

\newcommand{\avg}[1]{\left\langle #1 \right\rangle}
\newcommand{\BBB}{\color{blue}}

\newcommand{\res}{\mathbin{\vrule height 1.6ex depth 0pt width
0.13ex\vrule height 0.13ex depth 0pt width 1.3ex}}

\DeclareMathOperator*{\infess}{ess\,inf}                             
\DeclareMathOperator*{\supess}{ess\,sup}      

\DeclareMathOperator{\D}{D}
\let\d\relax
\DeclareMathOperator{\d}{d}

\DeclareMathOperator{\dist}{dist}

\DeclareMathOperator{\diam}{diam}

\DeclareMathOperator{\BV}{BV}
\DeclareMathOperator{\SBD}{SBD}
\DeclareMathOperator{\SBV}{SBV}
\DeclareMathOperator{\GSBV}{GSBV}
\DeclareMathOperator{\GSBD}{GSBD}
\DeclareMathOperator{\Lip}{Lip}

\DeclareMathOperator{\Dev}{Dev}

\makeatletter
\newcommand{\leqnomode}{\tagsleft@true\let\veqno\@@leqno}
\makeatother

\newtheorem{theorem}{Theorem}

\newtheorem{lemma}[theorem]{Lemma}
\newtheorem{prop}[theorem]{Proposition}
\newtheorem{corollary}[theorem]{Corollary}

\newtheorem{mainthm}{Theorem}           

\theoremstyle{definition}
\newtheorem{definition}{Definition}
\newtheorem{remark}[theorem]{Remark}

\theoremstyle{definition}
\newtheorem{step}{Step}

\title{Manifold-constrained free discontinuity problems and Sobolev approximation}
\date{}
\author{Federico Luigi Dipasquale, Bianca Stroffolini}

\newcommand{\Addresses}{{
  \bigskip
  \footnotesize
  Federico Luigi Dipasquale (Corresponding author), \textsc{Dipartimento di Matematica e Applicazioni ``Renato Caccioppoli'', Universit\`{a} degli studi di Napoli ``Federico II'', Via Cintia, Monte S. Angelo, 80126 Napoli, Italy}\par\nopagebreak
  \textit{E-mail address}: \texttt{federicoluigi.dipasquale@unina.it}

    \medskip
    
    Bianca Stroffolini, \textsc{Dipartimento di Matematica e Applicazioni ``Renato Caccioppoli'', Universit\`{a} degli studi di Napoli ``Federico II'', Via Cintia, Monte S. Angelo, 80126 Napoli, Italy}\par\nopagebreak
  \textit{E-mail address}: \texttt{bstroffo@unina.it}
}}

\begin{document}
\maketitle


\begin{abstract}
    We study the regularity of local minimisers of a prototypical free-discontinuity problem 
    involving both a manifold-valued constraint on the maps (which are defined on a bounded domain $\Omega \subset \R^2$) and a variable-exponent growth in the energy functional. To this purpose, we
    first extend to this setting the Sobolev approximation result for special function of bounded variation with small jump set originally proved by Conti, Focardi, and Iurlano \cite{CFI-ARMA, CFI-AIHP} for special functions of bounded deformation. Secondly, we use this extension to prove  regularity of local minimisers.
\end{abstract}

\section{Introduction}
Let $\Omega \subset \R^2$ be a bounded open set, 
$p : \Omega \to [1,+\infty)$ be a measurable function (that will be called a \emph{variable exponent} in the following course), 
and $\mathcal{M}$ be a compact, connected, smooth Riemannian manifold without boundary. 
In this paper, we deal with special functions of bounded variations from $\Omega$ 
into $\mathcal{M}$ whose approximate differential is integrable with respect to the variable exponent $p(\cdot)$
over $\Omega$ and whose jump set has finite $\mathcal{H}^1$-measure. The space of such functions 
will be denoted by the symbol $\SBV^{p(\cdot)}\left(\Omega, \mathcal{M}\right)$; we refer the reader 
to Section~\ref{sec:preliminaries} for its formal definition and as well as for basic material on 
functions that are integrable with respect to variable exponents.

In this work, we first prove an approximation result for maps of class 
$\SBV^{p(\cdot)}(B_\rho, \mathcal{M})$ with small jump set by functions that are also of Sobolev class 
in slightly smaller balls. (Here, $B_\rho$ denotes any ball in $\R^2$.)
Then, we apply such a result to the study of the regularity of local 
minimisers among $\S^{k-1}$-valued maps, where $\S^{k-1}$ denotes the unit sphere in $\R^k$ and $k \geq 2$, of an energy functional with nonstandard growth.
More precisely, we consider, for $u \in \SBV\left(\Omega, \S^{k-1}\right)$, the following energy functional: 
\begin{equation}\label{eq:functional}
    F(u, \Omega) := \int_\Omega \abs{\nabla u}^{p(x)}\,{\d}x + \mathcal{H}^1(J_u \cap \Omega).
\end{equation} 
The functional $F(u, \Omega)$ is finite exactly on $\SBV^{p(\cdot)}\left(\Omega, \S^{k-1}\right)$ and 
and we say that a function $u$ belonging to $\SBV^{p(\cdot)}\left(\Omega, \S^{k-1} \right)$ is a 
\emph{local minimiser (among $\S^{k-1}$-valued maps)} of $F$ if and only if 
    \[
        F(u, \Omega) \leq F(v, \Omega) 
    \]
for every $v \in \SBV\left(\Omega, \S^{k-1}\right)$ such that $\{ u \neq v \} \subset \subset \Omega$.
The functional $F$ may be seen as a prototypical energy involving both a manifold-valued constraint 
on the maps and a nonstandard growth (more precisely, $p(\cdot)$-growth) in the energy functional. 
Note that the functional $F$ reduces to the {\em $p(\cdot)$-energy}
\[
    w \mapsto \int_{\Omega \setminus \overline{J_u}} \abs{\nabla w}^{p(x)}\,{\d}x
\]
out of the closure of the 
jump set $J_u$ of $u$.
Since bounded special functions of bounded variation are of Sobolev class outside $\overline{J_u}$, we see 
that 
any local minimiser of $F$ 
is a local 
minimiser among $\S^{k-1}$-valued maps of the $p(\cdot)$-energy (i.e., a {\em $p(\cdot)$-harmonic map into $\S^{k-1}$}) in the open set $\Omega \setminus \overline{J_u}$. 
However, although $J_u$ is always negligible for 
the Lebesgue measure, there is no reason, {\em a priori}, for $\overline{J_u}$ to not be the whole $\overline{\Omega}$ and, as we shall see in the next paragraphs, proving that $\overline{J_u}$ is ``small'' (more precisely, \emph{essentially closed}) is indeed the most difficult step towards the 
regularity of local minimisers of $F$.

\paragraph{Main results}
In order to state our main results, let us set the basic notation, referring to Section~\ref{sec:preliminaries} for more terminology.
In all this work, unless stated otherwise, the following assumptions on the variable exponent 
$p : \Omega \to [1,+\infty)$ will be in force:
\begin{gather}
    \tag{${\rm p}_1$} p \mbox{ is log-H\"{o}lder continuous}; \\
    \tag{${\rm p}_2$} \mbox{for all } x \in \Omega, \quad 1 < p^- \leq p(x) \leq p^+ < 2, 
\end{gather}
where $p^- := \infess_{x \in \Omega}p(x)$ and $p^+ := \supess_{x \in \Omega}p(x)$. The assumption 
of log-H\"{o}lder conti-nuity is customary in the study of functionals with $p(\cdot)$-growth, c.f., 
eg., \cite{Zhikov3, AcerbiMingione}. We recall this notion along with its geometric meaning 
in Definition~\ref{def:log-Hol} and in Lemma~\ref{lemma:logHol} below. (We are going to comment on the assumption $p^+ < 2$ in the next paragraphs.)

On the manifold $\mathcal{M}$, we assume throughout it is a compact, connected, smooth Riemannian 
manifold without boundary, isometrically embedded in $\R^k$, for some $k \in \N$. 
Under these assumptions, it well-known that there exists a locally smooth retraction 
$\mathcal{P} : \R^k \setminus \mathcal{X} \to \mathcal{M}$, where $\mathcal{X}$ is a smooth complex of 
codimension $2$, with locally $q$-integrable gradient, for every $q \in [1,2)$. 
(See Section~\ref{sec:approx-M} for more details.) This retraction plays a key r\^{o}le in the proof 
of our first main result, Theorem~\ref{thm:A} below.
\begin{mainthm}\label{thm:A}
    Let $p : \Omega \to (1,+\infty)$ be a variable exponent satisfying (${\rm p}_1$), (${\rm p}_2$). 
    Let $\mathcal{M}$ be a compact, connected, smooth Riemannian manifold without boundary, 
    isometrically embedded in $\R^k$, for some $k \in \N$. There 
    exist universal constants $\xi$, $\eta > 0$ 
    such that for any $s \in (0,1)$ and any $u \in \SBV^{p(\cdot)}\left(B_\rho, \mathcal{M}\right)$ (where $B_\rho$ is any ball in $\R^2$) satisfying
    \begin{equation}\label{eq:small-jump-set-intro}
        \mathcal{H}^1\left(J_u \cap B_\rho\right) < \eta(1-s) \rho/2,
    \end{equation}
    there exists a function 
    $w \in W^{1,p(\cdot)}\left(B_{s\rho}, \mathcal{M}\right) \cap \SBV^{p(\cdot)}\left(B_{\rho}, \mathcal{M}\right)$ and a family $\mathcal{F}$ of balls such that the following holds. 
    The function $w$ coincides with $u$ a.e. outside of the union of the balls in the family $\mathcal{F}$. Such a union is contained in $B_{(1+s)\rho/2}$ and controlled in measure and perimeter by 
    $\xi$, $\eta$, 
    $\rho$, and $\mathcal{H}^1\left(J_u \cap B_\rho\right)$. 
    Moreover, $w$ has less jump than $u$ in 
    $B_\rho$ and $\int_{B_\rho} \abs{\nabla w}^{p(x)}\,{\d}x$ is controlled by $\rho$, 
    $p^-$, $p^+$, the log-H\"{o}lder constant of $p(\cdot)$, $\mathcal{M}$, $\mathcal{P}$, $k$, and 
    $\norm{\nabla u}_{L^{p{(\cdot)}}(B_\rho)}$.
\end{mainthm}
Theorem~\ref{thm:A} is an abridged version of a more precise and descriptive statement, 
Theorem~\ref{thm:approx-sbv-px} below. Theorem~\ref{thm:A} may be seen as an extension to our 
manifold-valued, variable-exponent setting of a result by Conti, Focardi, and Iurlano \cite{CFI-ARMA}, 
which was developed in the context of planar domains, constant exponents, and (unconstrained) $\SBD^p$ functions 
(i.e., for special functions of bounded deformation with approximate symmetric gradient in $L^p$ and 
jump set with finite $\mathcal{H}^{1}$-measure) 
 with small jump set, and later applied in \cite{CFI-AIHP} to the study of 
Griffith's type brittle fracture functionals. As we shall see in more detail later, the result is confined to the two dimensional setting and cannot be directly extended to higher dimensions.

The seminal idea for studying the regularity of local minimisers of free discontinuity problems for 
maps of class $\SBV$ is due to De~Giorgi, Carriero, and Leaci \cite{DeGiorgiCarrieroLeaci} and lies in 
showing that the jump set of any local minimiser $u$ is \emph{essentially closed}, i.e., satisfies
\[
    \mathcal{H}^1\left( \Omega \cap \left( \overline{J_u} \setminus J_u \right) \right) = 0.
\]
In the constant exponent case, once this is done, standard elliptic regularity yields 
that 
\[
    u \in C^1\left( \Omega \setminus \overline{J_u}\right).
\]
These ideas, firstly developed for scalar-valued functions in \cite{DeGiorgiCarrieroLeaci}, 
were extended to the case of $\S^{k-1}$-valued maps (and constant $p$) in \cite{CarrieroLeaci}. 
Very recently, they have been 
adapted to the scalar-valued, variable-exponent setting in~\cite{LScSoV} (for a larger class of convex functionals, of which~\eqref{eq:functional} is the prototype). In all these works, a major 
technical tool is the Sobolev-Poincar\'{e} inequality for $\SBV$-functions due, once again, to 
De~Giorgi, Carriero, and Leaci \cite{DeGiorgiCarrieroLeaci}. Here, we follow the approach 
devised by Conti, Focardi, and Iurlano \cite{CFI-ARMA, CFI-AIHP}, which avoids the use of 
truncations (fundamental in the proof of the Sobolev-Poincar\'{e} inequality 
in~\cite{DeGiorgiCarrieroLeaci}), by relying, instead, on the Sobolev 
approximation. With the aid of Theorem~\ref{thm:A}, we prove Theorem~\ref{thm:B} below, 
which is our second main result in this paper. Before stating the theorem, we have 
to introduce a strengthening on the assumption (${\rm p}_1$): 
\begin{equation}
    \tag{${\rm p}_1'$} p \mbox{ is strongly log-H\"{o}lder continuous},
\end{equation}
see Definition~\ref{def:strongly-log-hol} and Remark~\ref{rk:stong-log-Hol} about this condition 
and its r\^{o}le in the proof of Theorem~\ref{thm:B}.

\begin{mainthm}\label{thm:B}
    Let $\Omega \subset \R^2$ be a bounded, open set. Let $p(\cdot)$ be a variable exponent satisfying (${\rm p}_1'$), (${\rm p}_2$). Assume 
    $u \in \SBV^{p(\cdot)}\left(\Omega, \S^{k-1}\right)$ is any local minimiser of the 
    functional $F(\cdot,\Omega)$ defined by~\eqref{eq:functional}. 
    Then, the jump set $J_u$ of $u$ is essentially closed, i.e., 
    $\mathcal{H}^1\left(\Omega \cap \left(\overline{J_u} \setminus J_u\right)\right) = 0$. 
    Moreover, if in addition $p \in C^{0,\alpha}(\Omega)$ for some 
    $\alpha \in (0,1]$, then 
    there exists a relatively open set $\Omega_0 \subset \Omega \setminus \overline{J_u}$ and $\beta_0 \in (0,1)$, with $\mathcal{H}^{2-p^-}\left( \left( \Omega \setminus \overline{J_u} \right) \setminus \Omega_0\right) = 0$, such that $u \in C^{1,\beta_0}_{\rm loc}\left(\Omega_0, \S^{k-1} \right)$. The number $\beta_0$ depends only on $k$, $\mathcal{M}$, $\mathcal{P}$, $p^-$, $p^+$, $[p]_{0,\alpha}$, $\alpha$.
\end{mainthm}
In particular, $u \in C^1\left(\Omega_0 \setminus \overline{J_u}, \S^{k-1}\right)$, which is (almost) precisely 
the result expected from \cite{DeGiorgiCarrieroLeaci, CarrieroLeaci, CFI-AIHP} (in \cite{CFI-AIHP}, 
local $C^{1,\beta_0}$-regularity outside of $\overline{J_u}$ was proven in two dimension even in the $\SBD^p$ setting). In the particular case where $p$ is constant, we can take $\Omega_0 = \Omega$ in Theorem~\ref{thm:B} (see Remark~\ref{rk:regularity-p-const}), recovering classical results in \cite{CarrieroLeaci} {\em via} a different technique.

Before illustrating the ideas of the proofs of Theorem~\ref{thm:A} and Theorem~\ref{thm:B}, let us 
motivate the interest towards our results.

\paragraph{Background and motivations}
In recent years, there have been an incredible amount of papers accounting, under different perspectives, for approximations of functions of special bounded variation or even special bounded deformation \cite{CorToa, Chambolle, ChambolleContiFrancfort, FuDePra, CFI-ARMA, ChambolleContiIurlano, Crismale, CCS}. 

The arising of such variety of approximation results is due to their multiple applications in many 
different problems in the calculus of variations, 
connected (especially but not exclusively) with image segmentation 
and fracture mechanics. 
The most known examples are probably the  
variational problems associated with the nowadays classical 
Mumford-Shah functional and the Griffith static one. 

Among the above-quoted works, the most relevant to our purposes are \cite{CFI-ARMA, CFI-AIHP, CCS}.
In particular, in \cite{CFI-ARMA, CFI-AIHP}, 
relying on the Sobolev approximation, the authors were able to prove integral representation results 
and existence of strong minimisers for Griffith's functional (actually, for a wider class of 
functionals) defined over $\GSBD^p(\Omega)$, where $\Omega$ is a two-dimensional domain. The Sobolev 
approximation 
result has been extended to any dimension in \cite{CCS}. 
However, in~\cite{CCS}, differently from \cite{CFI-ARMA} and Theorem~\ref{thm:A} above, 
a small, exceptional set whose perimeter and volume are controlled by the size of the jump has to be 
removed from the domain. 
It is {\em a priori} unclear whether removing this small set is really necessary or not.  
Here, in Appendix~\ref{app:counterex} we show that 
there is no counterpart, in dimension higher than 2, 
of the construction of the Sobolev approximation in \cite{CFI-ARMA} and in the present paper
under assumption~\eqref{eq:small-jump-set-intro} alone. 
This suggests that the construction in \cite{CCS} is, in fact, optimal.

In the calculus of variations, energy functionals with $p(\cdot)$-growth in the gradient 
have been proposed in the modelling of materials which exhibit a strongly 
anisotropic behaviour
starting from the works~\cite{Zhikov1, Zhikov2}. 
In the Sobolev setting, the regularity of minimisers has been analysed to a certain degree of 
generality, both in the scalar and in the vector-valued, unconstrained case, see, e.g., \cite{CosciaMingione,AcerbiMingione}. 
Less literature is available for manifold-constrained maps. However,  
the regularity problem for $p(\cdot)$-harmonic maps has been recently considered in \cite{DeFilippis, ChlebickaDeFilippisKoch}.

In the last years, several works \cite{DeCLV, AlmiReggianiSolombrino, ScSoStr, LScSoV} 
have undertaken the study of energy functionals with nonstandard growth defined over spaces of maps 
of class $\SBV^{p(\cdot)}$ or even $\GSBV^{p(\cdot)}$ or $\GSBV^{\psi}$. 
The aims of these works range from lower semicontinuity results \cite{DeCLV, AlmiReggianiSolombrino} to 
integral representation theorems~\cite{ScSoStr}, up to regularity in the scalar-valued 
case~\cite{LScSoV}.

Spaces of functions of (special) bounded variation and values into a Riemannian 
manifold $\mathcal{M}$ have been recently studied in \cite{IgnatLamy, CanevariOrlandi}. 
The authors of~\cite{IgnatLamy, CanevariOrlandi} are mainly concerned with constructing liftings of 
such mappings from the manifold $\mathcal{M}$ to its universal cover $\widetilde{\mathcal{M}}$.
The interest towards such problem was initially stimulated by applications to the Landau-de Gennes theory of 
liquid crystals and to the Ginzburg-Landau model of superconductivity.  

The present work is, to the best of our knowledge, the first that considers the prototypical energy 
functional $F$ in~\eqref{eq:functional} for manifold-valued special functions of bounded variation 
with $L^{p(\cdot)}$-integrable approximate differential.
The functional $F$ in~\eqref{eq:functional} is itself a particular case of a more 
general one:
\begin{equation}\label{functional}
    \mathcal{F}(u, \Omega) := \int_\Omega  \abs{\nabla u}^{p(x)}\,{\d}x 
    + \int_{J_u} \phi_0(u^+,u^-, \nu) \, {\d}\H^1
\end{equation}
where $u\in \SBV^{p(x)}(\Omega, \mathcal{M})$ 
and $\phi_0$ is $\BV$-elliptic, see \cite[Definition~5.13]{AFP}, and bounded away from zero.  
Functionals like these appear in the theory of liquid crystals, for example nematics and, in this case, the manifold is isomorphic to $\S^1$. Another model, related to smectic thin films, can be found in a recent paper \cite{BallCanevariStro}. Here a free discontinuity problem is proposed in order to describe surface defects in a smectic thin film.  The free energy functional contains an interfacial energy,
which penalises dislocations of the smectic layers at the jump. The function space is a subspace of function of bounded variation $\SBV^2(\Omega, \mathbb{R}^2)$ with values in a suitable manifold $\mathcal{N}$.

\paragraph{Proofs of the main results: a sketch}
The proof of Theorem~A (better, of Theorem~\ref{thm:approx-sbv-px}) is contained in Section~\ref{sec:approx}. Essentially, it proceeds in two steps: 
first, we prove an analogous statement for unconstrained maps with values into $\R^k$ (Proposition~\ref{thm:approx-sbv-px-Rk}), following very 
closely the original argument in \cite[Theorem~2.1 and Proposition~3.2]{CFI-ARMA} and 
exploiting the assumption of log-H\"{o}lder continuity of the variable exponent. Then, we use the 
aforementioned retraction $\mathcal{P}$ to obtain Theorem~\ref{thm:A} from its unconstrained 
counterpart, by retraction onto the manifold of the image of the unconstrained approximating maps. 
The restriction on the dimension of domain in 
Theorem~\ref{thm:A} comes from the fact that, to craft the Sobolev approximation, we use the same 
construction as in \cite{CFI-ARMA}, 
which is strictly two dimensional and cannot be extended (without heavy modifications) to higher dimensions (see \cite{CCS} and 
Remark~\ref{rk:obstruction} and Appendix~\ref{app:counterex} below for more details on this point). 
The restriction $p^+ < 2$ in (${\rm p}_2$) is due 
instead to the usage of the retraction $\mathcal{P}$ in the proof of Theorem~\ref{thm:A} 
along with the fact that we merely require connectedness on $\mathcal{M}$ 
(so to allow for $\mathcal{M}$ to be a circle, an important case in potential 
applications -- for instance, to liquid crystals), see 
Lemma~\ref{lemma:retraction} and Remark~\ref{rk:integrability-proj} for more details. 
However, since we work in dimension 2, such a restriction is not a dramatic 
drawback, in the sense that this is the subcritical regime for Sobolev-Morrey's embedding 
and maps of class $W^{1,p(\cdot)}\left(\Omega, \mathcal{M}\right)$ are not automatically continuous
(neither in $\Omega$ nor in open subsets of $\Omega$ with positive measure). 
In other words, this is the regime in which all the essential complications in the study of 
the regularity of local minimisers of the functional $F$ in~\eqref{eq:functional} already show up, 
without the additional ones due to possibly large oscillations of the variable exponent. 

As alluded few lines above, 
the assumption of log-H\"{o}lder continuity of the variable exponent is extremely important in 
the proof of Theorem~\ref{thm:A} and this can be easily realised by looking at its geometric meaning 
(see Lemma~\ref{lemma:logHol}). Indeed, roughly speaking, such an assumption boils down to the 
possibility of locally ``freezing'' the variable exponent, up to a controlled error. 
Joint with the very precise estimates for the approximating Sobolev map coming from \cite{CFI-ARMA}, 
this yields a rather direct extension of the results of \cite{CFI-ARMA} to our 
variable-exponent setting, at least in the unconstrained case 
(compare the proofs of Proposition~\ref{thm:local-approx} and Proposition~\ref{thm:approx-sbv-px-Rk} 
with those of \cite[Theorem~2.1 and Proposition~3.2]{CFI-ARMA}).

The proof of Theorem~\ref{thm:B} proceeds in various steps, following the pioneering approach 
from \cite{DeGiorgiCarrieroLeaci, CarrieroLeaci}. The key point relies in proving a suitable 
power-decay of the energy in small balls with the radius of the ball. 
This is done in Theorem~\ref{thm:decay-lemma} by means of a classical argument by contradiction 
and a blow-up analysis, relying on assumption (${\rm p}_1'$). 
We adapt the strategy of \cite{CFI-AIHP} to exploit, in the 
blow-up analysis, Theorem~\ref{thm:A} (with $\mathcal{M} = \S^{k-1}$ and $k \geq 2$) instead of the 
classical Sobolev-Poincar\'{e} inequality in $\SBV$ (adapted to the variable-exponent framework in 
\cite{ScSoStr, LScSoV}).   


Once the decay lemma is obtained, another classical argument originating from 
\cite{DeGiorgiCarrieroLeaci} yields suitable density lower bounds for $F\left(u, B_\rho(x)\right)$, 
where $u$ is a local minimiser of $F$, 
$x \in \overline{J_u}$, and $\rho$ is small. In turn, such density lower bounds readily 
imply, by a standard argument in geometric measure theory, that $\overline{J_u}$ is 
essentially closed. This step requires more than log-H\"{o}lder continuity and indeed
we ask {\em strong} log-H\"{o}lder continuity in the statement of Theorem~\ref{thm:B}. 
To the purpose of proving essential closedness of the jump set, strong log-H\"older continuity 
turns out to be enough and, actually, also local minimality can be weakened to {\em quasi-minimality} 
(see Definition~\ref{def:quasi-min}). 
The full strength of H\"{o}lder continuity and of local minimality are instead 
needed to prove that 
$u \in C^{1,\beta_0}_{\rm loc}\left(\Omega_0, \S^{k-1}\right)$, where $\Omega_0$ is as in the 
statement of Theorem~\ref{thm:B}. 
Indeed, here we use the regularity result 
\cite[Theorem~1]{DeFilippis} (in the simplified form provided by Theorem~\ref{thm:defilippis} below) 
for the $p(\cdot)$-harmonic energy for manifold-valued maps, i.e., for local minimisers (among 
compactly supported perturbations) of the functional
\[
    W^{1,p(\cdot)}\left(\Omega, \mathcal{M}\right) \ni w \mapsto \int_\Omega \abs{\nabla w}^{p(x)}\,{\d}x, 
\]
which requires $p \in C^{0,\alpha}(\Omega)$, for some $\alpha \in (0,1]$, among its assumptions.

\paragraph{Organisation of the paper}
In Section~\ref{sec:preliminaries}, we establish notation and recall the basics facts about spaces of 
functions integrable with respect to variable exponents. 
In Section~\ref{sec:approx}, we prove Theorem~\ref{thm:A}. 
In Section~\ref{sec:existence}, we prove Theorem~\ref{thm:B}. The paper is completed by a series of 
appendices containing mostly technical material for which we were not able to find explicit proof in 
the literature. In Appendix~\ref{app:counterex}, we exhibit an example which shows that the 
approximation procedure in \cite{CFI-ARMA} and, in turn, in Section~\ref{sec:approx} cannot work in 
higher dimensional domains under the mere assumption of $\mathcal{H}^1$-smallness of the jump set.

\numberwithin{equation}{section}
\numberwithin{definition}{section}
\numberwithin{theorem}{section}

\section{Preliminaries}\label{sec:preliminaries}

\subsection{Notation}\label{sec:notation}

\begin{enumerate}[(i)]
    \item We use the symbol $\{x_n\}$ to denote a sequence, indexed by $n \in \N$, of elements $x_n$ of a certain 
set $E$. Usually, for the sake of a lighter notation, we do not relabel subsequences. 

\item
In inequalities like $A \lesssim B$, the symbol $\lesssim$ means that there exists a constant 
$C$, independent of $A$ and $B$, such that $A \leq CB$.
\item
We denote $B_\rho^n\left(x_0\right) := \left\{ x \in \R^n : \abs{x-x_0} < \rho \right\}$ the open ball 
of radius $\rho$ and centre $x_0$ in $\R^n$. Since we work almost exclusively in dimension $n=2$, we 
drop the superscript ``2'' for balls in $\R^2$. Often, the centre of the ball will be irrelevant and we shall omit it from the notation. Given a ball $B_\rho$, the symbol $B_{s\rho}$ denotes 
the ball concentric with $B_{\rho}$ and radius dilated by the factor $s > 0$.
\item
We denote by ${\mathcal{M}}$ a compact, connected, smooth Riemannian manifold without 
boundary, of dimension $m \geq 1$. 
Without loss of generality, we may always view at $\mathcal{M}$ as a compact, connected, 
$m$-dimensional smooth submanifold in $\R^k$, for some $k \in \N$. This can always be achieved, for instance, by means of Nash's isometric embedding theorem. If not specified otherwise, we will always assume to embed $\mathcal{M}$ in $\R^k$ via Nash's embedding. 
By compactness of $\mathcal{M}$, we can find $M > 0$, depending only on $\mathcal{M}$ and the choice 
of the isometric embedding, such that 
\begin{equation}\label{eq:M-in-a-ball}
    \mathcal{M} \subset \mathcal{B}^k_M := \left\{ y \in \R^k : \abs{y} < M \right\}.
\end{equation}

In the following course, when saying that a quantity ``depends on $\mathcal{M}$'', we will always mean 
that it depends on $\mathcal{M}$ and the chosen isometric embedding of $\mathcal{M}$ 
into $\R^k$. However, by compactness of $\mathcal{M}$, the choice of the isometric embedding is essentially irrelevant, in the sense that changing the embedding can result, at worst, in enlarging 
$k$ and the constants that depend on $\mathcal{M}$.
\item
In the special case $\mathcal{M}$ is a sphere, we always look at it 
as a submanifold of $\R^{m+1}$ the obvious way. In particular, we will denote 
$\S^{k-1}_t := \left\{ x \in \R^k : \abs{x} = t \right\}$ the $(k-1)$-dimensional 
sphere in $\R^k$ of radius $t > 0$ and centre the origin, endowed with the canonical metric. 
We set $\S^{k-1} := \S^{k-1}_1$. 
The canonical orthonormal basis of $\R^k$ will be denoted 
$\left\{ {\bf e}_1, \dots, {\bf e}_k \right\}$. 
\item
For any $\delta > 0$, we will denote 
\begin{equation}\label{eq:Omega-delta}
    \Omega_\delta := \left\{ x \in \Omega : \dist(x, \partial \Omega) > \delta \right\}.    
\end{equation}
the portion of the $\delta$-neighborhood of $\partial \Omega$ that lies in the interior of $\Omega$.
\end{enumerate}

\subsection{Variable exponent Lebesgue}\label{sec:Lp(x)}
In this section, we recall some basic facts about variable-exponent Lebesgue spaces. 
The reader can consult the monographs \cite{DHHR, C-UF} for more details.

A measurable function $p : \Omega \to [1,+\infty)$ will be called a 
\emph{variable exponent}. For every $A \subset \Omega$ (measurable and not 
empty) we define 
\[
	p^+_A := \supess_{x\in A} p(x) \quad \mbox{and} \quad 
	p^-_A := \infess_{x \in A} p(x).
\]
Of course, we can take $A = \Omega$, and in this case we write $p^+$ and $p^-$ 
in place of $p^+_\Omega$ and $p^-_\Omega$, respectively. We say that a variable 
exponent $p$ is {\em bounded} if $p^+ < +\infty$.  

The \emph{modular} $\varrho_{p(\cdot)}(u)$ of a measurable function 
$u : \Omega \to \R^m$ with respect to the variable exponent $p(\cdot)$ is 
defined by
\[
	\varrho_{p(\cdot)}(u) := \int_\Omega |u(x)|^{p(x)}\,{\d}x
\]
and the \emph{Luxembourg norm} (henceforth, simply \emph{norm}) of $u$ by
\[
	\|u\|_{L^{p(\cdot)}(\Omega)} := 
	\inf\{ \lambda > 0 : \varrho_{p(\cdot)}(u/\lambda)\leq 1\}.
\]
The \emph{variable exponent Lebesgue space} $L^{p(\cdot)}\left(\Omega, \R^k\right)$ is 
defined as the set of measurable functions $u : \Omega \to \R^k$ such that 
there exists $\lambda > 0$ so that $\varrho_{p(\cdot)}(u/\lambda) < +\infty$. 

We collect the properties of variable exponent Lebesgue space that we will use 
in the following propositions.
\begin{prop}
Let $p : \Omega \to [+1,\infty)$ be a variable exponent and assume that 
$p^+ < +\infty$.
Then:
\begin{enumerate}[(i)]
	\item $L^{p(\cdot)}\left(\Omega,\R^k\right)$ coincides with the set 
of measurable functions $u : \Omega \to \R^k$ such that $\varrho_{p(\cdot)}(u)$ is 
finite.
	\item $\|\cdot\|_{L^{p(\cdot)}(\Omega)}$ is a norm on 
	$L^{p(\cdot)}\left(\Omega, \R^k \right)$, under which $L^{p(\cdot)}\left(\Omega, \R^k\right)$ is a 
	Banach space.
	\item If $p^- > 1$, then $L^{p(\cdot)}\left(\Omega,\R^k\right)$ is uniformly convex and its 
	dual space is isomorphic to $L^{p'(\cdot)}\left(\Omega, \R^k\right)$, where $p'$ is the 
	variable exponent satisfying $\frac{1}{p} + \frac{1}{p'} = 1$ a.e. in $\Omega$.
	\item The following inequalities hold:
	\begin{equation}\label{eq:modular-norm}
	\begin{cases}
		\varrho_{p(\cdot)}(u)^{\frac{1}{p^+}} \leq \|u\|_{L^{p(\cdot)}(\Omega)} 
		\leq \varrho_{p(\cdot)}(u)^{\frac{1}{p^-}}, & \mbox{if } \|u\|_{L^{p(\cdot)}(\Omega)} > 1, \\
		\varrho_{p(\cdot)}(u)^{\frac{1}{p^-}} \leq \|u\|_{L^{p(\cdot)}(\Omega)} 
		\leq \varrho_{p(\cdot)}(u)^{\frac{1}{p^+}}, & \mbox{if } \|u\|_{L^{p(\cdot)}(\Omega)} \leq 1,
	\end{cases}
	\end{equation}
\end{enumerate}
\end{prop}

Modulars satisfy appropriate versions of Fatou's lemma, of the monotone 
convergence theorem, and of the dominated convergence theorem. Here we give a particular 
case, enough for our purposes, of a much more general result, 
c.~f. \cite[Lemma~3.2.8]{DHHR}.
\begin{prop}\label{prop:convergence-modulars}
    Let $p : \R^n \to [1,+\infty)$ be a variable exponent and let $\{f_h\}$ and $f$ 
    be measurable functions from $\Omega$ into $\R^k$, 
    for some $k \in \N$. Then,
    \begin{enumerate}[(i)]
        \item If $f_h \to f$ $\mathcal{L}^n$-almost everywhere in $\Omega$ as $h \to +\infty$, then 
        $\varrho_{p(\cdot)}(f) \leq \displaystyle\liminf_{h\to+\infty} \varrho_{p(\cdot)}(f_h)$. 
        \item If $\abs{f_h} \nearrow \abs{f}$ $\mathcal{L}^n$-almost everywhere in $\Omega$ as 
        $h \to +\infty$, then 
        $\varrho_{p(\cdot)}(f) = \displaystyle\lim_{h \to +\infty} \varrho_{p(\cdot)}(f_h)$.
        \item If $f_h \to f$ $\mathcal{L}^n$-almost everywhere in $\Omega$ as $h \to +\infty$ 
        and there exists $g \in L^{p(\cdot)}\left(\Omega, \R^k\right)$ such that 
        $\abs{f_h} \leq \abs{g}$ a.e., then 
        $f_h \to f$ in $L^{p(\cdot)}\left(\Omega, \R^k\right)$ as $h \to +\infty$.
    \end{enumerate}
\end{prop}

The following proposition extends to the variable exponent setting the classical 
embedding property of Lebesgue spaces on sets with finite (Lebesgue) measure.
\begin{prop}[{\cite[Corollary~3.3.4]{DHHR}}]\label{prop:embedding}
	Let $p, q : \Omega \to [1, +\infty)$ be bounded variable exponents on a 
	bounded set $\Omega \subset \R^n$. Then 
	$L^{p(\cdot)}(\Omega) \hookrightarrow L^{q(\cdot)}(\Omega)$ if and only if 
	$q \leq p$ a.e. in $\Omega$. The embedding constant is less than or equal 
	to the minimum between $2\left(1+ \mathcal{L}^n(\Omega)\right)$ and 
	$2 \max\left\{ \mathcal{L}^n(\Omega)^{\left(\frac{1}{q} - \frac{1}{p}\right)^{+}}, \mathcal{L}^n(\Omega)^{\left(\frac{1}{q} - \frac{1}{p}\right)^{-}} \right\}$.
\end{prop}

The following definition specify a \textbf{crucial} quantitative notion of continuity for 
variable exponents, called {\em log-H\"{o}lder continuity}. This condition was firstly introduced 
in~\cite{Zhikov3} as a way to avoid the Lavrentiev phenomenon for minimisers of
variational integrals and suddenly became customary in this context (c.f., e.g., \cite{AcerbiMingione}).
\begin{definition}\label{def:log-Hol}
We say that a variable exponent $p : \Omega \to [1,+\infty)$ satisfies the 
{\em log-H\"{o}lder condition} in $\Omega$ if and only if 
\begin{equation}\label{eq:log-Hol-def}
	\exists C_p > 0 : \,\,
	\forall x,y \in \Omega : 0 < \abs{x-y} \leq \frac{1}{2},\qquad
	\abs{p(x) - p(y)} \leq \frac{C_p}{- \ln\abs{x-y}}.
\end{equation}
Equivalently, $p : \Omega \to [1,+\infty)$ is log-H\"{o}lder continuous if 
and only if its modulus of continuity $\omega$ satisfies
	\[
		\limsup_{\rho \to 0} \omega(\rho) \log\left( \frac{1}{\rho} \right) < +\infty.
	\]
We call $C_p$ the \emph{log-H\"{o}lder constant} of $p$.
\end{definition}

The following lemma, an abridged version of \cite[Lemma~4.1.6]{DHHR} which 
is sufficient to our purposes, illustrates the geometrical meaning of 
log-H\"{o}lder continuity.
\begin{lemma}\label{lemma:logHol}
	Let $p : \Omega \to [1,+\infty)$ be a bounded, continuous variable exponent. 
	The following conditions are equivalent:
	\begin{enumerate}[(i)]
		\item $p$ is log-H\"{o}lder continuous;
		\item There exists a constant $\ell > 0$ such that, for all open balls 
		$B \subset \Omega$, there holds
		\begin{equation}\label{eq:log-hol}
			\mathcal{L}^n(B)^{p^-_B - p^+_B} \leq \ell.
		\end{equation}
	\end{enumerate}
\end{lemma}
\begin{remark}\label{rk:extension-p} 
A bounded, log-H\"{o}lder continuous variable exponent defined on a bounded set $\Omega \subset \R^n$ 
can always be extended to a bounded, log-H\"{o}lder continuous variable exponent $q$ which is 
defined on the whole of $\R^n$ and that satisfies $C_q \leq \left(p^-\right)^2C_p$, 
$\ell_q \leq \left(p^-\right) \ell$, $q^- = p^-$,  
$q^+ = p^+$, and~\eqref{eq:log-hol} for any open ball $B\subset \R^n$ see, e.g., 
\cite[Proposition~4.1.7]{DHHR}. For our purposes in this paper, there is no loss of generality in assuming from the very beginning 
that $p$ is defined on the whole $\R^n$.
\end{remark}

Later in this work we will need the following strengthening of the log-H\"{o}lder 
condition~\eqref{eq:log-Hol-def}.
\begin{definition}\label{def:strongly-log-hol}
We say that a variable-exponent $p(\cdot) : \Omega \to [1,+\infty)$ satisfies the 
{\em strong log-H\"{o}lder condition} in $\Omega$ if and only if its modulus of 
continuity $\omega$ satisfies
\begin{equation}\label{eq:strong-log-hol}
	\limsup_{\rho \to 0} \omega(\rho) \log\left( \frac{1}{\rho} \right) = 0.
\end{equation}
\end{definition}

\begin{remark}
    The extension to the whole $\R^n$ of a strongly log-H\"{o}lder continuous variable exponent 
    is still strongly log-H\"{o}lder continuous.
\end{remark}

\subsection{Variable exponent Sobolev spaces and $p(\cdot)$-harmonic maps}
Let $\Omega \subset \R^n$ be an open set and $p(\cdot) : \R^n \to \R$ be a variable exponent. As in the classical case, we define
\[
    W^{1,p(\cdot)}\left(\Omega, \R^k\right) := \left\{ u \in \left(W^{1,1} \cap L^{p(\cdot)}\right)\left(\Omega, \R^k \right) : \nabla u \in L^{p(\cdot)}\left(\Omega,\R^{k\times n}\right) \right\}.
\]
The norm 
\[
    \norm{u}_{W^{1,p(\cdot)}(\Omega)} := \norm{u}_{L^{p(\cdot)}(\Omega)} 
    + \norm{\nabla u}_{L^{p(\cdot)}(\Omega)},
\]
turns $W^{1,p(\cdot)}\left(\Omega,\R^k\right)$ into a Banach space, separable if $p^+ < +\infty$ and 
reflexive if $p^- > 1$ and $p^+ < +\infty$. 
We address the reader to~\cite{DHHR, C-UF} for full details about variable exponent Sobolev spaces.

As in the classical case, we define  
\[
    W^{1,p(\cdot)}\left(\Omega, \mathcal{M}\right) := \left\{ W^{1,p(\cdot)}\left(\Omega, \R^k\right) : u(x) \in \mathcal{M} \,\,\mbox{for a.e. } x \in \Omega \right\}.
\]

\begin{remark}\label{rk:traces}
Assume $\Omega \subset \R^n$ is a bounded, Lipschitz domain and $p$ is a bounded, log-H\"older 
continuous variable exponent satisfying $p^- > 1$. 
Then, the trace $u \vert_{\partial \Omega}$ of a function 
$u$ belonging to $W^{1,p(\cdot)}\left(\Omega, \R^k\right)$ is well-defined and it is, in particular, 
an element 
of $L^{p(\cdot)}\left(\partial\Omega,\R^k\right)$. Moreover, since $u \in W^{1,p(\cdot)}\left(\Omega,\R^k\right)$ implies, in particular, that $u$ belongs to $W^{1,1}\left( \Omega, \R^k\right)$, 
the trace $u \vert_{\partial \Omega}$ is specified at least $\mathcal{H}^{n-1}$-a.e. on 
$\partial \Omega$, see \cite[Section~12.1]{DHHR} for details.  
Therefore, it makes sense to say that $u \vert_{\partial \Omega}$ takes values in $\mathcal{M}$, by saying that $u \vert_{\partial \Omega}$ takes values in $\mathcal{M}$ if and only if $u \vert_{\partial \Omega}(y) \in \mathcal{M}$ for $\mathcal{H}^{n-1} \res \partial \Omega$-a.e. 
$y \in \partial \Omega$.
\end{remark}   

Next, following, e.g., \cite{DeFilippis}, we define a \emph{$p(\cdot)$-harmonic map} as a local 
minimiser, with respect to compactly supported perturbations, of the functional
\begin{equation}\label{eq:p(x)-energy}
    W^{1,p(\cdot)}\left( \Omega, \mathcal{M} \right) \ni w \mapsto \int_\Omega \abs{\nabla w}^{p(x)}\,{\d}x.
\end{equation}
The following partial regularity theorem for $p(\cdot)$-harmonic maps is a simplified version of 
a more general statement proven in~\cite{DeFilippis}. (We use the symbol $[x]$ to denote the integer part of $x$.)
\begin{theorem}[{\cite[Theorem~1]{DeFilippis}}]\label{thm:defilippis}
    Assume $p : \R^n \to (1,+\infty)$ is a variable exponent satisfying
    \begin{gather*}
        p \in C^{0,\alpha}(\R^n), \quad \mbox{for some }\alpha \in (0,1],\\
        1 < p^- \leq p(x) \leq p^+ < +\infty \quad \mbox{for all } x \in \Omega.
    \end{gather*}
    Let $\mathcal{M}$ be a compact, $([p^+]-1)$-connected, smooth Riemannian manifold, without boundary.
    Let $u \in W^{1,p(\cdot)}\left(\Omega,\mathcal{M}\right)$ be a $p(\cdot)$-harmonic map. Then, 
    there exists a relatively open set $\Omega_0 \subset \Omega$ such that $u \in C^{1,\beta_0}_{\rm loc}\left(\Omega_0, \mathcal{M}\right)$ for some $\beta_0 \in (0,1)$ depending only on $n$, $k$, 
    $\mathcal{M}$, $p^-$, $p^+$, $[p]_{0,\alpha}$, $\alpha$. Moreover, 
    $\mathcal{H}^{n-p^-}\left(\Omega \setminus \Omega_0 \right) = 0$.
\end{theorem}

\begin{remark}\label{rk:defilippis}
    We recall that, given an integer $j \geq 0$, a manifold $\mathcal{M}$ is said to be 
    \emph{$\ell$-connected} iff its first $j$ homotopy groups vanish identically, that is, iff 
    \[
        \pi_0(\mathcal{M}) = \dots = \pi_j(\mathcal{M}) = 0.
    \]
    We observe that, although \cite[Lemma~4]{DeFilippis} (which in turn is partly borrowed 
    from~\cite[Lemma~4.5]{Hopper}) is stated assuming $j \geq 1$ (i.e., 
    $\mathcal{M}$ simply connected), it holds as well for $j = 0$ 
    (c.f., e.g., \cite{HardtLin, BousquetPonceVanSchaftingen, CanevariOrlandi-TopI} and 
    Lemma~\ref{lemma:retraction} below). 
    This fact allows us to construct maps with all the properties required 
    by \cite[Lemma~5]{DeFilippis} without assuming $\mathcal{M}$ simply connected, c.f., e.g. Lemma~\ref{lemma:projection} below. Besides \cite[Lemma~4 and Lemma~5]{DeFilippis}, the assumption 
    that $\mathcal{M}$ is simply connected is never used in the proof of \cite[Theorem~1]{DeFilippis} 
    (and indeed removed in the follow up \cite{ChlebickaDeFilippisKoch} of \cite{DeFilippis}). 
    Consequently, \cite[Theorem~1]{DeFilippis} holds as well even if $\mathcal{M}$ is merely connected, 
    provided that $p^+ < 2$. In turn, Theorem~\ref{thm:defilippis} holds if $\mathcal{M}$ is merely 
    connected, provided that $p^+ < 2$, which are precisely the assumptions under which we work in this 
    paper.
\end{remark}

\subsection{The space $\SBV^{p(\cdot)}(\Omega,\mathcal{M})$}\label{sec:sbv-px}
Here, we collect some basic facts about $\SBV$ functions that are used throughout this paper 
and we define precisely the space of mappings of class $\SBV^{p(\cdot)}$ from an open set $\Omega \subset \R^n$ 
into a 
Riemannian manifold 
$\mathcal{M}$, where $p(\cdot)$ is a bounded, 
variable exponent.
We address the reader to \cite[Chapters~2,3,~4]{AFP} for details about $\SBV$ (and $\BV$) functions as 
well as for all the notions from Geometric Measure Theory (\cite[Chapter~2]{AFP}) that we will use in 
this work.

Let $\Omega \subset \R^n$ be an open set and $k \in \N$. We recall that 
the set $\SBV\left(\Omega,\R^k\right)$ is 
the linear subspace of $L^1\left(\Omega,\R^k\right)$ made up by those functions $u$ whose distributional 
gradient ${\D}u$ is a bounded Radon measure with no Cantor part. Endowed with the norm 
\[
    \norm{u}_{\BV\left(\Omega\right)} := \norm{u}_{L^1\left(\Omega\right)} + \abs{{\D}u}(\Omega),
\]
the space $\SBV\left(\Omega,\R^k\right)$ is a Banach space. 
For any $u \in \SBV\left(\Omega,\R^k\right)$, its distributional gradient ${\D}u$ splits as
\[
    {\D}u = \underbrace{(\nabla u )\mathcal{L}^2}_{:={\D^a}u} + 
    \underbrace{\left[(u^+ - u^-) \otimes \nu_u \right]\mathcal{H}^{n-1} \res J_u}_{:={\D^j}u},
\]
where $\nabla u$ denotes the approximate differential of $u$, $u^\pm$ the traces of 
$u$ on the two sides of the \emph{jump set} $J_u$, and $\nu_u$ is the normal field to $J_u$ 
(see, e.~g., \cite[Definition~3.67]{AFP}). We recall from \cite[Proposition~3.69 and Theorem~3.78]{AFP} 
that $J_u$ is a countably $\mathcal{H}^{n-1}$-rectifiable Borel set in $\Omega$.
%

Let $p \geq 1$ and
\[
    \SBV^p\left(\Omega,\R^k\right) := \left\{ u \in \SBV\left(\Omega,\R^k\right) : \nabla u \in L^p\left(\Omega,\R^{k\times n}\right) \mbox{ and } \mathcal{H}^{n-1}(J_u) < +\infty \right\}.
\]
Endowed with the restriction of the norm, $\SBV^p\left(\Omega,\R^k\right)$ is a Banach subspace 
of $\SBV\left(\Omega,\R^k\right)$
that proved useful in handling free discontinuity problems at least since~\cite{CarrieroLeaci}. 
On the relation between the spaces $\SBV$, $\SBV^p$, $W^{1,1}$, $W^{1,p}$, we recall that: for 
any open set $\Omega \subset \R^n$,
\begin{equation}\label{eq:sbv-and-no-jump-means-W11}
    u \in W^{1,1}\left( \Omega, \R^k\right) \iff u \in \SBV\left( \Omega, \R^k \right) \mbox{ and } 
    \mathcal{H}^{n-1}\left( J_u \right) = 0.
\end{equation}
If $\Omega$ is a bounded open set with Lipschitz boundary, then, by iterated application of Sobolev 
embedding,
\begin{equation}\label{eq:sbv-to-sobolev-p-const}
    u \in W^{1,p}\left( \Omega, \R^k\right) \iff 
    u \in \SBV^p\left( \Omega, \R^k \right) \mbox{ and } 
    \mathcal{H}^{n-1}\left( J_u \right) = 0.
\end{equation}

A straightforward generalisation of the space $\SBV^p\left(\Omega,\R^k\right)$ has been introduced 
in \cite{DeCLV}, by defining the space 
$\SBV^{p(\cdot)}\left(\Omega,\R^k\right)$, where $p(\cdot) : \Omega \to [1,+\infty)$ is any bounded 
variable exponent, as the subspace of those $u \in \SBV\left(\Omega,\R^k\right)$ whose jump set 
has finite $\mathcal{H}^{n-1}$-measure and whose approximate differential is 
$L^{p(\cdot)}$-integrable on $\Omega$.

In this work, we are mostly interested in maps taking values into a compact 
Riemannian manifold without boundary. 
By viewing $\mathcal{M}$ as a submanifold of $\R^k$, for some $k \in \N$ (for instance, by means 
of Nash's isometric embedding theorem), we can define 
\[
	\SBV^{p(\cdot)}\left(\Omega,\mathcal{M}\right) := 
	\left\{ u \in \SBV^{p(\cdot)}\left(\Omega,\R^k\right) : u(x) \in \mathcal{M}\, 
	\mbox{ for a.e. } x \in \Omega \right\}.
\]
This space enjoys closure and compactness properties analogous to those of the classical space 
$\SBV^p\left(\Omega,\R^k\right)$, see Appendix~\ref{app:cpt} for more details.

For later use, we notice that if $\Omega$ is any bounded open set, then 
\begin{equation}\label{eq:sobolev-sbv} 
    u \in W^{1,p(\cdot)}\left(\Omega, \mathcal{M}\right) \iff
    u \in \SBV^{p(\cdot)}\left(\Omega, \mathcal{M}\right) \mbox{ and }  
    \mathcal{H}^{n-1}\left(J_u\right) = 0  
\end{equation}
by compactness of $\mathcal{M}$ (which implies $u \in L^\infty\left(\Omega,\R^k\right)$). 
Moreover, relaxing the assumption of boundedness of the functions, 
but assuming instead that $\Omega$ is a bounded open set with sufficiently nice boundary 
(Lipschitz is enough) and that 
\begin{equation}\label{eq:p+p-star}
    p^- < n \quad \mbox{and} \quad p^+ < \left( p^- \right)^*,
\end{equation}
we have the continuous embedding 
$W^{1,p^-}\left(\Omega, \R^k \right) \hookrightarrow L^{p^+}\left(\Omega, \R^k\right)$, 
see \cite[Corollary~8.3.2]{DHHR}. Consequently, the following lemma holds.
\begin{lemma}\label{lemma:sbv-to-sobolev}
    Assume $\Omega \subset \R^n$ is a bounded open set with Lipschitz boundary and that 
    $p : \Omega \to [1,+\infty)$ is a log-H\"{o}lder continuous, variable exponent 
    satisfying~\eqref{eq:p+p-star}.  
    Then, 
    \[
    u \in W^{1,p(\cdot)}\left(\Omega, \R^k\right) \iff 
    u \in \SBV^{p(\cdot)}\left(\Omega, \R^k\right) \mbox{ and }  
    \mathcal{H}^{n-1}\left(J_u\right) = 0
    \]
\end{lemma}

\begin{proof}
    The direction $\Rightarrow$ is obvious. For the direction $\Leftarrow$, we argue as in the 
    classical case of constant exponents, exploiting assumption~\eqref{eq:p+p-star}: 
    by assumption and~\eqref{eq:sbv-and-no-jump-means-W11}, 
    $u \in W^{1,1}\left(\Omega,\R^k\right)$, which implies, by 
    Sobolev embedding, $u \in L^{1^*}\left( \Omega, \R^k\right)$. 
    If $1^* \geq p^-$, then $u \in W^{1,p^-}\left( \Omega, \R^k\right)$, which implies (in view 
    of~\eqref{eq:p+p-star}) $u \in L^{p^+}\left(\Omega,\R^k\right)$, hence 
    $u \in L^{p(\cdot)}\left(\Omega,\R^k \right)$ by 
    Proposition~\ref{prop:embedding}. Since we know by assumption that $\nabla u \in L^{p(\cdot)}\left(\Omega,\R^{k \times n}\right)$, it follows that 
    $u \in W^{1,p(\cdot)}\left(\Omega,\R^k\right)$. If 
    instead $1^* < p^-$, then $u \in W^{1,1^*}\left(\Omega,\R^k\right)$, and we can repeat the above step. After finitely many iterations, we obtain 
    $1^{**\dots*} \geq p^-$ and we conclude as in the above.
\end{proof}

\begin{remark}
    Assumption~\eqref{eq:p+p-star} is clearly satisfied if $n=2$ and $p^+ < 2$.
\end{remark}


\section{Sobolev approximation of functions in $\SBV^{p(\cdot)}(\Omega,\mathcal{M})$ with small jump set}\label{sec:approx}

In this section, we prove Theorem~\ref{thm:approx-sbv-px} below, 
which clearly entails Theorem~\ref{thm:A}.

\begin{theorem}\label{thm:approx-sbv-px}
    Let $p : \R^2 \to (1,\infty)$ be a bounded, log-H\"older-continuous,
	variable exponent satisfying $p^- > 1$ and $p^+ < 2$. Let $k \in \N$ and $\rho > 0$. 
    Assume $\mathcal{M}$ is a connected, compact Riemannian manifold without boundary, isometrically 
    embedded into $\R^k$.
	There exist universal constants $\xi$,  
	$\eta > 0$ such that 
    for any $s \in (0,1)$ and any 
	$u \in \SBV^{p(\cdot)}\left(B_{\rho}, \mathcal{M}\right)$ satisfying
	\begin{equation}\label{eq:small-jump-set}
		\mathcal{H}^1\left(J_u \cap B_\rho\right) < \eta(1-s) \frac{\rho}{2},
	\end{equation}
	the following holds. There are a countable family 
	$\mathcal{F} = \{B\}$ of closed balls, overlapping at most $\xi$ times, of 
	radius $r_B < (1-s)\rho$ and centre $x_B \in \overline{B_{s \rho}}$, and 
	a function $w \in \SBV^{p(\cdot)}\left(B_\rho, \mathcal{M}\right)$ such that
	\begin{enumerate}[(i)]
 		\item\label{item:approx-iv} $\cup_{\mathcal{F}} B \subset B_{\frac{1+s}{2}\rho }$ and 
		$\frac{1}{\rho}\sum_\mathcal{F} \mathcal{L}^2(B) + \sum_{\mathcal{F}} \mathcal{H}^1(\partial B) \leq \frac{2 \pi \xi}{\eta} \mathcal{H}^1(J_u \cap B_\rho)$. Moreover,
        \begin{equation}\label{eq:estimate-bad-balls}
            \sum_\mathcal{F} \mathcal{L}^2(B) \leq \min\left\{ \frac{2\pi\xi}{\eta} \rho \mathcal{H}^1\left(J_u \cap B_\rho\right), \pi \left( \frac{\xi}{\eta} \mathcal{H}^1\left(J_u \cap B_\rho\right) \right)^2 \right\}
        \end{equation}
		\item\label{item:approx-i} $w = u$ $\mathcal{L}^2$-a.e. on 
		$B_\rho \setminus \cup_{\mathcal{F}} B$.
		\item\label{item:approx-ii} $w \in W^{1,p(\cdot)}\left(B_{s \rho}; \mathcal{M}\right)$ and 
		$\mathcal{H}^1(J_w \setminus J_u) = 0$.
		\item\label{item:approx-iii} There holds
        \begin{equation}\label{eq:control-grad-w-M}
			\int_{B_{\rho}} \abs{\nabla w}^{p(x)} \,{\d}x \lesssim \left(1+\rho^2\right)
			\max\left\{\norm{\nabla u}_{L^{p(\cdot)}(B_{\rho})}^{p^-}, \norm{\nabla u}_{L^{p(\cdot)}(B_{\rho})}^{p^+}\right\}
		\end{equation}
		where the implicit constant in the right hand side is independent of $u$ and $w$ and depends 
        only on the quantities listed in Remark~\ref{rk:projection} and on the log-H\"{o}lder 
        constant of $p(\cdot)$. Moreover, if $p$ is constant, 
        then~\eqref{eq:control-grad-w-M} holds without the factor $(1+\rho^2)$ at right hand side.
	\end{enumerate}
\end{theorem}

The proof of Theorem~\ref{thm:approx-sbv-px} proceeds essentially in three steps. After having 
isometrically embedded $\mathcal{M}$ into some Euclidean space $\R^k$, for some $k \in \N$ (Step~1), 
we prove an analogous approximation result (Proposition~\ref{thm:approx-sbv-px-Rk}) for 
\emph{unconstrained} maps, i.e., for maps with values into $\R^k$ (Step~2). Then, we obtain 
Theorem~\ref{thm:approx-sbv-px} from Proposition~\ref{thm:approx-sbv-px-Rk} by means of a suitable 
retraction onto $\mathcal{M}$ (Step~3).

As mentioned in Section~\ref{sec:notation}, Step~1 can always accomplished by exploiting 
Nash's isometric embedding theorem. (However, since $\mathcal{M}$ is compact the choice of the 
embedding is essentially irrelevant, as it can only imply, at worst, enlarging the constant in~\eqref{eq:control-grad-w-M}.)
Step~2 is heavily based on arguments in \cite{CFI-ARMA} and worked out in Section~\ref{sec:approx-Rk} 
below. We complete the proof of Theorem~\ref{thm:approx-sbv-px} in Section~\ref{sec:approx-M}.

\begin{remark}\label{rk:bad-scaling}
    The factor $\left(1+\rho^2\right)$ at right hand side in~\eqref{eq:control-grad-w-M} is due to the 
    use of the embedding theorem for variable-exponent Lebesgue spaces 
    (Proposition~\ref{prop:embedding}) in the 
    proof of Proposition~\ref{thm:local-approx} below. For constant exponents, this is not needed and 
    one can instead  
    stick with the clever, optimal argument in \cite{CFI-ARMA} and obtain~\eqref{eq:control-grad-w-M} 
    without the extra-factor 
    $\left(1+\rho^2\right)$, c.f.~\eqref{eq:nabla-phi-u}, \eqref{eq:widetilde-phi-q-norm-grad-ineq} 
    below. 
    However, the dependence on $\rho$ at right hand side of~\eqref{eq:control-grad-w-M} is harmless for 
    our purposes in this paper, as we apply Theorem~\ref{thm:approx-sbv-px} only to maps defined on 
    a fixed, common ball (c.f. the blow-up analysis in the proof of Theorem~\ref{thm:decay-lemma}).
\end{remark}

\subsection{The unconstrained case}\label{sec:approx-Rk}
For the unconstrained case, 
we follow very closely the argument of \cite[Section~2 and Section~3]{CFI-ARMA}, developed in the 
setting of constant exponents and functions of class $\SBD^p$. As mentioned at the beginning of the 
section, the main result here is Proposition~\ref{thm:approx-sbv-px-Rk}. 
The key tool in the proof of Proposition~\ref{thm:approx-sbv-px-Rk} is 
Proposition~\ref{thm:local-approx} below. The statement of Proposition~\ref{thm:local-approx} is similar to that of the corresponding result in \cite{CFI-ARMA}, i.e., \cite[Theorem~2.1]{CFI-ARMA}, 
up to the complications arising because of the variable exponent. In particular, we shall need the 
following assumption~{(H)}:

\begin{equation}\leqnomode
    \tag{H} u \mbox{ is bounded} 
    \quad \mbox{{\em or}} 
    \quad p(\cdot) \mbox{ satisfies } p^- < n \,\,\mbox{and}\,\, p^+ \leq \left(p^-\right)^*. 
\end{equation}
This assumption is always satisfied in this work, because we deal with maps with values into a compact 
Riemannian manifold $\mathcal{M}$ and moreover, as mentioned in the Introduction and explained in more 
details in Section~\ref{sec:approx-M} below, to approximate them with maps with values into 
$\mathcal{M}$ (which is only assumed to be connected, in general) we need to assume $p^+ < 2$. 

Assumption~{(H)} ensures that every function of class $W^{1,1}$ with 
$L^{p(\cdot)}$-integrable gradient is also of class $W^{1,p(\cdot)}$ and, in turn, that  
functions belonging to $\SBV^{p(\cdot)}$ with no jump are also of class $W^{1,p(\cdot)}$ (see the 
discussion in Section~\ref{sec:sbv-px} and the end of Step~\ref{step:approximant-in-sbv-px} below for 
more details on this point). Finally, we notice that also in the constant-exponent case the regime 
$p<n$ is the most interesting one 
(and that the condition $p < p^*$ is trivially satisfied in such case). 

\begin{prop}\label{thm:local-approx}
	Let $p : \R^2 \to (1,+\infty)$ be a bounded, log-H\"older-continuous,
	variable exponent satisfying $p^- > 1$ and let $k \in \N$.  
    There exist a universal constant $\eta$ and a constant $\widetilde{c}_k$, depending only on $k$, 
	such that for any $r > 0$, if $J \in \mathscr{B}(B_{2 r})$ is any Borel set in $B_{2r}$ 
    that satisfies
	\begin{equation}\label{eq:H1-J}
		\mathcal{H}^1(J) < \eta(2r),
	\end{equation}
	then there exists $R \in (r,2r)$ for which the following holds.
    Let
	$u \in \SBV^{p(\cdot)}\left(B_{2r},\R^k\right)$ be any function such that 
	\[
	   \mathcal{H}^1\left(J_u \cap B_{2r} \setminus J\right) = 0,
	\] 
    and suppose, in addition, that assumption ${({\rm H})}$ holds.
	Then, there exists 
	${\phi}(u) \in \SBV^{p(\cdot)}\left(B_{2r},\R^k\right) \cap W^{1,p(\cdot)}\left(B_R, \R^k\right)$ 
	such that the following properties are satisfied.
	\begin{enumerate}[(i)]
		\item $\mathcal{H}^1(J_u \cap \partial B_R) = 0$.
		\item There holds 
		\begin{equation}\label{eq:q-norm-grad-phi}
		\int_{B_R} \abs{\nabla \phi(u)}^{q} \,{\d}x \leq 
	\widetilde{c}_k \int_{B_R} \abs{\nabla u}^{p^-} \,{\d}x,
		\end{equation}
		for every $q \in [1, p^-]$. Moreover,
		\begin{equation}\label{eq:p(x)-norm-grad-phi}
			\int_{B_R} \abs{\nabla \phi(u)}^{p(x)} \,{\d}x \lesssim 
			\left(1+ R^2\right)  \max\left\{ \norm{\nabla u}_{L^{p(\cdot)}(B_R)}^{p^+}, 
			\norm{\nabla u}_{L^{p(\cdot)}(B_R)}^{p^-} \right\},
		\end{equation}
		where the implicit constant at right hand side depends only on the log-H\"{o}lder constant 
        of $p$.
		\item There holds
		\begin{equation}\label{eq:dist-u-phi(u)}
			\norm{u - {\phi}(u)}_{L^1\left(B_R\right)} \lesssim R \abs{{\D}u}(B_R),
		\end{equation} 
		where the implicit constant at right hand side depends only on $k$.
		\item ${\phi}(u) = u$ a.e. on $B_{2r}\setminus \overline{B_R}$ and $\mathcal{H}^1\left(J_{\phi(u)} \cap \partial B_R\right)=0$.
		\item $\norm{\phi(u)}_{L^\infty(B_{2r})} \leq \norm{u}_{L^\infty(B_{2r})}$. 
	\end{enumerate} 
\end{prop}

\begin{proof}
	The strategy of the proof is borrowed from \cite[Theorem~2.1]{CFI-ARMA}, 
	which deals with maps of class $\SBD^p(B_{2r})$, for $1 \leq p < +\infty$ 
	a constant exponent. In rough words, the idea is to choose a ``good radius'' 
    $R \in (r,2r)$ for a ball $B_R$, concentric with $B_{2r}$, and to construct a grid  
	$B_R$ by triangles, refining towards the boundary of $B_R$, whose 
	vertices are chosen in such a way to satisfy essentially three properties: 
	{(i)} they are Lebesgue points of the function $u$ to approximate; 
	{(ii)} the one-dimensional restriction of $u$ to each edge of the 
	grid is of class $W^{1,1}$, and 
	{(iii)} the value of $u$ at the vertices is close to the local average. 
    (See, more precisely, properties (${\rm P}_1$)--(${\rm P}_5)$ below.)
    A crucial point is that, in 2D, such a grid can be constructed with {\em universally} 
    bounded geometry under the mere assumption~\eqref{eq:H1-J}.
    
	Once the grid has been built, one can then define a piecewise affine map ${\phi}(u)$ in $B_R$ by 
	prescribing that it agrees with $u$ at the vertices of each triangle of 
	the grid. Outside $B_R$, one sets $\phi(u) := u$, and it turns out that 
    the map $\phi(u)$ so defined satisfies all the properties in the statement. 
    
	Here, proceeding in a sequence of steps, we modify the approach of~\cite{CFI-ARMA} to 
	handle the complications arising because of the variable exponent,  at the same time benefiting 
    from some simplifications occurring because we deal with maps of class $\SBV$. 
	The main idea is to exploit the constancy of   
	$\nabla \phi(u)$ in each triangle of the grid to 
	make a sort of local ``$L^\infty-L^{p^-}$-interpolation'', to estimate locally the 
	$p(\cdot)$-modular of $\nabla \phi(u)$, c.f. Step~\ref{step:interpolation}. 
    Then, we obtain global estimates summing over all triangles of the grid.
	The price to pay is the slight loss of exponent at right hand side 
	in~\eqref{eq:p(x)-norm-grad-phi}. 
	Nevertheless, the estimate~\eqref{eq:p(x)-norm-grad-phi} will be enough to our purposes, 
	i.e., for the blow-up analysis in Theorem~\ref{thm:decay-lemma}.

	\setcounter{step}{0}
	
	\begin{step}[Construction of the grid]
	Following the argument in \cite[Theorem~2.1]{CFI-ARMA}, we can pick $R \in (r,2r)$ in such a way 
    that the following holds:
    \begin{equation}\label{eq:H1-null-J-cap-B_R}
        \mathcal{H}^1(J \cap \partial B_R) = 0.
    \end{equation}
    and
    \begin{equation}\label{eq:(2.3)-CFI-ARMA}
        \mathcal{H}^1\left( J \cap \left( B_{R} \setminus B_{R\left(1-2^{-h}\right)} \right) \right) < 10 \eta \delta_k, \quad \mbox{for every } h \in \N.
    \end{equation}
    (C.f.~\cite[(2.2) and (2.3)]{CFI-ARMA}.) 
    Once $R$ has been chosen as in the above, 
	we construct a triangular grid $\mathcal{T}$ for $B_R$, refining 
	towards the boundary of $B_R$, in a completely analogous way to as in first 
	part of the proof of \cite[Theorem~2.1]{CFI-ARMA}. We recall the main 
	properties of $\mathcal{T}$ and the main steps of its construction, 
	referring the reader to \cite[pp.~1340--1344]{CFI-ARMA} for full details.
	
	The construction of $\mathcal{T}$ starts in a completely geometrical way, 
    by considering the circles of radii 
	$R_h := R - \delta_h$, with $h \in \N$ and $\delta_h := R 2^{-h}$, and, for each 
    such circle, the $2^{h}$ distinguished points $x'_{h,j}$ given by 
	\[
		x'_{h,j} := R_h\left(\cos \frac{2\pi j}{2^{h}}, \sin \frac{2\pi j}{2^{h}}\right), \quad \mbox{for } j \in \left\{1,\dots,2^{h}\right\}.
	\]
	As in \cite[p.~1341]{CFI-ARMA}, we observe
	that any pair $x'$, $y'$ of neighbouring points in 
	$\mathcal{V}' := \left\{ x'_{h,j} \right\}_{h,j}$ 
	(the notion of ``neighbouring points'' in $\mathcal{V}'$ 
    coincides with the intuitive one and it 
    is formalised in \cite[pp.~1340-1341]{CFI-ARMA}) satisfies 
	\begin{equation}\label{eq:bound-dist-xy}
 		c_1 \delta_h \leq \abs{x' - y'} \leq c_2 \delta_h,
 	\end{equation}
 	where $c_1$, $c_2 > 0$ are universal constants. 
 	Connecting all such vertices, we obtain a grid 
 	$\mathcal{T}'$ with universally bounded geometry. 
  
    Now, we fix any function $u \in \SBV^{p(\cdot)}\left(B_{2r}, \R^k\right)$ such 
    that assumption (${\rm H}$) is satisfied (so, if $p$ does not satisfy \eqref{eq:p+p-star}, 
    we assume in addition that $u$ is bounded) 
    and such that $\mathcal{H}^1\left(J \setminus J_u\right) = 0$. 
    We are going to adapt the grid $\mathcal{T}'$ to the function $u$, 
    obtaining a good grid $\mathcal{T}$ on which 
 	constructing a piecewise linear approximation of $u$.
 	
 	To the claimed purpose, we follow the argument in \cite{CFI-ARMA} 
    and we start by taking any triangle $T'$ of the grid 
	$\mathcal{T}'$. We denote $x'$, $y'$, $z'$ be the 
	vertices of $T'$, and by $s_{x',y'}$, 
	$s_{x',z'}$, $s_{y',z'}$ the segments connecting them. We define 
	\begin{equation}\label{eq:def-Q}
		Q_{x',y'} := \left\{ \xi \in B_R : \dist\left(\xi, s_{x',y'}\right) < \abs{x' - y'}/(8c_2) \right\},
	\end{equation}
	and, in the same way, the sets 
	$Q_{x',z'}$, $Q_{y',z'}$. Next, we define 
 	\[
 		\alpha := c_1/(8c_2)
 	\]
 	and we consider, for each pair of indexes $(h,j)$, 
 	the ball $B( x'_{h,j}, \alpha \delta_h)$, of radius 
 	$\alpha \delta_h$ and center $x'_{h,j}$. Let $x'$, $y'$ any two 
    neighbouring points in $\mathcal{V}'$. 
    The sophisticated iterative procedure 
	in \cite[pp.~1341-1344]{CFI-ARMA} shows that one can determine two {\bf universal} constants 
    $\eta > 0$ (bounding from above the ratio $\mathcal{H}^1\left(J\right) / (2r)$) and 
    $\widetilde{c} > 0$ so 
    that there is a ``good'' set $\mathcal{G} \subset B\left( x', \alpha\delta_h \right) \times B\left(y', \alpha\delta_h \right)$ with positive measure and whose 
    points $(x,y)$ satisfy the following properties:
    \begin{itemize}
        \item[(${\rm P}_1$)] $u^\nu_z \in \SBV\left(s_{x,y}, \R^k\right)$.
        \item[(${\rm P}_2$)] $\mathcal{H}^0\left(J_{u^\nu_z}\right) = 0$, so that $u^\nu_z \in W^{1,1}\left(s_{x,y}, \R^k\right)$.
        \item[(${\rm P}_3$)] There holds 
        \begin{equation}\label{eq:ineq-CFI-ARMA-P3-type}
        \int_{s_{x,y}} \abs{(u^\nu_z)^\prime} \,{\d}t \leq  
        \frac{\widetilde{c}}{\delta_h} \int_{O_{x',y'}} \abs{\nabla u} \,{\d}x^\prime,
        \end{equation}
        where $O_{x',y'}$ denotes the convex envelope of the union 
        $B(x',\alpha\delta_h) \cup B(y',\alpha\delta_h)$.
        \item[(${\rm P}_4$)] For $\xi \in \{x,y\}$, there holds  
        \[
            \abs{u(\xi) - u_{x',y'}} \leq \frac{\widetilde{c}}{\delta_h} \abs{{\D}u}\left( Q_{x',y'} \right),
        \]
        where $u_{x',y'}$ denotes the average of $u$ in $Q_{x',y'}$.
        \item[(${\rm P}_5$)] $x$ and $y$ are Lebesgue point of $u$,
    \end{itemize}
    In the above, $u^\nu_z$ denotes the 1D-section of $u$ along the direction $\nu$ determined 
    by $x$, $y$, given by
    \[
        u^\nu_z := u(z+t \nu), \qquad z := \left( {\rm Id} - \nu \otimes \nu \right) x, \qquad 
        \nu := \frac{x-y}{\abs{x-y}}. 
    \]
    These properties are exactly the analogues of properties (${\rm P}_1$)--(${\rm P}_5$) 
    in \cite[Theorem~2.1]{CFI-ARMA} and they are proven with exactly the same reasoning as 
    in~\cite{CFI-ARMA}, just 
    replacing the symmetric gradient with the full gradient of $u$ and the rigid motions with averages.

    Arguing word-by-word as in \cite[p.~1344]{CFI-ARMA}, we can extract out of $\mathcal{G}$ a set of 
    points $\mathcal{V}$ which are good vertices for a new triangulation $\mathcal{T}$, which is 
    adapted to $u$ in the following sense:
	connecting all neighbouring points in $\mathcal{V}$ by edges (agreeing, 
	as in \cite{CFI-ARMA}, that two points in $\mathcal{V}$ are neighbours 
	if and only if they come from points that were neighbours in the previous 
	sense), we obtain a grid $\mathcal{T}$ by triangles $T$ whose vertices satisfy 
    (${\rm P}_1$)-(${\rm P}_5)$. 
    (In particular, the edges of $T$ do not intersect the jump set of $u$.) 
 	In addition, any pair of neighbouring vertices in $\mathcal{V}$
 	satisfy~\eqref{eq:bound-dist-xy}, and hence the angles of the triangles $T$ 
 	are uniformly bounded away from $0$ and $\pi$. 
 
 Moreover, for later purposes, we observe that the triangles $T'$ of $\mathcal{T}'$ and the triangles $T$ of $\mathcal{T}$ are 
 in a one-to-one correspondence. For any given $T \in \mathcal{T}$, 
 we denote $C_T$ the following convex envelope: 
 \[
 	C_T := {\rm conv}\left( \cup B\left(x', \alpha \delta_h\right)\right),
 \]
 where the union runs over the vertices $x'$ of the triangle 
 $T'$ in $\mathcal{T}'$ which correspond to $T$. Clearly, $O_{x',y'} \subset C_T$ for any 
 two vertices $x'$, $y'$ of $T'$ and $T \subset C_T$. 
 Moreover, following~\cite{CFI-ARMA}, we observe that there is a universal constant 
 $\kappa$ such that any $x \in B_R$ belongs to at most $\kappa$ 
 of the convex envelopes $C_T$. 
 In addition, we remark that, by construction, there exists two universal constants $\lambda_1$, $\lambda_2 > 0$ such that   
 \begin{equation}\label{eq:ratios}
 \begin{split}
     \lambda_1 &\leq \min\left\{\frac{\mathcal{L}^2\left(C_T\right)}{\mathcal{L}^2\left(Q_{x',y'}\right)}, \, 
     \frac{\mathcal{L}^2\left(T\right)}{\mathcal{L}^2\left(Q_{x',y'} \cap T \right)}, \, 
     \frac{\mathcal{L}^2\left(B\left(x',\alpha\delta_h\right)\right)}{\mathcal{L}^2\left(T\right)}, \,
     \frac{\mathcal{L}^2\left(O_{x',y'}\right)}{\mathcal{L}^2\left(C_T\right)}
     \right\}\\
     & \leq \max\left\{ \frac{\mathcal{L}^2\left(C_T\right)}{\mathcal{L}^2\left(Q_{x',y'}\right)}, \, 
     \frac{\mathcal{L}^2\left(T\right)}{\mathcal{L}^2\left(Q_{x',y'} \cap T \right)}, \, 
     \frac{\mathcal{L}^2\left(B\left(x',\alpha\delta_h\right)\right)}{\mathcal{L}^2\left(T\right)},\, 
     \frac{\mathcal{L}^2\left(O_{x',y'}\right)}{\mathcal{L}^2\left(C_T\right)}\right\}\leq \lambda_2 .
  \end{split}
 \end{equation}
 Let us set $\lambda := \max\{ \lambda_1,\lambda_2 \}$.
 \end{step}
 	
 \begin{step}[Definition and main properties of the approximating map]
 We define, as in \cite{CFI-ARMA}, a function ${\phi}(u)$ 
 by letting ${\phi}(u) = u$ in $B_{2r} \setminus \overline{B_R}$ while we take 
 ${\phi}(u)$ to be the piecewise affine function determined by $\mathcal{V}$ on $B_R$. More 
 precisely, in each triangle $T \in \mathcal{T}$, we define 
 ${\phi}(u)$ as the (uniquely determined) affine function obtained by 
 setting ${\phi}(u)(x) = u(x)$, for each vertex $x$ of $T$. In particular, on the edges of $T$, 
 ${\phi}(u)$ is the linear interpolation between the values 
 of $u$ at the vertices and $\nabla {\phi}(u)$ is a 
 constant $(k \times 2)$-matrix in each $T \in \mathcal{T}$. 
 As a consequence, given any vertices  $x$, $y$ of $T$, we have
 \begin{equation}\label{eq:CFI-trick}
    \phi(u)(x-y) = \phi(u)(x) - \phi(u)(y) = \left(\nabla \phi(u)\right)(x-y) = u(x) - u(y).
 \end{equation}
With the aid
    of the fact that $u^\nu_z \in W^{1,1}\left(s_{x,y}\right)$ and of~\eqref{eq:ineq-CFI-ARMA-P3-type}, 
    we can compute the elements of the matrix $\nabla \phi(u)$ and its $L^\infty$-norm in complete 
    analogy with \cite{CFI-ARMA}, obtaining the following counterparts of
	\cite[(2.12) and (2.13)]{CFI-ARMA}. For each $T \in \mathcal{T}$, every 
    pair $x$, $y$ of vertices of $T$, and every $j \in \{1,\dots,k\}$, by 
    the fundamental theorem of calculus on $s_{x,y}$ and~\eqref{eq:CFI-trick} we have  
    \[
        \nabla {\phi}(u) \nu \cdot {\bf e}_j  = 
        \frac{(\phi(u)(x) - \phi(u)(y))\cdot {\bf e}_j}{\abs{x-y}} = 
        \fint_{s_{x,y}} \left(u^\nu_z \cdot {\bf e}_j\right)^\prime \,{\d}t,
    \]
    hence, 
    by~\eqref{eq:ineq-CFI-ARMA-P3-type} and~\eqref{eq:ratios},
    \begin{equation}
        \abs{\nabla {\phi}(u) \nu \cdot {\bf e}_j } \leq 
        \fint_{s_{x,y}} \abs{\left(u^\nu_z\right)^\prime} \,{\d}t
        \stackrel{\eqref{eq:ineq-CFI-ARMA-P3-type},\eqref{eq:ratios}}{\leq} 
        \widetilde{c}\lambda\fint_{C_T} \abs{\nabla u}\,{\d}x. 
        \label{eq:ineq-type-CFI-ARMA-2.12}
    \end{equation}
    Thus, letting $x$, $y$ vary in the set of vertices of $T$ and $j$ run in $\{1,\dots,k\}$,
	\begin{equation}\label{eq:nabla-phi-u}
 		\norm{\nabla {\phi}(u)}_{L^\infty(T)} \leq 
 		\widetilde{c}\lambda \fint_{C_T} \abs{\nabla u}\,{\d}x,
 	\end{equation}
 	for each $T \in \mathcal{T}$. 
    From~\eqref{eq:nabla-phi-u} and Jensen's inequality, we obtain 
 	\begin{equation}\label{eq:widetilde-phi-q-norm-grad-ineq}
 		\int_T \abs{\nabla {\phi}(u)}^q \leq \widetilde{c} \lambda
 		\int_{C_T} \abs{\nabla u}^q\,{\d}x 
 	\end{equation}
 	for any $q \in [1,p^{-}]$ and any $T \in \mathcal{T}$. 
    Recalling that the sets $C_T$ overlap at most $\kappa$-times, and that $\kappa$ is a universal 
    constant, by taking the sum of~\eqref{eq:widetilde-phi-q-norm-grad-ineq} over $T \in \mathcal{T}$ 
    it follows from~\eqref{eq:sbv-to-sobolev-p-const} that $\phi(u)$ belongs to 
    $W^{1,q}\left(B_R, \R^k\right) \cap \SBV^q\left(B_{2r}, \R^k\right)$, for any $q \in [1,p^-]$.
    In addition, by construction of the grid, the $\BV$ trace of $u$ on $\partial B_R$
    and the Sobolev trace of $\phi(u)$ on $\partial B_R$ agree $\mathcal{H}^1$-a.e., hence no 
    additional jump is created by the above process (see \cite[p.~1346]{CFI-ARMA} for details). Consequently, 
    $\phi(u)$ has less jump than $u$ in $B_{2r}$.
 \end{step}

\begin{step}[{Proof of {(i)}, {(iii)}, {(iv)}, {(v)}}]
	\label{step:first-properties}
	Item~{(v)} is obvious. Item {(i)} follows immediately, because of~\eqref{eq:H1-null-J-cap-B_R} 
    and since $\mathcal{H}^1\left(J\setminus J_u\right) = 0$.  
    The first assertion in {(iv)} follows by definition of $\phi(u)$ and the argument 
	in \cite[p.~1346]{CFI-ARMA}. The second assertion in {(iv)} follows because, by construction,  
    $J_{\phi(u)} = J_u$ on $B_{2r} \setminus \overline{B_R}$ and, as observed in the previous step, 
    $u$ and $\phi(u)$ agree $\mathcal{H}^1$-a.e. on $\partial B_R$.
	
	To prove {(iii)}, we argue similarly to as in \cite{CFI-ARMA}. 
	Let $T'$ be any triangle of the original, geometrical grid 
	$\mathcal{T}'$ and let $x'$, $y'$, $z'$ be its 
	vertices, connected by the segments $s_{x',y'}$, 
	$s_{x',z'}$, $s_{y',z'}$.
    Let $Q_{x',y'}$ be defined as in~\eqref{eq:def-Q}. By
	Proposition~\ref{prop:poincare-BV-convex} and~\eqref{eq:bound-dist-xy}, we have
	\begin{equation}\label{eq:poincare-Qx'y'}
		\norm{u - u_{x',y'} }_{L^{1}(Q_{x',y'})} \leq 
		c_k \delta_h \abs{{\D}u}(Q_{x',y'}),
	\end{equation}
	where $u_{x',y'}$ denotes the average of $u$ in 
	$Q_{x',y'}$ and $c_k := c(k) c_2$, 
	where $c(k) \equiv C(2,k)$ is the constant, depending only on $k$, 
	provided by Proposition~\ref{prop:poincare-BV-convex} (for $n=2$) and 
	$c_2$ is the constant in~\eqref{eq:bound-dist-xy}. Moreover, by (${\rm P}_4$), 
    it follows that 
	\begin{equation}\label{eq:P4-type-inequality}
		\abs{u(\xi) - u_{\zeta',\xi'}} \lesssim \delta_h^{-1} \abs{{\D}u}(Q_{\zeta',\xi'})
		\quad \mbox{for any } \zeta', \xi' \in \{x',y',z'\}  \,\,\mbox{ with }\zeta'\neq \xi',
	\end{equation}
	where the implicit constant at right hand side is universal.
	 
	Next, for every $T \in \mathcal{T}$, let $T' \in \mathcal{T}'$ be the corresponding 
	triangle in the original grid, with vertices $x'$, $y'$, $z'$. By triangle inequality we have
	\begin{equation}\label{eq:dist-u-phi(u)-compu1}
		\int_T \abs{u - \phi(u)} \,{\d}x \leq 
		\int_T \abs{u - u_{x',y'}}\,{\d}x + 
		\int_T \abs{u_{x',y'} - \phi(u)}\,{\d}x.	
	\end{equation}
	About the first term above, we notice that, again  
	by triangle inequality, there holds
	\begin{equation}\label{eq:dist-u-phi(u)-compu2-bis}
		\int_T \abs{u - u_{x',y'}}\,{\d}x \leq 
		\int_T \abs{u - u_T} \,{\d}x + 
		\int_T \abs{u_T - u_{x',y'}} \,{\d}x,
	\end{equation}
    where $u_T$ denotes the average of $u$ in $T$.
	For the first term on the right in~\eqref{eq:dist-u-phi(u)-compu2-bis}, 
	by Proposition~\ref{prop:poincare-BV-convex} and \eqref{eq:bound-dist-xy}, 
	we get
	\begin{equation}\label{eq:dist-u-phi(u)-compu2}
		\norm{u - u_T}_{L^1(T)} \leq c_k \delta_h \abs{{\D}u}(T),
	\end{equation}
	where, again, $c_k = c(k) c_2$. 
	To estimate the second term in~\eqref{eq:dist-u-phi(u)-compu2-bis},  
	we first bound the constant $\abs{u_T - u_{x',y'}}$ 
	using~\eqref{eq:poincare-Qx'y'} and \eqref{eq:dist-u-phi(u)-compu2} 
	as follows:
	\[
	\begin{split}
		\mathcal{L}^2(Q_{x',y'} \cap T) \abs{u_T - u_{x',y'}} &= 
		\int_{Q_{x',y'} \cap T} \abs{u_T - u_{x',y'}} \,{\d}x\\
		&\leq \int_{Q_{x',y'}} \abs{u_{x',y'}-u} \,{\d}x + 
		\int_{T} \abs{u - u_T} \,{\d}x \\
		&\leq c_k \delta_h \abs{{\D}u}(Q_{x',y'}) 
		+ c_k \delta_h \abs{{\D}u}(T).
	\end{split}
	\]
	Thus,
	\[
		\int_T \abs{u_{x',y'}-u_T} \,{\d}x \leq 
		\frac{c_k \mathcal{L}^2(T)}{\mathcal{L}^2(Q_{x',y'} \cap T)} 
		\left( \delta_h \abs{{\D}u}(Q_{x',y'}) + \delta_h \abs{{\D}u}(T) \right).
	\]
    Thus, from~\eqref{eq:ratios}, 
	\begin{equation}\label{eq:dist-u-phi(u)-compu3}
    	\int_T \abs{u_{x',y'}-u_T} \,{\d}x \leq 
		c_k \lambda \delta_h \left(  \abs{{\D}u}(Q_{x',y'}) + \abs{{\D}u}(T) \right) 
		\leq \widetilde{c}_k \delta_h \abs{{\D}u}(C_T),
	\end{equation}	 
	where $\widetilde{c}_k := c_k \lambda$ is a constant depending only on $k$. 
    Combining~\eqref{eq:dist-u-phi(u)-compu2-bis}, \eqref{eq:dist-u-phi(u)-compu2}, 
    and~\eqref{eq:dist-u-phi(u)-compu3}, we obtain
    \begin{equation}\label{eq:dist-u-phi(u)-compu4}
        \int_T \abs{u-u_{x',y'}}\,{\d}x \leq \widetilde{c}_k \delta_h \abs{{\D}u}\left(C_T\right),
    \end{equation}
    where $\widetilde{c}_k$ depends only on $k$.
	
	Concerning the second term in~\eqref{eq:dist-u-phi(u)-compu1}, 
	we proceed in two steps. First, we observe that the function 
	$u_{x',y'} - {\phi}(u)$ is an affine function, and hence its 
	modulus achieve its maximum at one of the vertices of $T$. Denoting $\xi$ 
	such a vertex, we have
	\[
		\int_T \abs{u_{x',y'} - {\phi}(u)}\,{\d}x \leq 
		\mathcal{L}^2(T) \norm{u_{x',y'} - {\phi}(u)}_{L^\infty(T)} 
		= \mathcal{L}^2(T) \abs{u_{x',y'} - {\phi}(u)(\xi)} 
	\]
	Recalling that $u(\xi) = \phi(u)(\xi)$, if $\xi = x$ or $\xi =y$ then, 
	by~\eqref{eq:P4-type-inequality},
	\begin{equation}\label{eq:dist-phi(u)-compu-2.3}
		\mathcal{L}^2(T) \abs{u_{x',y'} - u(\xi)} \leq 
		\widetilde{c} \delta_h \abs{{\D}u}(Q_{x',y'}) 
        \leq \widetilde{c} \delta_h \abs{{\D}u}(C_T),
	\end{equation}
	If instead $\xi = z$, then we write
	\[
		\abs{u_{x',y'} - u(z)} \leq \abs{u_{x',y'} - u_{x',z'}} + \abs{u_{x',z'}-u(z)}
	\]
	The second term above is again estimated by~\eqref{eq:P4-type-inequality}, which yields
	\begin{equation}\label{eq:dist-phi(u)-compu-2.4}
		\abs{u_{x',z'}-u(z)} \lesssim \delta_h^{-1} \abs{{\D}u}(Q_{x',z'}) 
		\lesssim c \delta_h^{-1} \abs{{\D}u}(C_T),
	\end{equation}
	where the constant $c$ is universal. Next, we estimate the 
	constant $\abs{u_{x',y'} - u_{x',z'}}$ as follows: 
	\[
	\begin{split}
		\mathcal{L}^2(B(x',\alpha\delta_h))\abs{u_{x',y'} - u_{x',z'}} &= 
		\int_{B(x',\alpha\delta_h)} \abs{u_{x',y'} - u_{x',z'}} \,{\d}x \\
		&\leq 
		\int_{B(x',\alpha\delta_h)} \abs{u_{x',y'} - u} \,{\d}x + 
		\int_{B(x',\alpha\delta_h)} \abs{u - u_{x',z'}} \,{\d}x \\
		&\leq 
		\int_{Q_{x',y'}} \abs{u_{x',y'} - u} \,{\d}x + 
		\int_{Q_{x',z'}} \abs{u - u_{x',z'}} \,{\d}x \\
		&\lesssim \delta_h \left( \abs{{\D}u}(Q_{x',y'}) + \abs{{\D}u}(Q_{x',z'})\right)
		\lesssim \delta_h \abs{{\D}u}(C_T),
	\end{split}
	\]
    where in the last line we used again Proposition~\ref{prop:poincare-BV-convex} and 
    the implicit constant depends only on $k$. Therefore, using again~\eqref{eq:ratios}
	\begin{equation}\label{eq:dist-phi(u)-compu-2.5}
		\int_T \abs{u_{x',y'} - u_{x',z'}}  \,{\d}x = \mathcal{L}^2(T)\abs{u_{x',y'} - u_{x',z'}} 
		\lesssim \delta_h \abs{{\D}u}(C_T),
	\end{equation}
	up to constant depending only on $k$. 
	Combining~\eqref{eq:dist-phi(u)-compu-2.3}, 
	\eqref{eq:dist-phi(u)-compu-2.4}, and~\eqref{eq:dist-phi(u)-compu-2.5}, 
	we obtain
	\begin{equation}\label{eq:dist-phi(u)-compu3}
		\int_T \abs{u_{x',y'} - {\phi}(u)}\,{\d}x \lesssim 
		\delta_h {\abs{{\D}u}(C_T)},
	\end{equation}
	where the implicit constant at right hand side depends only on $k$. 
    Combining~\eqref{eq:dist-u-phi(u)-compu4} and~\eqref{eq:dist-phi(u)-compu3}, 
    by~\eqref{eq:dist-u-phi(u)-compu1} we get
    \[
        \int_T \abs{u - \phi(u)}\,{\d}x \lesssim \delta_h \abs{{\D}u}\left(C_T\right),
    \]
    where the implicit constant at right hand side depends only on $k$. 
    Therefore, taking the sum of the above inequalities as $T$ ranges over $\mathcal{T}$ 
    (recalling once again that the sets $C_T$ have finite overlap), we obtain~\eqref{eq:dist-u-phi(u)}.
\end{step}
		
\begin{step}[$L^\infty$-$L^{p^-}$ interpolation]\label{step:interpolation}
	We notice that if $\Omega \subset \R^n$ is any measurable set, 
	$f \in (L^{p^-} \cap L^\infty)\left(\Omega,\R^k\right)$ is arbitrary, and 
	$p(\cdot)$ is any bounded, variable exponent, with 
	$1 \leq p^- \leq p(x) \leq p^+ < \infty$, then for a.e. $x \in \Omega$ there holds	
 \begin{equation}\label{eq:interpolation-1}
	\abs{f(x)}^{p(x)} \leq \norm{f}_{L^\infty(\Omega)}^{p(x)-p^-} \abs{f(x)}^{p^-}.
 \end{equation}
 Hence,
 \begin{equation}\label{eq:interpolation-2}
 	\int_\Omega \abs{f(x)}^{p(x)}\,{\d}x \leq 
 	\int_\Omega \norm{f}_{L^\infty(\Omega)}^{p(x)-p^-} \abs{f(x)}^{p-} \,{\d}x
 	\leq \max\left\{ 1,  \norm{f}^{p^+-p^{-}}_{L^\infty(\Omega)}\right\} \int_\Omega \abs{f(x)}^{p^-}\,{\d}x.
 \end{equation}
\end{step}
	
\begin{step}[Conclusion]\label{step:approximant-in-sbv-px}
 For each $T \in \mathcal{T}$, we apply~\eqref{eq:interpolation-2} with 
 $\Omega = T$ and $f = \nabla {\phi}(u)$ (which is 
 constant in $T$, being ${\phi}(u)$ affine in $T$). 
 By~\eqref{eq:nabla-phi-u}, \eqref{eq:interpolation-2} and 
 H\"older's inequality, we obtain 
 \[
 	\int_T \abs{\nabla {\phi}(u)}^{p(x)}\,{\d}x \leq 
 	\max\left\{ 1, \norm{\nabla {\phi}(u)}^{p^+-p^-}_{L^\infty(T)} \right\}
 	\int_{C_T} \abs{\nabla u}^{p^-}	\,{\d}x,
 \]  
 where $C_T$ is defined as in Step~\ref{step:first-properties}. 
 On the other hand, using again~\eqref{eq:nabla-phi-u}, 
 \begin{equation}\label{eq:interpolation-3}
 	\norm{\nabla {\phi}(u)}^{p^+-p^-}_{L^\infty(T)} \lesssim 
 	\max\left\{ 1, \left( \fint_{C_T} \abs{\nabla u} \right)^{p^+ - p^-} \right\} \norm{\nabla u}_{L^{p^-}(C_T)}^{p^-}
 \end{equation}
 where the implicit constant at right hand side is universal. 
 By~\eqref{eq:interpolation-2} and~\eqref{eq:interpolation-3},
 \begin{equation}\label{eq:interpolation-4}
 \begin{split}
 \int_T \abs{\nabla {\phi}(u)}^{p(x)}\,{\d}x &\lesssim 
 \max\left\{ \norm{\nabla u}_{L^{p^-}(C_T)}^{p^-}, \mathcal{L}^2(C_T)^{p^-_T - p^+_T} \norm{\nabla u}_{L^{p^-}(C_T)}^{p^+} \right\} \\
  &\lesssim \max\left\{ \norm{\nabla u}_{L^{p^-}(C_T)}^{p^-}, \ell \norm{\nabla u}_{L^{p^-}(C_T)}^{p^+} \right\}
 \end{split}
 \end{equation}
 where the last inequality follows from the log-H\"{o}lder condition~\eqref{eq:log-hol} (applied to 
 the smallest ball $B_T$ containing $C_T$) and $\ell$ is precisely the constant in~\eqref{eq:log-hol}. 
Thus, by~\eqref{eq:interpolation-4} and 
Proposition~\ref{prop:embedding},
\begin{equation}\label{eq:interpolation-5}
\int_T \abs{\nabla {\phi}(u)}^{p(x)} \lesssim
2 \ell (1 + \mathcal{L}^2(B_{2r})) \max\left\{ \norm{\nabla u}_{L^{p(\cdot)}(C_T)}^{p^-}, \norm{\nabla u}_{L^{p(\cdot)}(C_T)}^{p^+} \right\},
\end{equation}
up to a universal constant.
 The inequality~\eqref{eq:interpolation-5} clearly implies~\eqref{eq:p(x)-norm-grad-phi} at the local level.
 Taking the sum of the inequalities~\eqref{eq:interpolation-5} as $T$ varies 
 in $\mathcal{T}$ (recalling again that the sets $C_T$ have finite overlap), 
 we obtain~\eqref{eq:p(x)-norm-grad-phi}, completing the proof of (ii). 
 By assumption~{(H)} and either~\eqref{eq:sobolev-sbv} or Lemma~\ref{lemma:sbv-to-sobolev}, 
 it follows that, 
 \[
 	\int_{B_{R}} \abs{{\phi}(u)}^{p(x)}\,{\d}x < +\infty, 
 \]
 and therefore, by \eqref{eq:p(x)-norm-grad-phi},  
 ${\phi}(u) \in \SBV^{p(\cdot)}\left(B_{2r},\R^k\right) \cap W^{1,p(\cdot)}\left(B_R, \R^k\right)$. 
 This concludes the proof.
\end{step}

\end{proof}

\begin{remark}
On the account of the inequalities~\eqref{eq:modular-norm}, 
inequality~\eqref{eq:p(x)-norm-grad-phi} readily implies 
\begin{equation}\label{eq:p(x)-norm-grad-phi-bis}
	\begin{split}
		\int_{B_R} \abs{\nabla \phi(u)}^{p(x)} \,{\d}x \lesssim 
		\left(1+R^2\right)\max&\left\{ \int_{B_R} \abs{\nabla u}^{p(x)} \,{\d}x, \left( \int_{B_R} \abs{\nabla u}^{p(x)} \,{\d}x \right)^{\frac{p^+}{p^-}},\right.\\
        &\left.\left( \int_{B_R} \abs{\nabla u}^{p(x)} \,{\d}x \right)^{\frac{p^+}{p^-}}\right\}, 
	\end{split}
\end{equation}
where the implicit constant depends only on the log-H\"{o}lder constant of $p$. 
Consequently, since $\phi(u) = u$ in $B_{2r} \setminus \overline{B_R}$, the same 
inequality holds with $2r$ in place of $R$.
\end{remark}

\begin{remark}
    The inequalities~\eqref{eq:p(x)-norm-grad-phi} and~\eqref{eq:p(x)-norm-grad-phi-bis} are far from 
    being sharp, both because of the loss of exponent in Step~\ref{step:interpolation} above and 
    because we did not try to optimise with respect to the constants at right hand side. Nonetheless, 
    they are enough for our purposes in this paper.
\end{remark}

\begin{remark}\label{rk:obstruction}
As in \cite{CFI-ARMA}, Proposition~\ref{thm:local-approx} and, consequently, 
Proposition~\ref{thm:approx-sbv-px-Rk} and Theorem~\ref{thm:approx-sbv-px} below, 
is valid only in 2D. This is due 
to the fact only in 2D a condition like~\eqref{eq:H1-J} is strong enough to 
ensure that the jump set can be avoided by a grid made by simplexes with 
universally bounded edges and angles. 
To obtain a similar outcome 
in the higher dimensional setting, it seems to be really necessary to exclude a 
``small'' set from the domain, along the lines of \cite{CCS}. 
See Appendix~\ref{app:counterex} for more details on these points.
\end{remark}

Next, along the lines of the argument in \cite{CFI-ARMA}, we use 
Proposition~\ref{thm:local-approx} and the $\SBV^{p(\cdot)}$ compactness theorem 
(Corollary~\ref{cor:compactness-SBVpx}) to prove 
Proposition~\ref{thm:approx-sbv-px-Rk} below (which is the counterpart 
of \cite[Proposition~3.2]{CFI-ARMA}). Before stating the theorem, we recall 
from \cite{CFI-ARMA} the following useful lemma.
\begin{lemma}[{\cite[Lemma~3.1]{CFI-ARMA}}]\label{lemma:CFI-ARMA-3.1}
	Let $s \in (0,1)$, $\eta \in (0,1)$, and $\rho > 0$. 
	Let $J$ be a $\mathcal{H}^1$-rectifiable Borel set in $B_{\rho}$ 
	such that $\mathcal{H}^1(J) < \eta(1-s)\rho/2$. 
	Then, for $\mathcal{H}^1$-almost every $x \in J \cap B_{s\rho}$ there 
	exists a radius $r_x \in (0,(1-s)\rho/2)$ such that
	\begin{gather}
		\mathcal{H}^1(J \cap \partial B_{r_x}(x)) = 0, \label{eq:lemma-3.1-CFI-ARMA-1}\\
		\eta r_x \leq \mathcal{H}^1\left(J \cap B_{r_x}(x)\right) \leq 
		\mathcal{H}^1\left(J \cap B_{2r_x(x)}\right) < 2\eta r_x \label{eq:lemma-3.1-CFI-ARMA-2}.
	\end{gather}
\end{lemma}

We recall the a set $E \subset \R^n$ is \emph{$\mathcal{H}^d$-rectifiable} if and only if it is 
countably $\mathcal{H}^d$-rectifiable and $\mathcal{H}^d(E) < +\infty$, 
see \cite[Definition~2.57]{AFP}.
For the reader's convenience, we provide a detailed proof of Lemma~\ref{lemma:CFI-ARMA-3.1}, because 
this gives us the occasion to fix some small drawbacks of that in \cite{CFI-ARMA} and to provide 
some additional details.
\begin{proof}
    Since $J$ is a Borel set with finite $\mathcal{H}^1$-measure, the restricted measure 
    $\mathcal{H}^1 \res J$ is a Radon measure in $\R^2$.
    Fix any $x \in J \cap B_{2s\rho}$. Then, for all $\lambda \in ((1-s)\rho,2(1-s)\rho)$ 
    except at most countably many we have 
    $\mathcal{H}^1\left(J \cap \partial B_{\lambda_x / 2^k}(x)\right) = 0$
    for all $k \in \N$ (see, e.~g., \cite[Example~1.6.3]{AFP}), i.e., 
    for which~\eqref{eq:lemma-3.1-CFI-ARMA-1} holds. 
    Choose any such $\lambda_x \in ((1-s)\rho,2(1-s)\rho)$ 
    and define
    \[
        r_x := \max\left\{ \frac{\lambda_x}{2^k} : k \in \N, \,\,\mathcal{H}^1\left(J\cap B_{\lambda_x/2^k}(x)  \right) \geq \eta \frac{\lambda_x}{2^k} \right\}.
    \]
    Since $J$ is a Borel set in $\R^2$ and it is $\mathcal{H}^1$-rectifiable, 
    the $\mathcal{H}^1$-density of $J$ is $1$ at $\mathcal{H}^1$-almost every $x \in J$ 
    (for instance, by Besicovitch-Marstrand-Mattila theorem, see, e.~g., \cite[Theorem~2.6.3]{AFP}). 
    Clearly, $J \cap B_{2 s\rho}$ is $\mathcal{H}^1$-rectifiable as well.
    Therefore, the set 
    \[
        I_x := \left\{ \lambda_x/2^k : k \in \N, \,\,\mathcal{H}^1\left(J\cap B_{\lambda_x/2^k}(x)  \right) \geq \eta \lambda_x/2^k \right\} 
    \]
    is nonempty for $\mathcal{H}^1$-almost every $x$. Indeed, otherwise we could find a set 
    $E \subset J \cap B_{2s\rho}$ with $\mathcal{H}^1(E) > 0$ such that
    \[
        \forall k \in \N, \qquad \mathcal{H}^1\left( J\cap B_{\lambda_x/2^k}(x) \right) < \eta \lambda_x/2^k.
    \]
    Thus, we would have 
    \[
        \lim_{k \to +\infty} \frac{\mathcal{H}^1\left( J\cap B_{\lambda_x/2^k}(x)\right)}{\lambda_x/2^k} \leq \eta.
    \]
    Since $\eta < 1$, we would obtain that, for any $x \in E$, the
    $\mathcal{H}^1$-density of $J \cap B_{s\rho}$ at $x$ is less than $1$, a contradiction.
    Thus, for $\mathcal{H}^1$-a.e. 
    $x \in J \cap B_{2s \rho}$, the set $I_x$ has a least upper bound, which is actually a 
    maximum, because the only possible accumulation point of $I_x$ is 0. 
    If the max in the definition of $r_x$ is attained for $k \geq 2$, 
    then~\eqref{eq:lemma-3.1-CFI-ARMA-2} holds by definition. 
    If instead $r_x$ is attained for $k = 1$, then $2 r_x = \lambda_x$ 
    and~\eqref{eq:lemma-3.1-CFI-ARMA-2} holds as well, otherwise 
    \[
       \mathcal{H}^1\left( J \cap B_{2s\rho} \right) \geq 
       \mathcal{H}^1\left( J \cap B_{2 r_x}(x)\right) 
        = \mathcal{H}^1\left( J \cap B_{\lambda_x}(x)\right) \geq \eta \lambda_x > 
        \eta(1-s)\rho,
    \]
    a contradiction. Thus, \eqref{eq:lemma-3.1-CFI-ARMA-2} is verified. This concludes the proof.
\end{proof}

\begin{prop}\label{thm:approx-sbv-px-Rk}
	Let $p : \R^2 \to (1,\infty)$ be a bounded, log-H\"older-continuous,
	variable exponent satisfying $p^- > 1$ and let $k \in \N$.
	There exist universal constants $\xi$,  
	$\eta > 0$ 
    such that the following holds.
    Provided that Assumption~{(${\rm H}$)} is satisfied, then
    for any $s \in (0,1)$ and any $u \in \SBV^{p(\cdot)}\left(B_{\rho}, \R^k\right)$ satisfying
	\[
		\mathcal{H}^1(J_u \cap B_\rho) < \eta(1-s) \frac{\rho}{2},
	\]
	there are a countable family 
	$\mathcal{F} = \{B\}$ of closed balls, overlapping at most $\xi$ times, of 
	radius $r_B < (1-s)\rho$ and centre $x_B \in \overline{B_{s \rho}}$, and 
	a function $w \in \SBV^{p(\cdot)}\left(B_\rho, \R^k\right)$ such that
	\begin{enumerate}[(i)]
  		\item\label{Rk-approx-item-i} $\cup_{\mathcal{F}} B \subset B_{\frac{1+s}{2}\rho }$ and 
		$\frac{1}{\rho}\sum_\mathcal{F} \mathcal{L}^2(B) + \sum_{\mathcal{F}} \mathcal{H}^1(\partial B) \leq \frac{2\pi\xi}{\eta} \mathcal{H}^1(J_u \cap B_\rho)$. Moreover,
        \begin{equation}\label{eq:estimate-bad-balls-Rk}
            \sum_\mathcal{F} \mathcal{L}^2(B) \leq \min\left\{ \frac{\pi\xi}{\eta} \rho \mathcal{H}^1(J_u \cap B_\rho), \pi \left( \frac{\xi}{\eta} \mathcal{H}^1(J_u \cap B_\rho) \right)^2 \right\}
        \end{equation}
		\item\label{Rk-approx-item-ii} $\mathcal{H}^1\left(J_u \cap \cup_{\mathcal{F}} \partial B\right) = \mathcal{H}^1\left(\left(J_u \cap B_{s \rho}\right) \setminus \cup_{\mathcal{F}} B\right) = 0$.
		\item\label{Rk-approx-item-iii} $w = u$ $\mathcal{L}^2$-a.e. on 
		$B_\rho \setminus \cup_{\mathcal{F}} B$.
		\item\label{Rk-approx-item-iv} $w \in W^{1,p(\cdot)}\left(B_{s \rho}; \R^k\right)$ and 
		$\mathcal{H}^1(J_w \setminus J_u) = 0$.
		\item\label{Rk-approx-item-v} For each $B \in \mathcal{B}$, one has $w \in W^{1,p(\cdot)}\left(B, \R^k\right)$, 
		with 
		\begin{equation}\label{eq:-nabla-w-q}
			\int_B \abs{\nabla w}^q \,{\d}x \lesssim 
			\int_B \abs{\nabla u}^q \,{\d}x 
		\end{equation}
		for any $q \in [1,p^-]$, where the implicit constant at right hand side is universal. Moreover,
		\begin{equation}\label{eq:nabla-w-px}
			\int_B | \nabla w |^{p(x)} \,{\d}x \lesssim \left(1+r^2\right)
			\max\left\{\norm{\nabla u}_{L^{p(\cdot)}(B)}^{p^-}, \norm{\nabla u}_{L^{p(\cdot)}(B)}^{p^+}\right\}
		\end{equation}
		where the implicit constant in the right hand side depends only on the log-H\"{o}lder constant 
        of $p$.
		\item\label{Rk-approx-item-vi} If, in addition, $u \in L^\infty\left(B_\rho; \R^k\right)$, then 
        $w \in L^\infty\left(B_\rho; \R^k\right)$ with
		\[
			\| w\|_{L^\infty(B_\rho)} \leq \|u\|_{L^\infty(B_\rho)}.
		\] 
	\end{enumerate}
\end{prop}

\begin{proof}
With Proposition~\ref{thm:local-approx} at hand, the proof is completely 
analogous to that of \cite[Proposition~3.2]{CFI-ARMA}. 
We sketch the key steps to point out the main differences, addressing the 
reader to \cite{CFI-ARMA} for the missing details. 

Let $s \in (0,1)$ be fixed and set $J = J_u$. Then, $J$ is a Borel set in 
$B_\rho$, $\mathcal{H}^1$-rectifiable, and such that 
$\mathcal{H}^1(J \cap B_\rho) < \eta (1-s)\rho/2$. 
According to Lemma~\ref{lemma:CFI-ARMA-3.1}, the set $E$ of points 
$x \in J_u \cap B_{s\rho}$ at which at least one 
between~\eqref{eq:lemma-3.1-CFI-ARMA-1} 
or~\eqref{eq:lemma-3.1-CFI-ARMA-2} does not hold is a 
$\mathcal{H}^1$-null subset of 
$J_u \cap B_{s\rho}$.
Then, removing $E$ from $J_u \cap B_{s\rho}$, we have 
$\mathcal{H}^1((J_u \cap B_{s\rho}) \setminus E) = 0$ and, moreover, a function 
$r : (J_u \cap B_{s\rho}) \setminus E \to (0,(1-s)\rho/2)$ is defined everywhere 
on $(J_u \cap B_{s\rho}) \setminus E$. By Besicovitch covering theorem 
\cite[Theorem~2.18 and Theorem~2.17]{AFP} applied to 
$(J \cap B_{s\rho}) \setminus E$, we 
can find a universal number $\xi \in \N$ of countable families $\mathcal{F}'_j$ of disjoint closed 
balls $\left\{B_j^i\right\}_{i \in \N}$ (for $j \in \{1,\dots,\xi\}$), of 
radius $r_{B_j^i} < (1-s)\rho/2$ and centre 
$x_{B_j^i} \in \overline{B_{s\rho}}$, 
that cover $\mathcal{H}^1$-almost all $J_u \cap B_{s \rho}$. 
Upon setting $\mathcal{F} := \cup_{j=1}^\xi \mathcal{F}'_j$, we have 
$\mathcal{H}^1((J_u\cap B_{s\rho}) \setminus \cup_{\mathcal{F}} B) = 0$, 
and, again as in the proof of \cite[Proposition~3.2]{CFI-ARMA},
\[
	\sum_{B \in \mathcal{F}} \mathcal{H}^1(\partial B) \leq 
	\frac{2\pi \xi}{\eta}\mathcal{H}^1(J_u \cap B_{s\rho}).
\]
Moreover, since all the radii $r_B$ of the balls $B$ in $\mathcal{F}$ are 
smaller than $\rho$, we get  
$\sum_{\mathcal{F}} r_B^2 \leq \rho \sum_{\mathcal{F}} r_B$, yielding 
$\mathcal{L}^2\left(\cup_{\mathcal{F}} B\right) \leq \frac{\pi \xi}{\eta}\rho \mathcal{H}^1\left( J_u \cap B_{2s\rho} \right)$. Furthermore, 
since $\Sigma_{\mathcal{F}} r_B^2 \leq \left( \Sigma_{\mathcal{F}} r_B \right)^2$, 
\eqref{eq:estimate-bad-balls-Rk} follows, too. By~\eqref{eq:lemma-3.1-CFI-ARMA-1}, 
there holds $\mathcal{H}^1(J_u \cap \cup_{\mathcal{F}} \partial B) = 0$, and this 
completes the proof of~\ref{Rk-approx-item-ii}.

On each ball $B^i_1$ of the first family $\mathcal{F}'_1$, we consider the 
function $\phi^i_1(u)$ given by Proposition~\ref{thm:local-approx}. For each 
$h \in \N$, we define
\[
	w^h_1 := \begin{cases}
	\phi^i_1(u) & \mbox{on } B_h^i, \,\, i\leq h, \\
	u &\mbox{otherwise}.
	\end{cases}
\] 
Properties {(ii)} and {(iv)} in Proposition~\ref{thm:local-approx} imply that 
\[
	w^h_1 \in \SBV^{p(\cdot)}\left(B_{s\rho}, \R^k\right), \quad 
	w^h_1 \in W^{1,p(\cdot)}\left(\cup_{i \leq h} B_1^i\right), \quad
	w^h_1 =u \quad \mathcal{L}^2\mbox{-a.e. on } B_{\rho} \setminus 
	\cup_{i \in \N} B_1^i ,
\]
and
\begin{equation}\label{eq:H1-w_h=0}
	\mathcal{H}^1(J_{w_1^h} \setminus J_u) = 0\,.
\end{equation}
In addition, by {(ii)} in Proposition~\ref{thm:local-approx}, it also follows that, 
for any $h \in \N$, 
\begin{equation}\label{eq:glob-approx-1}
	\int_{B_\rho} \abs{\nabla w^h_1}^{p(x)}\,{\d}x \lesssim \left(1+\rho^2\right)
	\max\left\{ \norm{\nabla u}_{L^{p(\cdot)}(B_\rho)}^{p^-}, \norm{\nabla u}_{L^{p(\cdot)}(B_\rho)}^{p^+} \right\}
\end{equation}
up to a constant depending only on the log-H\"{o}lder constant of $p(\cdot)$. 
Moreover, for any $h \in \N$ there holds
\[
	\abs{{\D}w_1^h}(B_\rho) \leq \abs{{\D}u}\left( B_\rho \setminus \cup_{i\leq h} B_1^i \right) 
	+ \widetilde{c} \int_{\cup_{i\leq h} B_1^i} \abs{\nabla u}\,{\d}x,
\]
for a universal constant $\widetilde{c} > 0$. 
Furthermore, by~{(iii)} of Proposition~\ref{thm:local-approx}, 
\[
	\norm{w^h_1 - w^{h-1}_1}_{L^1(B_\rho)} 
	= \norm{w^h_1-u}_{L^1(B^h_1)} 
	\leq c_k \rho \abs{{\D}u}(B_1^h),
\]
and, for any $j \geq h \geq 1$,
\[
	\norm{w^h_1 - w^j_1}_{L^1(B_\rho)} \leq c_k \rho \abs{{\D}u}\left(B_1^h\right),
\]
where $c_k > 0$ depends only on $k$.
Thus, as in \cite[Proposition~3.2]{CFI-ARMA}, we see that the 
sequence $\left\{w^h_1\right\}$ converges to $w_1$ in $L^1\left(B_{\rho},\R^k\right)$ 
as $h \to +\infty$, where 
\[
	w_1 := 
	\begin{cases}
		\phi(u) & \mbox{on } \cup_{i \in \N} B^i_1, \\
		u & \mbox{otherwise},
	\end{cases}
\]
and $\phi(u) := \sum_{i \in \N } \phi^i_1(u) \chi_{B_1^i}$. By {(ii)}, {(iii)} 
in Proposition~\ref{thm:local-approx}, it follows that 
$\phi(u) \in \SBV^{p(\cdot)}\left(B_{\rho}, \R^k\right) \cap W^{1,p(\cdot)}\left(\cup_{\mathcal{F}_1'} B, \R^k\right)$. 
Therefore, by \eqref{eq:glob-approx-1}, \eqref{eq:H1-w_h=0}, and Corollary~\ref{cor:compactness-SBVpx},
we have $w_1 \in \SBV^{p(\cdot)}\left(B_{\rho},\R^k\right)$. Exactly as in 
\cite[end of p.~1349]{CFI-ARMA}, since 
\[
	\mathcal{H}^1\left( J_{w^h_1} \cap \cup_{i \in \N} B^i_1 \right) = 
	\mathcal{H}^1\left(J_u \cap \cup_{i \geq h+1} B_1^i \right), 
\]
we conclude that
\[
	\mathcal{H}^1 \left( J_{w_1} \cap \cup_{i \in \N} B^i_1 \right) \leq 
	\liminf_{h\to +\infty} \mathcal{H}^1 \left( J_{w^h_1} \cap \cup_{i \in \N} B^i_1 \right) = 0
\]
and therefore, by~\eqref{eq:sobolev-sbv} or Lemma~\ref{lemma:sbv-to-sobolev} (depending on which of 
the two conditions in Assumption~{(${\rm H}$)} is verified), it follows 
that $w_1 \in W^{1,p(\cdot)}\left(\cup_{i \in \N} B^i_1\right)$.
In addition, by 
construction, $w_1 = u$ $\mathcal{L}^2$-a.e. on $B_{\rho} \setminus \cup_{\mathcal{F}_1'} B$ and $\mathcal{H}^1(J_{w_1} \setminus J_u) = 0$.

Iterating the previous construction for any integer $l$ with $1 < l \leq \xi$, 
considering the $l$-th family $\mathcal{F}'_l$, the sequence 
\[
	w^h_l := \begin{cases}
		\phi^i(w_{l-1}) & \mbox{on } B^i_l,\, i\leq h,\\
		w_{l-1} & \mbox{otherwise},
	\end{cases}
\]
and its $L^1$-limit
\[
	w_l := \begin{cases}
		\phi(w_{l-1}) & \mbox{on } \cup_{i \in \N} B^i_l, \\
		w_{l-1} & \mbox{otherwise},
	\end{cases} 
\]
with $\phi(w_{l-1}) := \sum_{i \in \N} \phi^i(w_{l-1}) \chi_{B_l^i}$,
we see that there holds 
\[
	w_l \in \SBV^{p(\cdot)}\left(B_\rho,\R^k\right), \quad  
	w_l \in W^{1,p(\cdot)}\left(\cup_{j \leq l} \cup_{\mathcal{F}'_j} B, \R^k\right), \quad  
	w_l = w_{l-1} \,\, \mathcal{L}^2\mbox{-a.e. on } B_\rho \setminus \cup_{\mathcal{F}'_l} B, 
\]
and
\[  
	\mathcal{H}^1(J_{w_l} \setminus J_u) = 0.
\] 
Now, setting $w := w_\xi$, we have $w \in \SBV^{p(\cdot)}\left(B_\rho, \R^k\right)$, 
$w \in W^{1,p(\cdot)}\left(\cup_{\mathcal{F}} B, \R^k\right)$, $w = u$ $\mathcal{L}^2$-a.e. 
on $B_\rho \setminus \cup_{\mathcal{F}} B$, and $\mathcal{H}^1(J_w \setminus J_u) = 0$.
Thus, the conclusion follows exactly as in \cite[Proposition~3.2]{CFI-ARMA}.
\end{proof}

\begin{remark}
    In the special case in which $p$ is constant, then~\eqref{eq:nabla-w-px} reduces 
    to~\eqref{eq:-nabla-w-q} (in particular, as observed in Remark~\ref{rk:bad-scaling}, the factor 
    $1+r^2$ at right hand side in~\eqref{eq:nabla-w-px} is not present in this case). 
\end{remark}

\begin{remark}
    As an obvious consequence of items {(i)} and {(iii)} of Theorem~\ref{thm:approx-sbv-px-Rk} and 
    of~\eqref{eq:modular-norm}, it follows that
    \begin{equation}\label{eq:modular-w}
    \begin{split}
        \int_{B_{\rho}} &\abs{\nabla w}^{p(x)}\,{\d}x \lesssim 
        \left(1+\rho^2\right) \max\left\{\norm{\nabla u}_{L^{p(\cdot)}(B_{\rho})}^{p^-}, \norm{\nabla u}_{L^{p(\cdot)}(B_{\rho})}^{p^+}  \right\} \\
        &\lesssim \left(1+\rho^2\right) \max\left\{ \int_{B_{\rho}} \abs{\nabla u}^{p(x)}\,{\d}x, \left(\int_{B_{\rho}} \abs{\nabla u}^{p(x)}\,{\d}x\right)^\frac{p^-}{p^+}, \left(\int_{B_{\rho}} \abs{\nabla u}^{p(x)}\,{\d}x\right)^\frac{p^+}{p^-}  \right\},
    \end{split}
    \end{equation}
    where the implicit constant at right hand side depends only on the log-H\"{o}lder constant of 
    $p(\cdot)$.
\end{remark}

\subsection{The constrained case}\label{sec:approx-M}

In this section, we prove Theorem~\ref{thm:approx-sbv-px}. Our argument is reminiscent of that 
of \cite[Theorem~6.2]{HardtLin} and it combines  
Proposition~\ref{thm:approx-sbv-px-Rk} and the following topological property (see, e.g., 
\cite[Lemma~6.1]{HardtLin}, \cite[Proposition~2.1]{BousquetPonceVanSchaftingen}, \cite[Proposition~4.5]{Hopper}, \cite[Lemma~9]{CanevariOrlandi-TopII}).
We recall that, given a topological space $A$ and 
$B \subseteq A$ a subset of $A$, a \emph{retraction} of $A$ onto $B$ is a continuous map $\varrho : A \to B$ satisfying 
$\varrho(z) = z$ for any $z \in B$.
\begin{lemma}\label{lemma:retraction}
    If $\mathcal{M}$ is a compact $m$-dimensional smooth submanifold of $\R^k$ with 
    $\pi_0(\mathcal{M}) = \pi_1(\mathcal{M}) = \dots = \pi_j(\mathcal{M}) = 0$ then there exist a 
    compact $(k-j-2)$-dimensional smooth complex $\mathcal{X}$ in $\R^k$ and a locally smooth 
    retraction $\mathcal{P} : \R^k \setminus \mathcal{X} \to \mathcal{M}$ so that
    \begin{equation}\label{eq:nabla-retraction}
        \abs{\nabla \mathcal{P}(y)} \leq \frac{C}{\dist\left(y, \mathcal{X}\right)}
    \end{equation}
    for any $y \in \R^k \setminus \mathcal{X}$ and some constant $C=C(\mathcal{M}, k, \mathcal{X})$.
    Moreover, in a neighborhood of $\mathcal{M}$, $\mathcal{P}$ is smooth of constant rank $m$.
\end{lemma}

\begin{remark}
    By construction, $\mathcal{X}$ is strictly away from $\mathcal{M}$. Apparently, 
    Lemma~\ref{lemma:retraction} was firstly proven by Hard and Lin \cite[Lemma~6.1]{HardtLin}, with 
    $\mathcal{X}$ a Lipschitz polyhedron and $\mathcal{P}$ locally Lipschitz retraction. 
    Later, it was realised that the same argument, at price of very minor complications, allows to 
    construct $\mathcal{X}$ as a smooth complex and $\mathcal{P}$ as a locally smooth retraction, 
    c.f., e.g., \cite{BousquetPonceVanSchaftingen, CanevariOrlandi-TopI}.
\end{remark}

\begin{remark}\label{rk:integrability-proj}
Up to a bounded change of metric in $\R^k$, \eqref{eq:nabla-retraction} implies
    \begin{equation}\label{eq:retraction}
        \int_{B_R} \abs{\nabla \mathcal{P}(x)}^p\,{\d}x < +\infty \quad \mbox{whenever}\quad 
        1 \leq p < j+2 \quad \mbox{and}\quad 0 < R < +\infty,
    \end{equation}
see, e.g., \cite[Lemma~6.1]{HardtLin} or \cite[Lemma~14]{CanevariOrlandi-TopI}.
\end{remark}

We import almost verbatim a useful remark from~\cite{CanevariOrlandi-TopI}.
\begin{remark}\label{rk:TFI}
    Notice that $\mathcal{P}_a \vert_{\mathcal{M}} : y \in \mathcal{M} \mapsto \mathcal{P}(y-a)$, for 
    $\abs{a}$ small enough, defines a smooth family of maps $\mathcal{M} \to \mathcal{M}$ such that 
    $\mathcal{P}_0 \vert_{\mathcal{M}} = \mathcal{P} \vert_{\mathcal{M}} = {\rm Id}_{\mathcal{M}}$. Therefore, the implicit function theorem implies that 
    $\mathcal{P}_a \vert_{\mathcal{M}}$ has a smooth inverse 
    $\left( \mathcal{P}_a \vert_{\mathcal{M}}\right)^{-1} : \mathcal{M} \to \mathcal{M}$ for 
    $\abs{a}$ sufficiently small (depending only on $\mathcal{M}$ and $\mathcal{P}$).
\end{remark}

The following observation, explicitly pointed out in the proof of \cite[Lemma~14]{CanevariOrlandi-TopI} 
(but already implicitly used in \cite[Theorem~6.2]{HardtLin}, \cite[Lemma~2.3]{HKL}), will be useful as 
well.
\begin{lemma}\label{lemma:HKL}
    For any positive numbers $\sigma$, $\Lambda$, any $v \in L^\infty\left(\Omega,\R^k\right)$ with 
    $\norm{v}_{L^\infty(\Omega)} \leq \Lambda$, any measurable function $g : \Omega \to [0,+\infty)$ 
    and any Borel function $f : \R^k \to [0,+\infty)$, there holds
    \begin{equation}\label{eq:fubini}
        \int_{\mathcal{B}_\sigma^k}\left( \int_\Omega g(x) f(v(x)-y) \,{\d}x \right)\,{\d}y \leq 
        \int_\Omega g(x)\,{\d}x \int_{B^k_{\sigma+\Lambda}} f(z)\,{\d}z
    \end{equation}
\end{lemma}
\begin{proof}
    Just apply Fubini's theorem and make the change of variable $z := v(x) - y$ in the integral with 
    respect to $y$.
\end{proof}

Next, we combine Remark~\ref{rk:traces} and  Lemma~\ref{lemma:retraction} to obtain the following 
auxiliary result. (A related statement is proven in~\cite[Lemma~5]{DeFilippis}.)
\begin{lemma}\label{lemma:projection}
    Let $\Omega \subset \R^n$ be a bounded open set with Lipschitz boundary. 
    Let $j \geq 0$ be an integer, let $k \in \N$ and
    assume that $\mathcal{M}$ is a compact, $j$-connected, smooth Riemannian manifold without boundary, 
    isometrically embedded into $\R^k$.
    Let
    $p : \Omega \to (1,+\infty)$ be a bounded, 
    variable exponent, satisfying 
    \[
        1 < p^- \leq p(x) \leq p^+ < j + 2
    \]
    for every $x \in \Omega$.
    Then, for any $w \in \left(L^\infty \cap W^{1,p(\cdot)}\right)\left(\Omega, \R^k \right)$ so that  
    $w \vert_{\partial \Omega}$ takes values in $\mathcal{M}$, there exists 
    $\widetilde{w} \in W^{1,p(\cdot)}\left(\Omega, \mathcal{M}\right)$ such that $w = \widetilde{w}$ 
    a.e. on $\{ x \in \Omega : w(x) \in \mathcal{M} \}$ and $w = \widetilde{w}$ on $\partial \Omega$ 
    (in the sense of traces). Moreover, 
    \begin{equation}\label{eq:px-gradient-projection}
        \int_{\Omega} \abs{\nabla \widetilde{w}}^{p(x)}\,{\d}x \leq C \int_{\Omega} \abs{\nabla w}^{p(x)}\,{\d}x,
    \end{equation}
    where $C > 0$ is a constant independent of $w$ and $\widetilde{w}$ and that depends only on the 
    quantities listed in Remark~\ref{rk:projection}. 
\end{lemma}

The proof of Lemma~\ref{lemma:projection} follows very closely arguments and computations in that 
of~\cite[Theorem~6.2]{HardtLin}, taking advantage of several observations 
in~\cite{CanevariOrlandi-TopI}. 
For the reader's convenience, we provide full details. 
During the proof, we will keep track of the precise dependencies of the 
constant $C$ in~\eqref{eq:px-gradient-projection}, gathering them in Remark~\ref{rk:projection} below.

\begin{proof}
    Let $\mathcal{P} : \R^k \setminus \mathcal{X} \to \mathcal{M}$ be given by 
    Lemma~\ref{lemma:retraction}. For $a \in \R^k$, define $\mathcal{X}_a := \left\{ y + a : y \in \mathcal{X} \right\}$ and $\mathcal{P}_a : \R^k \setminus \mathcal{X}_a \to \mathcal{M}$ by 
    $\mathcal{P}_a(y) : y \in \R^k \setminus \mathcal{X}_a \mapsto \mathcal{P}(y-a)$.  
    By Lemma~\ref{lemma:retraction}, the map $\mathcal{P}_a$ is well-defined and locally smooth and,  
    by Remark~\ref{rk:TFI}, we may find $\sigma > 0$ so small that the restricted map 
    $\mathcal{P}_a \vert_{\mathcal{M}} : \mathcal{M} \to \mathcal{M}$ is a diffeomorphism for all
    $a \in \mathcal{B}^k_\sigma := \left\{ y \in \R^k : \abs{y} < \sigma \right\}$. 
    Up to further reducing $\sigma$ below a threshold value 
    $\bar{\sigma}$ (depending only on $\mathcal{P}$ and $\mathcal{M}$), 
    if necessary, the inverse function theorem implies that
    \[
        \lambda := \sup_{a \in \mathcal{B}_\sigma} \Lip\left(\mathcal{P}_a \vert_{\mathcal{M}}\right)^{-1}
    \]
    is a finite number depending only on $\mathcal{P}$ and $\mathcal{M}$.
    
    Now, we notice that the set
    \[
        N := \left\{ (x,y) \in \Omega \times \R^k : w(x) - y \in \mathcal{X} \right\} 
    \]
    is measurable and $\mathcal{H}^{n+k-j}(N) = 0$, because each slice $N \cap \left( \{x\} \times \R^k \right) = w(x) - \mathcal{X}$ ha dimension $k-j-2$. By Fubini's theorem, it follows that 
    $\mathcal{H}^n\left(N \cap \left( \Omega \times \{a\}\right) \right) = 0$ for a.e. $a \in \R^k$, 
    so $\mathcal{P}_a\circ w$ is well-defined for a.e. $a \in \R^k$ and it is a measurable function. 
    Moreover, by the chain rule, for a.e. $(x,y) \in \Omega \times \R^k$ we have
    \[
        \abs{\nabla\left(\mathcal{P}_a \circ w\right)} = \abs{\nabla\left(\mathcal{P} \circ (w(x)-a) \right)} = \abs{(\nabla \mathcal{P})(w(x)-a)}\abs{\nabla w(x)}.
    \]
    Applying Lemma~\ref{lemma:HKL} with $v = w$, $f = \abs{\nabla \mathcal{P}}^{p(x)}$,  
    $g = \abs{\nabla w}^{p(x)}$, $\Lambda = \norm{w}_{L^\infty(\Omega)}$ and $\sigma > 0$ as above, 
    we obtain
    \[
        \begin{split}
            \int_{B_\sigma^k} \int_\Omega \abs{\nabla(\mathcal{P}_a \circ w)(x)}^{p(x)}\,{\d}x\,{\d}a &\leq \int_{\Omega} \abs{\nabla w(x)}^{p(x)} \int_{\mathcal{B}^k_{\sigma+\Lambda}} \abs{(\nabla \mathcal{P}_a)(w(x))}^{p(x)}\,{\d}a\,{\d}x \\
            &\leq \int_{\Omega} \abs{\nabla w(x)}^{p(x)} \int_{\mathcal{B}^k_{\sigma+\Lambda}} \abs{(\nabla \mathcal{P})(w(x)-a)}^{p(x)}\,{\d}a\,{\d}x \\
            &\leq \int_{\Omega} \abs{\nabla w(x)}^{p(x)} \left( \int_{\mathcal{B}_\sigma^k} \abs{(\nabla \mathcal{P})(y)}^{p(x)}\,{\d}y\right)\,{\d}x \\
            &\leq C' \int_{\Omega} \abs{\nabla w(x)}^{p(x)}\,{\d}x < +\infty,
        \end{split}
    \]
    where in the last line we used the obvious inequality 
    $\abs{\nabla \mathcal{P}(y)}^{p(x)} \leq 1 + \abs{\nabla \mathcal{P}(y)}^{p^+}$ (which holds 
    for almost every $(x,y) \in \Omega \times \R^k$) and 
    Remark~\ref{rk:projection}. Thus, the constant $C'$ depends only on $p^+$, 
    $\mathcal{M}$ and $\mathcal{P}$ (through $j$, $k$, 
    and  
    $\int_{B_{\bar{\sigma}}} \abs{\nabla \mathcal{P}}^{p^+} \,{\d}y$, where the latter is finite 
    by~\eqref{eq:retraction}, since $p^+ < j+2$). 
    Therefore, by Chebyshev inequality, we may choose $a \in \mathcal{B}_\sigma^k$ so that
    \begin{equation}\label{eq:nabla-Pa-w}
        \int_{\Omega} \abs{\nabla (\mathcal{P}_a \circ w)(x)}^{p(x)}\,{\d}x \leq 
        C' \abs{\mathcal{B}_\sigma^{k}}^{-1} 
        \int_{\Omega} \abs{\nabla w(x)}^{p(x)}\,{\d}x,
    \end{equation}
    and the right hand side does not depend on $a$.
    We conclude that the map 
    \[
        \widetilde{w} := \left(\mathcal{P}_a \vert_{\mathcal{M}}\right)^{-1} \circ \mathcal{P}_a \circ w
    \]
    belongs to $W^{1,p(\cdot)}(\Omega, \mathcal{M})$ and it agrees with $w$ a.~e. 
    on the set $\{ x \in \Omega : w(x) \in \mathcal{M} \}$ and on $\partial \Omega$ in the sense of traces (because $\mathcal{P}_a$ is a retraction onto $\mathcal{M}$ and $\left(\mathcal{P}_a \vert_{\mathcal{M}}\right)^{-1} \circ \mathcal{P}_a$ is the identity on $\mathcal{M}$). 
    Finally, ~\eqref{eq:px-gradient-projection} follows from~\eqref{eq:nabla-Pa-w}, 
    with $C := \max\left\{\lambda^{p^+},\lambda^{p^-}\right\}  C' \abs{\mathcal{B}_\sigma^{k}}^{-1} $ depending only on $p^+$, $p^-$, 
    $\mathcal{M}$, $\mathcal{P}$, 
    and $k$. 
\end{proof}

\begin{remark}\label{rk:projection}
    The constant $C$ in~\eqref{eq:px-gradient-projection} depends only on $p^+$, $p^-$, 
    $k$, $\mathcal{M}$, and $\mathcal{P}$.
\end{remark}

We are now ready for the proof of Theorem~\ref{thm:approx-sbv-px}.

\begin{proof}[Proof of Theorem~\ref{thm:approx-sbv-px}]
    Given $u \in \SBV^{p(\cdot)}\left(B_\rho, \mathcal{M} \right)$ satisfying~\eqref{eq:small-jump-set}, for any fixed $s \in (0,1)$, let $v \in W^{1,p(\cdot)}\left(B_{s\rho}, \R^k\right) \cap \SBV^{p(\cdot)}\left(B_{\rho}, \R^k\right)$ and $\mathcal{F}$ be, respectively, the  
    function and the family of balls provided by Proposition~\ref{thm:approx-sbv-px-Rk}.
    Since $u$ takes values in $\mathcal{M}$, by \eqref{eq:M-in-a-ball} it follows  
    $\norm{u}_{L^\infty(B_{\rho})} \leq M$ and
    by item~\ref{Rk-approx-item-vi} of 
    Proposition~\ref{thm:approx-sbv-px-Rk} we have, in addition,  
    $\norm{v}_{L^\infty(B_{\rho})} \leq M$. 
    Set $v_s := v \vert_{B_{s \rho}}$. 
    Clearly, $v_s \in W^{1,p(\cdot)}\left(B_{s\rho}, \R^k\right)$, 
    $\norm{v_s}_{L^\infty(B_{s\rho})} \leq M$, and 
    $v_s \vert_{\partial B_{s \rho}}$ takes values in $\mathcal{M}$. 
    By Lemma~\ref{lemma:projection} with $j = 0$, 
    a function $w_s \in W^{1,p(\cdot)}\left(B_{s\rho}, \mathcal{M}\right)$ exists such that 
    $w_s = v_s$ a.~e. on the set $\{ x \in \Omega : v_s(x) \in \mathcal{M} \}$ and such that 
    \begin{equation}\label{eq:same-traces-on-Bs}
        w_s \vert_{\partial \Omega} = v_s \vert_{\partial \Omega} = v \vert_{\partial \Omega}.
    \end{equation}
    (In particular, $w_s \vert_{\partial \Omega}$ takes values in $\mathcal{M}$.) Moreover, 
    by~\eqref{eq:px-gradient-projection},
    \begin{equation}\label{eq:px-grad-w_s}
        \int_{B_{s\rho}} \abs{\nabla w_s}^{p(x)}\,{\d}x \leq C \int_{B_{s\rho}} \abs{\nabla v_s}^{p(x)}\,{\d}x,
    \end{equation}
    where, by Remark~\ref{rk:projection}, $C > 0$ is a constant that depends only on $p^+$, $p^-$, 
    $k$, $\mathcal{M}$, $\mathcal{P}$, 
    where $\mathcal{P}$ is the 
    retraction provided by Lemma~\ref{lemma:retraction}. 
    Upon setting
    \[
        w := 
        \begin{cases}
            u & \mbox{in } B_\rho \setminus \overline{B_{s\rho}}\\
            w_s & \mbox{in } B_{s\rho},
        \end{cases}
    \]
    we see that $w \in W^{1,p(\cdot)}\left( B_{s \rho}, \mathcal{M} \right) \cap \SBV^{p(\cdot)}\left(B_\rho, \mathcal{M} \right)$ and that $\mathcal{H}^1(J_w \setminus J_u) = 0$ (indeed, 
    in view of~\eqref{eq:same-traces-on-Bs}, $\mathcal{H}^1(J_w \setminus J_v) = 0$). Moreover, 
    $w \vert_{\partial \Omega} = w_s \vert_{\partial \Omega}$.
    
    Now, we claim that the function $w$ and the family $\mathcal{F}$ are a function and a family 
    of balls satisfying all the properties required by the statement. 
    Indeed, item~\ref{item:approx-iv} is 
    obvious (as it already follows from Proposition~\ref{thm:approx-sbv-px-Rk}). 
    Item~\ref{item:approx-i} holds because $u = v$ a.~e. on $B \setminus \cup_{\mathcal{F}} B$, 
    so that in particular $v_s = u$ on $B_{s\rho} \setminus \cup_{\mathcal{F}} B$, and $w_s$ 
    coincides with $v_s$ when the latter takes values in $\mathcal{M}$ (in particular, on the set 
    $\{u = v\} \cap B_{s\rho}$). 
    Item~\ref{item:approx-ii} is an immediate consequence of the above discussion.
    Finally, item~\ref{item:approx-iii} follows from \eqref{eq:px-grad-w_s}, \eqref{eq:nabla-w-px}, 
    and the definition of $w$.
\end{proof}

\begin{remark}
    Since generalised functions of bounded deformation are not closed under composition with (locally) 
    Lipschitz or even smooth mappings, there is no immediate extension of Theorem~\ref{thm:approx-sbv-px} to the case of functions belonging, for instance, to $\SBD^{p(\cdot)}\left( B_\rho, \S^1 \right)$.
\end{remark}

\begin{remark}
    For simplicity, we 
    stated Theorem~\ref{thm:approx-sbv-px} assuming that $\mathcal{M}$ is merely 
    connected. However, a quick inspection of the proof shows that it is valid, with identical 
    proof, if we assume $\mathcal{M}$ is $j$-connected and we replace the constraint $p^+ < 2$ 
    with $p^+ < j+2$.
\end{remark}

\section{Regularity of local minimisers}\label{sec:existence}

In this section, we prove Theorem~\ref{thm:B}. We will follow ideas that originated 
in~\cite{DeGiorgiCarrieroLeaci} in the scalar case and that were later adapted to the context of 
sphere valued maps in \cite{CarrieroLeaci}. Differently 
from~\cite{DeGiorgiCarrieroLeaci,CarrieroLeaci}, we will {\em not} use the Sobolev-Poincar\'{e} 
inequality for $\SBV$ functions and we will {\em not} make use of medians or truncations. 
Instead, we will follow the approach of~\cite{CFI-AIHP} and employ the Sobolev approximation from 
Section~\ref{sec:approx}. As explained in the Introduction, the key point lies in proving that the jump 
set of local minimisers is essentially closed. 

Let $\Omega \subset \R^2$ be a bounded open set and let $k \in \N$. In this section, 
we shall need the following variants of the integral functional $F$ in~\eqref{eq:functional}: 
for $A \subset \R^2$ a Borel set, we consider
\begin{equation}\label{eq:F}
    F(u,c,A) := \int_A \abs{\nabla u}^{p(x)}\,{\d}x + c \mathcal{H}^1(J_u \cap A),
\end{equation}
defined on $\SBV_{\rm loc}\left(A, \R^k\right)$, where $c > 0$. Clearly, $F(u,c,A)$ reduces to the functional~\eqref{eq:functional} for $A = \Omega$ and $c = 1$.

\begin{definition}\label{def:loc-min}
    Let $c > 0$ and $\Omega \subset \R^2$ be a bounded, open set. 
    We say that a map $u \in \SBV^{p(\cdot)}\left(\Omega, \S^{k-1}\right)$ is a 
    \emph{local minimiser (among maps with values into $\S^{k-1}$)} 
    of the functional $F(\cdot, c, \Omega)$ iff there holds
    \[
        F(u,c,\Omega) \leq F(v, c, \Omega)
    \]
    for any $v \in \SBV^{p(\cdot)}\left( \Omega, \S^{k-1}\right)$ such that 
    $\{ v \neq u \} \subset\subset \Omega$.
\end{definition}
For every open set $A \subset \R^2$, every $t > 0$ and every 
$u \in \SBV_{\rm loc}\left(\Omega, \R^k\right)$ such that 
$\abs{u} = t$ a.e. in $\Omega$, we define 
\[
    \Phi(u,c,A,t) := \inf\left\{ F(u,c,A) : v \in \SBV_{\rm loc}\left(\Omega,\R^k\right), \,\, v = u \,\mbox{ in } \Omega \setminus \overline{A}, \,\, \abs{v} = t \,\mbox{ a.e. in } \Omega \right\}.
\]
If $\Phi(u,c,A,t) < +\infty$, we define the \emph{deviation from minimality} of $u$ as
\[
    \Dev(u,c,A,t) := F(u,c,A) - \Phi(u,c,A,t).
\]
Equivalently, $\Dev(u,c,A,t)$ can be characterised as the smallest number $\kappa \in [0,+\infty]$ 
so that 
\[
    F(u,c,A) \leq F(v,c,A) + \kappa
\]
for all $v \in \SBV_{\rm loc}\left(\Omega,\R^k\right)$ satisfying $\{v \neq u\} \subset\subset A$ 
and such that $\abs{v} = t$ a.e. in $\Omega$. 

For convenience, when $c = 1$, we set
\begin{equation}\label{eq:F-bis}
    F(u,A) := F(u,1,A), \quad \Dev(u,A,t) := \Dev(u,1,A,t), \quad \Dev(u,A) := \Dev(u,1,A,1).
\end{equation}

Similarly to as in~\cite[Section~7]{AFP}, we establish the following definition.
\begin{definition}\label{def:quasi-min}
    A function $u \in \SBV_{\rm loc}\left(\Omega,\S^{k-1}\right)$ is a \emph{quasi-minimiser} 
    (among maps with values into $\S^{k-1}$) of the 
    functional $F(\cdot, \Omega)$ in $\Omega$ if there exists 
    a constant $\kappa' \geq 0$ such that for all $x \in \Omega$ and all balls 
    $B_\rho(x) \subset \Omega$,
    \begin{equation}\label{eq:Dev-quasimin}
        \Dev(u, B_\rho(x)) \leq \kappa' \rho^2.
    \end{equation}
\end{definition}

By a standard comparison argument, we get an immediate upper bound on the local energy of quasi-minimisers, that will be used in the proof of Theorem~\ref{thm:density} below.
\begin{lemma}[Energy upper bound]\label{lemma:energy-upper-bound}
    Let $u$ be any quasi-minimiser of $F(\cdot, \Omega)$ in $\Omega$. 
    Then, for all $x \in \Omega$ and for all balls $B_\rho(x) \subset \Omega$, we have
    \begin{equation}\label{eq:F-upper-bound}
        F(u, B_\rho(x)) \leq 2 \pi \rho + \kappa' \rho^2.     
    \end{equation}
\end{lemma}

\begin{proof}
    The proof is almost identical to that of \cite[Lemma~7.19]{AFP}; the only difference is that we 
    are compelled to consider competitors with values into $\S^{k-1}$. (See also \cite[Lemma~3.2]{CarrieroLeaci} and the remark thereafter.)
    Let $u \in \SBV_{\rm loc}(\Omega, \S^{k-1})$ be a quasi-minimiser of $F(\cdot,\Omega)$ in 
    $\Omega$, take any $\rho' > 0$ with $\rho' < \rho$ and put 
    $v := u \chi_{B_\rho(x) \setminus B_{\rho'}(x)} + {\bf e}_k \chi_{\overline{B_{\rho'}}(x)}$.
    Then, $v \in \SBV_{\rm loc}(B_\rho(x), \S^{k-1})$ and, by quasi-minimality of $u$,
    \[
        F(u,B_{\rho'}(x)) \leq \mathcal{H}^1\left(J_v \cap \overline{B_{\rho'}}\right) + \Dev(u,B_\rho(x)) \stackrel{\eqref{eq:Dev-quasimin}}{\leq} 2 \pi \rho + \kappa' \rho^2.
    \]
    The conclusion follows by letting $\rho' \to \rho$.
\end{proof}

To obtain that $J_u$ is essentially closed, we again argue along the lines of 
\cite{DeGiorgiCarrieroLeaci, CarrieroLeaci, CFI-AIHP, LScSoV}. The heart of matter is proving, in 
any set $\Omega_\delta \subset \Omega$ defined as in~\eqref{eq:Omega-delta},
inequality~\eqref{eq:density-lower-bound} below.

\begin{theorem}[Density lower bounds]\label{thm:density}
    Let $p : \R^2 \to (1,+\infty)$ be a variable exponent satisfying~{(${\rm p}_1'$)} 
    and~{(${\rm p}_2$)}.
    Let $\delta > 0$ and $\Omega_\delta$ be defined by~\eqref{eq:Omega-delta}. 
    There exist $\theta_\delta$ and $\rho_\delta$ depending only on $p^-$, $p^+$, and $\delta$ with the 
    property that, if $u \in \SBV^{p(\cdot)}\left(\Omega,\S^{k-1}\right)$ is a quasi-minimiser of $F$ 
    in $\Omega$, then
    \begin{equation}\label{eq:density-lower-bound}
        F(u, B_\rho(x)) > \theta_\delta \rho
    \end{equation}
    for all balls $B_\rho(x) \subset \Omega_\delta$ with centre $x \in \overline{J_u}$ and radius 
    $\rho < \rho_\delta' := \frac{\rho_\delta}{\kappa'}$, where $\kappa'$ is the constant in~\eqref{eq:F-upper-bound}. 
    Moreover,
    \begin{equation}\label{eq:closed-jump-set}
        \mathcal{H}^1\left(\Omega \cap \left(\overline{J_u} \setminus J_u\right)\right) = 0.
    \end{equation}
\end{theorem}
Once again, the idea of the proof of Theorem~\ref{thm:density} goes back to De~Giorgi, Carriero, 
and Leaci \cite{DeGiorgiCarrieroLeaci}. It is based on proving a power-decay property 
of the energy in small balls with respect to the radius of the balls which, after a 
classical iteration argument, yields the conclusion. Given such a decay 
property, that in our case is provided by Theorem~\ref{thm:decay-lemma} below, 
the proof is purely ``algebraic'' and it does not depend in any way on 
the target manifold of the considered maps. However, some modifications are needed 
with respect to the argument in \cite{DeGiorgiCarrieroLeaci, AFP} and also 
with respect to the recent adaptation to the variable-exponent setting provided 
by~\cite[Theorem~4.7]{LScSoV} because we need to avoid recurring to medians and 
truncations. 
The key technical tool that allows us to do so is provided by Proposition~\ref{prop:not-in-Ju} 
in Appendix~\ref{app:not-in-Ju}, which gives us a sufficient condition to conclude that a given point 
does not belong to the jump set of a map $u$ and that is proved using only the Sobolev approximation 
results from Section~\ref{sec:approx}. 
Once Proposition~\ref{prop:not-in-Ju} is obtained, 
the proof of Theorem~\ref{thm:density} follows by the 
classical argument in~\cite{DeGiorgiCarrieroLeaci}, exactly as in \cite{AFP, LScSoV}.
%
%

The bulk of this section is instead devoted to proving 
Theorem~\ref{thm:decay-lemma}, a task which requires to modify in a nontrivial 
way both the arguments in \cite{LScSoV} and those in \cite{CarrieroLeaci}.

\begin{theorem}[Decay lemma]\label{thm:decay-lemma}
	Let $p : \R^2 \to (1,+\infty)$ be a variable exponent satisfying~{(${\rm p}_1'$)} 
    and~{(${\rm p}_2$)},
    let $\delta > 0$, and let 
    $\Omega_\delta$ be defined as in~\eqref{eq:Omega-delta}.
	There exists a constant $C_{\delta} = C(p^-,p^+,\delta) > 0$ with the following 
	property. For every 
	$\tau \in (0,1)$ there exist $\eps = \eps(\tau,\delta)$ and 
	$\theta = \theta(\tau,\delta)$ in $(0,1)$ such that, if 
	$u \in \SBV(\Omega, \S^{k-1})$ satisfies, for all $x \in \Omega_\delta$ and all 
	$\sigma < \eps^2$ such that $B_\sigma(x) \subset \subset \Omega_\delta$, 
	the conditions 
	\[
		F(u, B_\sigma(x)) \leq \eps \sigma, \qquad 
		\Dev(u, B_\sigma(x)) \leq \theta F(u, B_\sigma(x)),
	\]
	then 
	\begin{equation}\label{eq:decay}
		F(u, B_{\tau \sigma}(x)) \leq C_{\delta} \tau^2 F(u,B_\sigma(x)).
	\end{equation}	 
\end{theorem}
The idea of the proof of Theorem~\ref{thm:decay-lemma} is classical and rooted 
in an argument by contradiction. However, the combined non-standard growth of the 
energy functional and the fact that we deal with maps with values into spheres 
make it somewhat trickier than usual. To handle the first complication, i.e., 
the non-standard growth of the functional, we borrow some ideas from the recent 
work~\cite{LScSoV}, where the functional $F$ is considered on scalar-valued functions of 
class $\SBV^{p(\cdot)}$. 
Note that our argument requires, as in \cite{LScSoV}, the variable exponent to satisfy the  
\emph{strong log-H\"{o}lder condition}, 
c.f.~\eqref{eq:unif-conv-p^0} and Remark~\ref{rk:stong-log-Hol} below. 
  
To take the sphere-valued constraint into account, we reason as in \cite[Section~3]{CarrieroLeaci} but, 
instead of using medians, truncations, and the Sobolev-Poincar\'{e} inequality 
for $\SBV$ functions, we exploit the Sobolev approximation results from 
Section~\ref{sec:approx}. 
Nonetheless, for the proof of Theorem~\ref{thm:decay-lemma} we shall need several auxiliary results in 
the spirit of~\cite{CarrieroLeaci}.
\begin{lemma}[{\cite[Lemma~3.2]{CarrieroLeaci}}]\label{lemma:CaLe-3.2}
    Let $u \in \SBV\left(B_r, \R^k\right)$. For every $c > 0$ and $t > 0$, the functions
    \[
        \rho \mapsto F(u,c, B_\rho) \qquad \mbox{and} \qquad \rho \mapsto \Dev(u,c,B_\rho, t)
    \]
    are non-decreasing in $(0,r)$. 
\end{lemma}
The proof of Lemma~\ref{lemma:CaLe-3.2} is straightforward and, as in~\cite{CarrieroLeaci}, 
left to the reader.

\begin{lemma}\label{lemma:hyperplane}
    Let $B \subset \R^2$ be a ball of radius $s \in [1/2,1)$.
    For any $h \in \N$,  
    let $t_h > 0$, $w_h : B \to \S^{k-1}_{t_h}$ be a measurable function, and set
    $\lambda_h^{(1/2)} := \abs{\avg{w_h}_{B_{1/2}}}$. 
    Assume that $t_h \to +\infty$ as $h \to +\infty$ and that
    \begin{enumerate}[(i)]
        \item There holds 
        \begin{equation}\label{eq:CaLe-type-lemma-1}
           \lim_{h \to +\infty}  t_h\abs{1 - \frac{\lambda_h^{(1/2)}}{t_h}} = d,
        \end{equation}
        for some constant $d \geq 0$.
        \item There exists $w_\infty : B \to \R^k$ such that 
        $w_h - \lambda^{(1/2)}_h \to w_\infty$ a.e. in $B$ as 
        $h \to +\infty$.
    \end{enumerate}
    Then, there exists a hyperplane $\Pi \subset \R^k$ such that 
    $w_\infty(x) \in \Pi$ for $\mathcal{L}^2$-almost every $x \in B$.
\end{lemma}

\begin{proof} 
  The proof is along the lines of that of the corresponding assertion in the proof of
  \cite[Theorem~3.6]{CarrieroLeaci}.
    Up to rotations in the target space, we may assume that for any 
    $h \in \N$ there holds 
    $\avg{w_h}_{B_{1/2}} = \lambda^{(1/2)}_h {\bf e}_k$, where ${\bf e}_k$ 
    denotes the last vector of the canonical basis of $\R^k$. 
    Following \cite{CarrieroLeaci}, we observe that the trivial identity 
    $\abs{w_h - \lambda_h^{(1/2)} + \lambda_h^{(1/2)}}^2 = t_h^2$ yields 
    \[
        \frac{\abs{w_h - \lambda_h^{(1/2)}}^2}{t_h} 
        + 2\left(w_h - \lambda_h\right)\cdot \frac{\lambda_h^{(1/2)}}{t_h} 
        = t_h\left(1-\frac{\lambda_h^{(1/2)}}{t_h}\right) \left( 1+\frac{\lambda_h^{(1/2)}}{t_h}\right).
    \]
    Hence, by~\eqref{eq:CaLe-type-lemma-1} (which implies $\displaystyle\lim_{h \to +\infty} \frac{\lambda_h^{(1/2)}}{t_h} = 1$) and assumption~{(ii)} we obtain 
    $0 + 2 w_0 \cdot {\bf e}_k = 2d$ in $B_1$, so that 
    \[
        d = \avg{w_0 \cdot {\bf e}_k}_{B} = \lim_{h \to +\infty}
        \left( \avg{w_h \cdot {\bf e}_k}_{B} - \lambda_h^{(1/2)} \right) = 0,
    \]
    hence $w_0 \cdot {\bf e}_k = 0$ a.e. in $B$. The conclusion follows.
\end{proof}

For later use, we record here the following trivial linear algebra lemma.
\begin{lemma}\label{lemma:included-hyperplanes}
    Let $\left\{ \Pi^{(s)} \right\}_{s \in [1/2,1)}$ be a family of (proper) hyperplanes in $\R^{k}$, such that $\Pi^{(s)} \leq \Pi^{(t)}$ for any $s$, $t \in [1/2,1)$ with $s \leq t$. Then, there 
    exists $\overline{s} \in [1/2,1)$ such that $\Pi := \Pi^{(\overline{s})}$ contains all the hyperplanes $\Pi^{(s)}$.
\end{lemma}

\begin{proof}
    By assumption, the 
    hyperplane $\Pi^{(1/2)}$ is a vector subspace, of dimension $r \leq k-1$, contained in all 
    the hyperplanes $\Pi^{(s)}$. Let $\left\{ a_1, \dots, a_r\right\}$ be $r$ affinely 
    independent vectors in $\Pi^{(1/2)}$. If $r = k-1$ or $r$ coincide with the maximal 
    dimension of the hyperplanes $\Pi^{(s)}$, we are done. So, let us assume $r < k-1$ 
    and that there is $s > 1/2$ such that $\Pi^{(1/2)}$ is a proper subspace of $\Sigma^{(s)}$. 
    Let $r^{(s)} \leq k-1$ be the dimension of $\Pi^{(s)}$. Then, in $\Sigma^{(s)}$ we can find 
    $r^{(s)}-r$ vectors $a_{r+1}, \dots, a_{r^{(s)}}$ which are affinely independent both each 
    other and from $a_1,\dots,a_r$. Again, if $r^{(s)} = k-1$, we are done, otherwise we repeat the 
    process. In at most $k-r-1$ steps, we get the conclusion. 
\end{proof}

Finally, we recall the following $L^\infty$-$L^1$ estimate for $p$-harmonic functions, see \cite{Uhlenbeck, Hamburger}.
We report here a statement which follows as a special case of a more general 
result proved for differential forms. 

\begin{lemma}\label{lemma:L-infty-L1-estimate}
    Let $u \in W^{1,p}\left(\Omega,\R^k\right)$ be $p$-harmonic in $\Omega$. Let $B_r$ be a ball 
    with $r \in (0,1]$ such that $B_r \subset \Omega$. Then there exists a positive 
    constant $c_0$ which depends only on $n$, $k$, and $p$ such that
    \begin{equation}\label{eq:L-infty-L1-estimate}
        \sup_{B_\frac{r}{2}} \abs{\nabla u}^p \leq C_0 
        \fint_{B_r} \abs{\nabla u}^p\,{\d}y.
    \end{equation}
\end{lemma}

We are now ready to prove Theorem~\ref{thm:decay-lemma}. As already mentioned, the argument combines 
ingredients from~\cite{CarrieroLeaci}, \cite{CFI-AIHP}, and~\cite{LScSoV}.
\begin{proof}[Proof of Theorem~\ref{thm:decay-lemma}]
    Fix $\delta > 0$. 
	It is enough to assume $\tau \in (0,1/2)$ (otherwise, just take $C_{\delta} = 4$). 
	By contradiction, assume that~\eqref{eq:decay} does not hold. Then, there 
	exist sequences
	\[
		\{u_h\} \subset \SBV\left(\Omega,\S^{k-1}\right); \quad 
		\{\eps_h\}, \,\, \{\theta_h\} \subset (0,1),\,\,\{\sigma_h\} \subset (0,+\infty); 
		\quad \{x_h\} \subset \Omega_\delta
	\]
	such that
	\begin{align}
		&\lim_{h \to +\infty} \eps_h = \lim_{h \to +\infty} \theta_h = 0; \label{eq:decay-1.1} \\
		\forall h \in \N, \quad &\sigma_h \leq \eps_h^2 \quad \mbox{and} \quad	
		B_{\sigma_h}(x_h) \subset \subset \Omega_\delta \label{eq:decay-1.2},
	\end{align}
	and, moreover, 
	\begin{gather}
		F(u_h, B_{\sigma_h}(\sigma_h)) \leq \eps_h \sigma_h, \label{eq:decay-2.1}\\
		\Dev(u_h, B_{\sigma_h}(x_h)) \leq \theta_h F(u_h, B_{\sigma_h}(x_h)),\label{eq:decay-2.2} \\
		F(u_h, B_{\tau \sigma_h(x_h)}) > C_1 \tau^2 F(u_h, B_{\sigma_h}(x_h))\label{eq:decay-2.3}, 
	\end{gather}
	where $C_1 = 2 C_0$ and $C_0$ is the 
    constant at right hand side in~\eqref{eq:L-infty-L1-estimate}.
    Furthermore, up to the extraction of a (not relabeled) subsequence, we may assume that $\{x_h\}$ 
    converges to some $x_0 \in \overline{\Omega_\delta}$. Let us define
    \[
        p_0 := p(x_0).
    \]
	
	We now proceed step-by-step, following the argument in \cite{LScSoV}.
	\setcounter{step}{(0)}
	\begin{step}[Scaling, Part~1]
		For all $h \in \N$, we define the translated maps and exponents
		\begin{align*}
            & \widetilde{u}_h : B_1 \to \S^{k-1}, \quad \widetilde{u}_h(y) := u_h(x_h + \sigma_h y),  \\
			& p_h : B_1 \to (1,+\infty), \quad p_h(y) := p(x_h + \sigma_h y), 
		\end{align*}
        where $y \in B_1$. We notice that for any $h \in \N$ and any $y \in B_1$
        \[
            1 < p^- \leq p^-_h \leq p_h(y) \leq p^+_h \leq p^+ < 2.
        \]
        Moreover, each rescaled variable exponent $p_h(\cdot)$ is still strongly log-H\"{o}lder 
        continuous and this implies that the sequence $\{p_h\}$ converges to $p_0$ uniformly 
        in $B_1$ as $h \to +\infty$. Indeed, upon setting
        \[
            p^0_h := p_h(0) = p(x_h),
        \]
        we have
        \begin{equation}\label{eq:unif-conv-p^0}
            \sup_{y \in B_1}\abs{p_h(y) - p_0} \leq \sup_{y \in B_1}\left\{ \abs{p_h(y) - p_h^0} + \abs{p_h^0 - p_0} \right\} = \omega(\sigma_h) + {\rm o}(1),
        \end{equation}
        which tends to zero as $h\to +\infty$ thanks to~\eqref{eq:strong-log-hol}.
        Concerning the translated maps $\widetilde{u}_h$, we observe that 
        $\int_{B_1} \abs{\nabla \widetilde{u}_h}^{p_h(y)}\,{\d}y \leq \eps_h$, so that, in particular
        \begin{gather}
            \sup_{h \in \N} \int_{B_1} \abs{\nabla \widetilde{u}_h}^{p_h(y)}\,{\d}y \leq 1, \label{eq:properties-u_h-tilde-1} \\
            \forall h \in \N, \quad \mathcal{H}^1\left(J_{\widetilde u_h} \cap B_1\right) \leq \eps_h
            \label{eq:properties-u_h-tilde-2}.
        \end{gather}
        where the second equality follows from~\eqref{eq:decay-2.1} and the obvious fact that 
        $\sigma_h \mathcal{H^1}\left( J_{\widetilde{u}_h} \cap B_1 \right) = \mathcal{H}^1\left( J_u \cap B_{\sigma_h}(x_h) \right)$ for any $h \in \N$ (c.~f. \cite[Remark~7.13]{AFP}).
        
        Next, we set 
		\[
			\gamma_h := \frac{1}{\eps_h}, \qquad t_h := \frac{(\sigma_h \gamma_h)^{1/p_h^0}}{\sigma_h}.
		\]
		We define the rescaled functions
        \[
            v_h : B_1 \to \R^k, \quad v_h := t_h \widetilde{u}_h,
        \]
        and we observe that, for all $h \in \N$,
		\[
			\abs{v_h} = t_h \quad \mbox{a.e. in } B_1,
		\]
		so that $v_h \in \SBV^{p(\cdot)}\left(B_1,\S^{k-1}_{t_h}\right)$.
        From~\eqref{eq:decay-1.2} and the definition of 
		$\gamma_h$, there holds
		\[
			t_h \geq \gamma_h^{2-1/p_h^0},
		\]
		so that (being $p_h^0 \geq p^- > 1$), 
		\[
			t_h \to +\infty \quad \mbox{as} \quad h \to +\infty.
		\]
	\end{step}
	
	\begin{step}[Scaling, Part~2]
		For any $h \in \N$, we define 
		\[
			F_h(v,\gamma_h,B_\sigma) := 
			\int_{B_\sigma} 
			\abs{\nabla v}^{p_h(y)}\,{\d}y 
			+ \gamma_h \mathcal{H}^1(J_v \cap B_\sigma).
		\]
		From the definitions of $p_h$ and $v_h$, Equations~\eqref{eq:decay-2.1}, \eqref{eq:decay-2.2}, 
		\eqref{eq:decay-2.3} become, respectively,
		\begin{align}
			& F_h(v_h,\gamma_h,B_1) \leq 1, \label{eq:F_h} \\
			& \Dev_h(v_h,\gamma_h,B_1,t_h) \leq \theta_h, \label{eq:Dev_h} \\
			& F_h(v_h,\gamma_h,B_\tau) > C_1 \tau^2 F_h(v_h,\gamma_h,B_1) \label{eq:F_h-absurd},
		\end{align}
        where $\Dev_h$ denotes the deviation from miniminality related to $F_h$.
		By~\eqref{eq:properties-u_h-tilde-2}, there holds
		\begin{equation}\label{eq:H1-v_h}
			\mathcal{H}^1(J_{v_h} \cap B_1) \leq \eps_h,
		\end{equation}
		and, since $\gamma_h \sigma_h \leq 1$, we also have
		\begin{equation}\label{eq:modular-v_h}
			\int_{B_1} \abs{\nabla v_h}^{p_h(y)}\,{\d}y = \gamma_h\sigma_h 
			\int_{B_1} \abs{(\nabla u_h)(x_h + \sigma_h y)}^{p_h(y)}\,{\d}y \leq 1.
		\end{equation}
	\end{step}

    \begin{step}[Approximation and lower bound]
        For any $s \in [1/2,1)$, thanks to~\eqref{eq:properties-u_h-tilde-2}, we can use 
        Theorem~\ref{thm:approx-sbv-px} to associate with any $\widetilde{u}_h$ 
        (for any $h$ sufficiently large), a function $\widetilde{z}_h^{(s)}$ belonging to 
        $W^{1,p_h(\cdot)}\left(B_s,\S^{k-1}\right) \cap \SBV^{p_h(\cdot)}\left(B_1,\S^{k-1}\right)$ and a family of balls 
        $\mathcal{F}^{(s)}_h$ satisfying, in particular, 
        \[
            \int_{B_1} \abs{\nabla \widetilde{z}_h^{(s)}}^{p_h(y)}\,{\d}y \lesssim 
            \max\left\{ \norm{\widetilde{u}_h}_{L^{p_h(\cdot)}(B_1)}^{p_h^-}, \norm{\widetilde{u}_h}_{L^{p_h(\cdot)}(B_1)}^{p_h^+} \right\}, 
        \] 
        where the implicit constant at right hand side depends only the log-H\"{o}lder constant of 
        $p(\cdot)$ 
        and, by~\eqref{eq:properties-u_h-tilde-2} and~\eqref{eq:estimate-bad-balls}, 
        \begin{equation}
            \mathcal{L}^2\left( \cup_{\mathcal{F}^{(s)}_h} B \right) 
            \stackrel{\eqref{eq:properties-u_h-tilde-2}, \eqref{eq:estimate-bad-balls}}{\lesssim} \left( \mathcal{H}^1(J_{\widetilde{u}_h} \cap B_1) \right)^2 \leq \eps_h^2 \longrightarrow 0.
        \label{eq:control-measure-bad-balls-sh}
        \end{equation}
        as $h \to +\infty$. Moreover, 
        \begin{equation}\label{eq:same-ws-wt-ae}
            \forall s, t \in [1/2,1), \qquad \widetilde{z}_h^{(s)} = \widetilde{z}_h^{(t)} = \widetilde{u}_h \quad \mbox{a.e. on } B_1 \setminus \left( \cup_{\mathcal{F}^{(s)}_h} B \bigcup \cup_{\mathcal{F}^{(t)}_h} B \right).
        \end{equation}
        
        Next, for any $s \in [1/2,1)$ and any sufficiently large $h \in \N$ as in the above, we denote 
        \[
            \widetilde\lambda_h^{(s)} := \abs{\avg{\widetilde{z}_h^{(s)}}_{B_s}}
        \]
        and we let 
        \begin{equation}\label{eq:def-z_h}
            z_h^{(s)} := t_h \widetilde{z}_h^{(s)}.
        \end{equation}
        We observe that $z_h^{(s)} \in W^{1,p_h(\cdot)}\left(B_s,\S^{k-1}_{t_h}\right) \cap \SBV^{p_h(\cdot)}\left(B_1,\S^{k-1}_{t_h}\right)$ 
        for any $h$ and any $s \in [1/2,1)$ and that each map $z_h^{(s)}$ coincides with the 
        Sobolev approximation of $v_h$ provided by Theorem~\ref{thm:approx-sbv-px}. In particular, 
        \begin{equation}\label{eq:same-zs-zt-ae}
            \forall s, t \in [1/2,1), \qquad z_h^{(s)} = z_h^{(t)} = v_h \quad \mbox{a.e. on } B_1 \setminus \left( \cup_{\mathcal{F}^{(s)}_h} B \bigcup \cup_{\mathcal{F}^{(t)}_h} B \right)
        \end{equation}
        and, by~\eqref{eq:control-measure-bad-balls-sh} and~\eqref{eq:estimate-bad-balls}, for any 
        $s \in [1/2,1)$ we have
        \begin{equation}
            \mathcal{L}^2\left( \cup_{\mathcal{F}^{(s)}_h} B \right) 
            \stackrel{\eqref{eq:control-measure-bad-balls-sh}, \eqref{eq:estimate-bad-balls}}{\lesssim} \left( \mathcal{H}^1(J_{\widetilde{u}_h} \cap B_1) \right)^2 = \left( \mathcal{H}^1(J_{v_h} \cap B_1) \right)^2 \leq \eps_h^2 \longrightarrow 0 
        \label{eq:control-measure-bad-balls}
        \end{equation}
        as $h \to +\infty$.
        Moreover, for any $s \in [1/2,1)$, there holds
        \begin{equation}\label{eq:def-lambda-h-s}
            \avg{z_h^{(s)}}_{B_s} = t_h \avg{\widetilde{z}_h^{(s)}}_{B_s}, \quad \mbox{so that} \quad 
            \lambda_h^{(s)} := \abs{\avg{{z}_h^{(s)}}_{B_s}} = t_h \widetilde\lambda_h^{(s)}.
        \end{equation}
        We claim that, for any $s \in [1/2,1)$,
        \begin{equation}\label{eq:CaLe-trick-conditions}
            \lim_{h \to +\infty} \frac{\lambda_h^{(s)}}{t_h} = 1, \qquad 
            \lim_{h \to +\infty} t_h^{(p^-)^*} \abs{1 - \frac{\lambda_h^{(s)}}{t_h}}^{(p^-)^*} = d,
        \end{equation}
        for some $d \geq 0$.
        Indeed, the second equality entails the first and we have,  
        for any $s \in [1/2,1)$,
        \begin{equation}
        \begin{split}
             t_h^{(p^-)^*} &\abs{1 - \frac{\lambda_h^{(s)}}{t_h} }^{(p^-)^*} = t_h^{(p^-)^*} \abs{ 1 - {\widetilde{\lambda}_h^{(s)}}}^{(p^-)^*} = t_h^{(p^-)^*} \left(\fint_{B_s} \abs{1 - {\widetilde{\lambda}_h^{(s)}}}\,{\d}x \right)^{(p^-)^*} \\
              &= t_h^{(p^-)^*} \left(\fint_{B_s} \abs{\abs{\widetilde{z}_h^{(s)}} - \abs{\avg{\widetilde{z}_h^{(s)}}_{B_s}}} \,{\d}x\right)^{(p^-)^*} 
              \leq t_h^{(p^-)^*} \fint_{B_s} \abs{ \widetilde{z}_h^{(s)} - \avg{\widetilde{z}_h^{(s)}}_{B_s}}^{(p^-)^*}\,{\d}x \\
              &\leq \frac{t_h^{(p^-)^*}}{\pi} \norm{\nabla \widetilde{z}_h^{(s)}}_{L^{p^-}(B_s)}^{(p^-)^*} 
              \stackrel{\eqref{eq:def-z_h}}{\lesssim} \norm{\nabla z_h^{(s)}}_{L^{p^-}(B_s)}^{(p^-)^*} 
              \lesssim \norm{\nabla z_h^{(s)}}_{L^{p(\cdot)}(B_1)}^{(p^-)^*} \stackrel{\eqref{eq:modular-v_h}, \eqref{eq:modular-w}}{\lesssim} 1,
        \end{split}
        \label{eq:CaLe-trick}
        \end{equation}
        where we used that $\abs{\widetilde{z}^{(s)}_h} = 1$ a.e., Jensens's inequality, the 
        classical Sobolev-Poincar\'{e} inequality, and Proposition~\ref{prop:embedding}. 
        By~\eqref{eq:CaLe-trick}, it follows, in particular, that 
        \begin{equation}\label{eq:sobolev-poincare}
            \norm{z_h^{(s)} - \avg{z_h^{(s)}}_{B_s}}_{L^{q}(B_1)} \lesssim 
            \norm{\nabla z_h^{(s)}}_{L^{p^-}(B_s)} \lesssim 1
        \end{equation}
        for any $q \in [1, (p^-)^*]$, up to a constant independent of $s$. 
        Furthermore, we observe that 
        \[
            \abs{\lambda^{(s)}_h - \lambda^{(1/2)}_h} \lesssim 1.
        \]
        Indeed, by~\eqref{eq:def-lambda-h-s} and~\eqref{eq:CaLe-trick-conditions},
        \[
            \abs{\lambda^{(s)}_h - \lambda^{(1/2)}_h} \stackrel{\eqref{eq:def-lambda-h-s}}{=} t_h  
            \abs{\widetilde{\lambda}^{(s)}_h - \widetilde{\lambda}^{(1/2)}_h} \leq 
            t_h\left( \abs{\widetilde{\lambda}^{(s)}_h - 1} + \abs{1-\widetilde{\lambda}^{(1/2)}_h} \right) \stackrel{\eqref{eq:CaLe-trick-conditions}}{\lesssim} 1.
        \]
        Therefore, for any $s \in [1/2,1)$, the sequence 
        $\left\{ z^{(s)}_h - \lambda^{(1/2)}_h \right\}$ is 
        bounded in $W^{1,p^-}\left(B_s, \R^k \right)$ and we can find a subsequence and a 
        limiting map $z^{(s)}$ such that $\left\{ z^{(s)}_h - \lambda^{(1/2)}_h \right\}$ converges to 
        some $z^{(s)}$ weakly in $W^{1,p^-}\left(B_s\right)$, strongly in $L^{p^-}(B_s)$, 
        and pointwise almost everywhere in $B_s$.  By~\eqref{eq:same-zs-zt-ae} 
        and~\eqref{eq:control-measure-bad-balls}, there holds $z^{(s)} = z^{(t)}$ almost everywhere in $B_s$ for all $s$, $t \in [1/2,1)$ with $s \leq t$. Hence, we can define a map 
        $z \in W^{1,p^-}_{\rm loc}\left(B_1,\R^k\right)$ by setting $z := z^{(s)}$ in $B_s$ for any 
        $s \in [1/2,1)$. In particular, we have $\nabla z^{(s)}_h \rightharpoonup \nabla z$ in 
        $L^1(B_s)$. 
        Moreover, the uniform convergence of $\{p_h\}$ to $p_0$ yields, by De~Giorgi's 
        semicontinuity theorem\footnote{The reference 
        \cite{DeGiorgi} is somewhat hard to find. The reader may equally well refer to 
        the more general Ioffe's theorem, e.~g. \cite[Theorem~5.8]{AFP}.} \cite[Corollary~4.6]{DeGiorgi} and 
        since $\nabla z_h^{(s)} = \nabla v_h$ a.e. in $B_1 \setminus \cup_{\mathcal{F}^{(s)}_h} B$,
        \begin{equation}\label{eq:DeGiorgi}
        \begin{split}
            \int_{B_s} \abs{\nabla z}^{p_0}\,{\d} y 
            &\leq \liminf_{h\to+\infty}\int_{B_{s}} \abs{ \left(\nabla z_h^{(s)} \right) \chi_{{B_1} \setminus \cup_{\mathcal{F}_h^{(s)}} B }}^{p_h(y)}\,{\d}y = \liminf_{h\to+\infty}\int_{B_{s} \setminus \cup_{\mathcal{F}_h^{(s)}} B } \abs{ \nabla z_h^{(s)} }^{p_h(y)}\,{\d}y \\
        &= \liminf_{h\to+\infty}\int_{B_{s} \setminus \cup_{\mathcal{F}_h^{(s)}} B } \abs{ \nabla v_h }^{p_h(y)}\,{\d}y  
        \leq \liminf_{h\to+\infty}\int_{B_s} \abs{ \nabla v_h }^{p_h(y)}\,{\d}y  \\
        &\leq \liminf_{h \to +\infty} F_h(v_h, \gamma_h, B_s)
        \end{split}
        \end{equation}
   As~\eqref{eq:DeGiorgi} holds for any $s \in [1/2,1)$, we can let $s \uparrow 1$ in the above 
   inequalities and find both
        \begin{equation}\label{eq:DeGiorgi-bis}
            \int_{B_1} \abs{\nabla z}^{p_0}\,{\d} y \leq \liminf_{h\to+\infty} \int_{B_1} \abs{\nabla v_h}^{p_h(y)}\,{\d}y,            
        \end{equation}
        and
        \begin{equation}\label{eq:LSC}
            \int_{B_1} \abs{\nabla z}^{p_0}\,{\d} y \leq \liminf_{h\to+\infty} F_h(v_h, \gamma_h, B_1).
        \end{equation}
        In addition, Lemma~\ref{lemma:hyperplane} applied to the sequence  
        $\left\{ z^{(s)}_h - \lambda^{(1/2)}_h \right\}$ in $B = B_s$, for any $s \in [1/2,1)$, 
        yields that each map $z^{(s)}$ takes value in a (proper) hyperplane $\Pi^{(s)}$, that we 
        may assume to be contained in the hyperplane $\{ x_k = 0 \}$. 
        Moreover, by Lemma~\ref{lemma:included-hyperplanes}, we can 
        find a maximal hyperplane $\Pi$ containing all the $\Pi^{(s)}$. 
        Finally, since $p^- < 2$ and $p^+ < (p^-)^*$, we have $z \in W^{1,p_0}(B_1, \Pi)$, 
        by~\eqref{eq:DeGiorgi}, Poincar\'{e}'s inequality and the fact that $p_0 \leq p^+$.
    \end{step}

    \begin{step}[Convergence of the energy and $p_0$-minimality of $z$]
    We now improve the lower bound in the previous step by showing that 
    $z$ is a local minimiser of the $p_0$-energy with respect to compactly supported perturbations.

    We start by noticing that, since the function $s \mapsto F_h(v_h, \gamma_h, B_s)$ is increasing and 
    uniformly bounded for $s \in (0,1]$, by Helly's selection theorem and possibly passing to a 
    (not relabelled) subsequence, we may assume that the limit
    \begin{equation}\label{eq:helly}
        \lim_{h \to +\infty} F_h(v_h, \gamma_h, B_s) =: \alpha(s)
    \end{equation}
    exists for all $s \in (0,1]$, where the function $s \mapsto \alpha(s)$ 
    is non-decreasing. (Hence, continuous for all but at most countably many $s \in (0,1]$.) 
    
    Let $v \in W^{1,p_0}(B_1, \Pi)$ be such that 
    $\left\{ v \neq z \right\} \subset \subset B_1$. 
    Let $\{v^\eps\} \subset W^{1,\infty}\left(B_1,\R^k\right)$ be any sequence of Lipschitz functions 
    strongly converging to $v$ in $W^{1,p_0}(B_1)$. Take $\rho', \rho \in (0,1)$, with $\rho' < \rho$ 
    and $\rho$ a continuity point of $\alpha(\cdot)$. For any $h$ large enough, by \eqref{eq:H1-v_h} we 
    have $\mathcal{H}^1\left(v_h \cap B_1 \right) < \eta(1-\rho)$, where $\eta$ is the universal 
    constant in Theorem~\ref{thm:approx-sbv-px}, and therefore the approximation $z^{(\rho)}$ of $v_h$, 
    constructed as in the previous step, is well-defined.
    Let $\zeta \in C_c^\infty(B_\rho, [0,1])$ be 
    a cut-off function with $\zeta \equiv 1$ in $B_{\rho'}$ and 
    $\abs{\nabla \zeta} \leq \frac{2}{\rho-\rho'}$. We may assume that $\left\{ v \neq z \right\} \subset \subset B_\rho$ and define (for $h$ large enough)
    \[
        \overline{w}_h^\eps := \zeta\left( v^\eps + \lambda^{(1/2)}_h \right) + (1-\zeta) z_h^{(\rho)}.
    \]
    By definition, $\overline{w}_h^\eps \in W^{1,p_h(\cdot)}\left(B_\rho, \R^k\right)$, 
    $\overline{w}_h^\eps \in \SBV^{p_h(\cdot)}\left(B_1, \R^k\right)$, and 
    $\overline{w}_h^\eps = z_h^{(\rho)}$ in $B_1 \setminus B_\rho$ (so that, in particular,
    $\overline{w}_h^\eps = v_h$ a.e. in $B_1 \setminus B_{\frac{1+\rho}{2}}$). Moreover, 
    there holds $w_h^\eps = z_h^{(\rho)}$ in the sense of (Sobolev) traces on $\partial B_\rho$.
    By writing
    \[
        \overline{w}_h^\eps = z_h^{(\rho)} + \zeta\left( \left(v^\eps - v \right) + (v - z) + z + \lambda^{(1/2)}_h - z_h^{(\rho)} \right),
    \]
    we see that, for any small enough $\eps > 0$ and any $h \in \N$ large enough,
    \[
        \overline{w}_h^\eps \cdot {\bf e}_k \geq t_h(1-\eps) \quad \mbox{in } B_\rho.
    \]
    Indeed, $v \cdot {\bf e}_k = z \cdot {\bf e}_k = 0$ because $v$ and $z$ take values in $\Pi$, 
    $v^\eps(y) \to v(y)$ at a.e. $y \in B_1$ as $\eps \to 0$ by strong convergence, and 
    $\abs{z + \lambda^{(1/2)}_h - z_h^{(\rho)}} \to 0$ a.e. in $B_\rho$ as $h \to +\infty$, 
    which follows as a consequence of the weak convergence of 
    $\left\{z_h^{(\rho)} - \lambda^{(1/2)}_h \right\}$ 
    to $z$ in $W^{1,p^-}(B_\rho)$ as $h \to +\infty$ and of the Sobolev-Poincar\'e 
    inequality~\eqref{eq:sobolev-poincare}. This also implies that $z_h^{(\rho)} \cdot {\bf e}_k \geq t_h(1-\eps)$ a.~e. for $h$ large enough, which yields the claim. 
    
    As a consequence of the above discussion, for any $\eps >0$ small 
    enough and any $h \in \N$ large enough, we can define
    \[
        w_h^\eps := t_h \frac{\overline{w}_h^\eps}{\abs{\overline{w}_h^\eps}}
    \]
    and we have $w_h^\eps \in W^{1,p_h(\cdot)}\left(B_\rho,\S_{t_h}^{k-1}\right)$, 
    $w_h^\eps \in \SBV^{p_h(\cdot)}\left(B_1, \S_{t_h}^{k-1}\right)$ and $w_h^\eps = z_h^{(\rho)}$ a.e. 
    in $B_1 \setminus B_\rho$. Thus, 
    \[
        \abs{\nabla w_h^\eps} \leq \frac{1}{1-\eps} \abs{\nabla \overline{w}^\eps_h} 
        \quad \mbox{a.e. in } B_1
    \]
    and, in turn,
    \begin{equation}\label{eq:modular-w_h}
    \begin{split}
        &\int_{B_\rho} \abs{\nabla w_h^\eps}^{p_h(y)}\,{\d}y \leq \frac{1}{(1-\eps)^{p^+}} \int_{B_\rho} \abs{\nabla \overline{w}^\eps_h}^{p_h(y)}\,{\d}y \\
        &\leq \frac{1}{(1-\eps)^{p^+}} \left\{ \int_{B_{\rho'}} \abs{\nabla v^\eps}^{p_h(y)} \,{\d}y + \int_{B_\rho \setminus B_{\rho'}} \abs{\nabla z_h^{(\rho)}}^{p_h(y)}\,{\d}y + \int_{B_\rho \setminus B_{\rho'}} \abs{\nabla \overline{w}^\eps_h}^{p_h(y)}\,{\d}y \right\} \\
        &\leq \frac{1}{(1-\eps)^{p^+}} \left\{ \int_{B_{\rho'}} \abs{\nabla v^\eps}^{p_h(y)} \,{\d}y + \int_{B_\rho \setminus B_{\rho'}} \abs{\nabla z_h^{(\rho)}}^{p_h(y)}\,{\d}y \right. \\
        &\left. + 4^{p^+} \left[ \underbrace{\int_{B_\rho \setminus B_{\rho'}} \abs{\nabla v^\eps}^{p_h(y)} \,{\d}y}_{(\text{I}_{h,\eps})}  + \underbrace{\int_{B_{\rho} \setminus B_{\rho'}} \abs{\nabla z_h^{(\rho)}}^{p_h(y)}\,{\d}y }_{(\text{II}_h)}+ \underbrace{\frac{1}{(\rho - \rho')} \int_{B_\rho \setminus B_{\rho'}} \abs{v^\eps + \lambda_h^{(1/2)} - z_h^{(\rho)}}^{p_h(y)}\,{\d}y}_{(\text{III}_h)} \right]  \right\}.
    \end{split}
    \end{equation}
    For any fixed $\eps > 0$, the sequence $\left\{ \abs{\nabla v^\eps}^{p_h(\cdot)} \right\}$ is equibounded in $B_1$ (because $v^\eps$ is Lipschitz and $\{p_h\}$ converges uniformly), and therefore
    \begin{equation}\label{eq:upsc-0}
        \lim_{h \to +\infty}{(\text{I}_{h,\eps})} = \int_{B_{\rho'}} \abs{\nabla v^\eps}^{p_0}\,{\d}y
    \end{equation}
    Consequently, by the strong convergence $v^\eps \to v$ in $W^{1,p_0}(B_1)$,
    \begin{equation}\label{eq:upsc-0-bis}
        \lim_{\eps \to 0} \lim_{h \to +\infty}{(\text{I}_{h,\eps})} = \int_{B_{\rho'}} \abs{\nabla v}^{p_0}\,{\d}y.
    \end{equation}
    Moreover, since $\int_{B_1} \abs{\nabla v_h}^{p_h(y)} \,{\d}y \leq 1$ by~\eqref{eq:modular-v_h}, 
    recalling~\eqref{eq:F_h}, item~\ref{item:approx-iii} of Theorem~\ref{thm:approx-sbv-px} 
    and~\eqref{eq:modular-norm}, we obtain that 
    \begin{equation}\label{eq:upsc-1}
        \int_{B_\rho \setminus B_{\rho'}} \abs{\nabla z_h^{(\rho)}}^{p_h(y)}\,{\d}y 
        \lesssim \left(\int_{B_\rho \setminus B_{\rho'}} \abs{\nabla v_h}^{p_h(y)}\,{\d}y \right)^{\frac{p^-}{p^+}}
        \lesssim \left( F_h(v_h,\gamma_h,B_\rho \setminus B_{\rho'}) \right)^{\frac{p^-}{p^+}}
    \end{equation}
    for any $h \in \N$, where the implicit constants depend only on the log-H\"{o}lder constant of $p$. 
    
    On the other hand, by the very definition of deviation from minimality, there holds
    \begin{equation}\label{eq:upsc-2}
        F_h(v_h,\gamma_h, B_\rho) - \Dev_h(v_h,\gamma_h,B_\rho,t_h) \leq 
        F_h\left(w_h^{(\rho)}, \gamma_h, B_\rho\right) 
    \end{equation}
    for all $h \in \N$. Thus, passing to the limit superior in~\eqref{eq:upsc-2} and 
    using~\eqref{eq:Dev_h} 
    (recalling that, by Lemma~\ref{lemma:CaLe-3.2}, the function $s \mapsto \Dev_h(\cdot, \cdot, \cdot, B_s, \cdot)$ is increasing, for any $h \in \N$), 
    we obtain
    \begin{equation}\label{eq:upsc-3}
        \limsup_{h \to +\infty} F_h(v_h,\gamma_h, B_\rho) \leq 
        \limsup_{h \to +\infty} F_h\left(w_h^{(\rho)}, \gamma_h, B_\rho \right) 
    \end{equation}
    By~\eqref{eq:sobolev-poincare}, it follows that
    \[
        \lim_{h\to +\infty}{(\text{III}_h)} = 0,
    \]
    and by~\eqref{eq:upsc-1} 
    \[
        \limsup_{h \to +\infty} {(\text{II}_h)} =  c\left( \alpha(\rho) - \alpha(\rho') \right),
    \]
    where the constant $c$ is universal.
    Therefore, by~\eqref{eq:modular-w_h}, \eqref{eq:upsc-0}, \eqref{eq:upsc-3}, 
    and the last two estimates, 
    \[
        \limsup_{h \to +\infty} F_h(v_h,\gamma_h, B_\rho) \leq
        \frac{1}{(1-\eps)^{p^+}}\left\{ \int_{B_{\rho'}} \abs{\nabla v^\eps}^{p_0}\,{\d}y + 4^{p^+} \int_{B_\rho \setminus B_{\rho'}} \abs{\nabla v^\eps}^{p_0} \,{\d}y + c(\alpha(\rho)-\alpha(\rho'))\right\},
    \]
    and letting $\rho' \to \rho$ (recalling that $\rho$ is a continuity point of 
    $\alpha(\cdot)$), 
    \[
        \limsup_{h \to +\infty} F_h(v_h,\gamma_h, B_\rho) \leq  \frac{1}{(1-\eps)^{p^+}}\int_{B_\rho} \abs{\nabla v^\eps}^{p_0}\,{\d}y.
    \]
    As this holds for any $\eps > 0$ (small enough), we can let $\eps \to 0$ and 
    thanks to~\eqref{eq:upsc-0-bis} we obtain
    \begin{equation}\label{eq:z-minimality}
        \int_{B_\rho} \abs{\nabla z}^{p_0}\,{\d}y \leq \lim_{h \to +\infty} F_h(v_h,\gamma_h, B_\rho) 
        = \alpha(\rho) \leq \int_{B_\rho} \abs{\nabla v}^{p_0}\,{\d}y,
    \end{equation}
    Since we can take, in particular, $v = z$, we obtain
    \begin{equation}\label{eq:z=alpha}
        \int_{B_\rho} \abs{\nabla z}^{p_0} \,{\d}y = \alpha(\rho)
    \end{equation}
    for all but countably many $\rho \in (0,1)$. In addition, since the left hand side 
    of~\eqref{eq:z=alpha} is a continuous function of $\rho \in (0,1]$, we have that $\alpha(\cdot)$ is actually 
    a continuous function of $\rho \in (0,1]$. Thus, the above argument holds for all $\rho \in (0,1)$ 
    and we conclude that $z$ minimises the $p_0$-energy locally in $B_1$ with respect to compactly 
    supported perturbations.
\end{step}

	\begin{step}[Conclusion]\label{step:decay-conclusion}
        By~\eqref{eq:L-infty-L1-estimate}, we have
        \[
            \sup_{B_\tau} \abs{\nabla z}^p\,{\d}y \leq C_0 \tau^2 \int_{B_1} \abs{\nabla z}^p\,{\d}y
        \]
        and from this, \eqref{eq:helly}, \eqref{eq:z=alpha}, and~\eqref{eq:LSC}, we deduce
        \begin{equation}
        \begin{split}
            \lim_{h \to +\infty} F_h(v_h, \gamma_h, B_\tau) &\stackrel{\eqref{eq:helly},\eqref{eq:z=alpha}}{=} \int_{B_\tau} \abs{\nabla z}^{p_0}\,{\d}y \leq \sup_{B_\tau} \abs{\nabla z}^{p_0} \tau^2 \mathcal{L}^2(B_1)  \\
            & \leq C_0 \tau^2 \int_{B_1} \abs{\nabla z}^{p_0}\,{\d}y \stackrel{\eqref{eq:LSC}}{=} C_0 \tau^2 \lim_{h \to +\infty} F_h(v_h, \gamma_h, B_1) \\
            & < C_1 \tau^2 \lim_{h \to +\infty} F_h(v_h, \gamma_h, B_1)
        \end{split}
        \label{eq:decay-contradiction}
        \end{equation}
        a contradiction to~\eqref{eq:F_h-absurd}.
	\end{step} 
\end{proof}

\begin{remark}\label{rk:stong-log-Hol}
    We remark that, in the proof of Theorem~\ref{thm:decay-lemma}, we used the strongly 
    log-H\"{o}lder condition only to ensure the uniform convergence 
    of the rescaled variable exponents $p_h(\cdot)$ to a constant. 
\end{remark}

We are now in the position to prove Theorem~\ref{thm:density}. 
\begin{proof}[Proof of Theorem~\ref{thm:density}]
    The proof is similar to that of \cite[Theorem~4.7]{LScSoV}, which adapts to 
    the variable exponent setting the classical argument from \cite{DeGiorgiCarrieroLeaci} 
    (see also \cite[Theorem~7.21]{AFP}). 
    We work out the main steps, addressing the reader to the 
    aforementioned references for the missing details.
    \setcounter{step}{0}
    \begin{step}[Small energy in a ball implies power-decay in smaller balls]\label{density-lower-bounds} 
    Fix $\delta > 0$ and $\tau \in (0,1)$ such 
    that $\sqrt{\tau} \leq 1 / C_\delta$ and set $\eps_\delta := \eps(\tau, \delta)$, where $C_\delta$ 
    and $\eps(\tau,\delta)$ are given by Theorem~\ref{thm:decay-lemma}. 
    Fix, in addition, 
    $\sigma \in (0,1)$ such that 
    \[
        \sigma \leq \frac{\eps_\delta}{C_\delta(2 \pi +1 )},
    \]
    and set
    \[
        \rho_\delta := \min\left\{ 1, \eps(\sigma,\delta)^2, \eps_\delta \tau^2 \theta(\tau,\delta), \eps_\delta \sigma \theta(\sigma,\delta) \right\},
    \]
    where $\theta(\tau,\delta)$, $\eps(\sigma,\delta)$ and $\theta(\sigma,\delta)$ are the numbers provided by 
    Theorem~\ref{thm:decay-lemma} corresponding to the pairs $(\tau,\delta)$ and $(\sigma,\delta)$, respectively. Arguing exactly as in \cite[Theorem~4.7]{LScSoV}, from Theorem~\ref{thm:decay-lemma} 
    we obtain that for any $\rho < \rho_\delta':= \frac{\rho_\delta}{\kappa'}$ and any ball 
    $B_\rho \subset\subset \Omega_\delta$ the condition
    \begin{equation}\label{eq:density-compu1}
        F(u,B_\rho(x)) \leq \eps(\sigma,\delta) \rho
    \end{equation}
    implies
    \begin{equation}\label{eq:density-compu2}
        \forall h \in \N, \qquad F\left(u, B_{\sigma \tau^h \rho}(x) \right) \leq \eps_\delta \sigma^\frac{h}{2} \left(\sigma \tau^h \rho \right).    
    \end{equation}
    \end{step}

    \begin{step}[Local lower bounds for the energy in balls centred at jump points]
        Assume that~\eqref{eq:density-compu1} holds for some $x \in \Omega_\delta$, with 
        $B_\rho(x) \subset\subset \Omega_\delta$ and $\rho < \rho_\delta'$. 
        From~\eqref{eq:density-compu2}, we easily obtain that 
        \[
            \lim_{\rho \to 0} \frac{1}{\rho}F(u, B_{\rho}(x)) = 0, 
        \]
        hence Proposition~\ref{prop:not-in-Ju} and Remark~\ref{prop:not-in-Ju} 
        (applied with $\mathcal{M} = \S^{k-1}$) imply that $x \not\in J_u$. Consequently, 
        the lower bound
        \begin{equation}\label{eq:density-compu3}
            F(u, B_\rho(x)) > \eps(\sigma,\delta) \rho
        \end{equation}
        must hold for every $x \in \Omega_\delta \cap J_u$ and actually for every 
        $x \in \Omega_\delta \cap \overline{J_u}$ (because $\eps(\sigma,\delta)$ does not depend on the particular choice of $x \in \Omega_\delta \cap J_u$). Set
        \[
            \theta_\delta := \eps(\sigma,\delta).
        \]
    \end{step}

    \begin{step}[$J_u$ is locally essentially closed]
        We now prove that
        \begin{equation}\label{eq:density-compu6}
            \mathcal{H}^1\left( \Omega_\delta \cap \left( \overline{J_u} \setminus J_u \right) \right) = 0
        \end{equation}
        for any $\delta > 0$.
        To this purpose, it is enough to prove that the difference between $\Omega_\delta \cap \left( \overline{J_u} \setminus J_u \right)$ and a $\mathcal{H}^1$-null set is itself a $\mathcal{H}^1$-null set. To do this, we argue as in the last part of the proof of \cite[Theorem~4.7]{LScSoV} and we consider the set 
        \[
            \Sigma_\delta := \left\{ x \in \Omega_\delta : \limsup_{r \to 0} \int_{B_r} \abs{\nabla u}^{p(y)}\,{\d}y > 0 \right\}.
        \]
        Since $\abs{\nabla u}^{p(\cdot)} \in L^1\left(\Omega\right)$, it follows that 
        $\mathcal{H}^1\left(\Sigma_\delta\right) = 0$ (see, e.g., 
        \cite[Section~2.4.3, Theorem~2.10]{EvansGariepy}). On the other hand, if 
        $x \in \Omega_\delta \cap \left(\overline{J_u} \setminus \Sigma_\delta\right)$, 
        from~\eqref{eq:density-compu3} we obtain
        \[
            \Theta^*(J_u, x) := \limsup_{r \to 0} \frac{\mathcal{H}^1\left(J_u \cap B_r(x)\right)}{2 \pi r} \geq 
            \theta_\delta > 0.
        \]
        Define the Radon measure $\mu := \mathcal{H}^1 \res J_u$ and consider the Borel set 
        $E_\delta := \Omega_\delta \cap \left( \overline{J_u} \setminus J_u \right) \setminus \Sigma_\delta$. Since $\Theta^*(J_u, x) \geq \theta_\delta$ for every $x \in E_\delta$, we deduce that 
        \[
            \mu\left( \left(\Omega_\delta \cap \left(\overline{J_u} \setminus J_u \right) \setminus \Sigma_\delta \right) \right) 
            \geq \theta_\delta \mathcal{H}^1\left(\Omega_\delta \cap \left(\overline{J_u} \setminus J_u \right) \setminus \Sigma_\delta \right). 
        \]
        Since $\mathcal{H}^1\left(\Sigma_\delta\right) = 0$, $\theta_\delta > 0$, and 
        \[
            \mu\left( \left(\Omega_\delta \cap \left(\overline{J_u} \setminus J_u \right)\setminus \Sigma_\delta \right) \right)  = 
            \mathcal{H}^1(\underbrace{J_u \cap (\Omega_\delta \cap \left(\overline{J_u} \setminus J_u \right) \setminus \Sigma_\delta}_{=\varnothing} )) = 0,
        \]
        we obtain~\eqref{eq:density-compu6}. 
    \end{step}

    \begin{step}[Conclusion]
    The essential closedness of $J_u$, i.e., \eqref{eq:closed-jump-set}, 
    now follows by an easy argument by contradiction. Assume 
    $\mathcal{H}^1\left(\Omega \cap \left(\overline{J_u} \setminus J_u\right)\right) > 0$ and 
    let $\delta_h \downarrow 0$ be any decreasing sequence (so that $\left\{\Omega_{\delta_h}\right\}$ 
    is an ascending sequence of sets such that $\Omega = \cup_{h=1}^\infty \Omega_{\delta_h}$). 
    Then, by the $\sigma$-subadditivity of $\mathcal{H}^1$ and~\eqref{eq:density-compu6},
    \[
        0 < \mathcal{H}^1\left(\Omega \cap \left(\overline{J_u} \setminus J_u\right)\right) 
        \leq \sum_{h=1}^\infty \mathcal{H}^1\left(\Omega_{\delta_h} \cap \left(\overline{J_u} \setminus J_u\right)\right) = 0,
    \]
    a contradiction.
    \end{step}
\end{proof}

We are now ready for the proof of Theorem~\ref{thm:B}.

\begin{proof}[Proof of Theorem~\ref{thm:B}]
    Let $u \in \SBV^{p(\cdot)}\left(\Omega, \S^{k-1}\right)$ be any local minimiser 
    of the functional $F(\cdot,\Omega)$ given by~\eqref{eq:functional}. 
    Then, $u$ is a quasi-minimiser of $F(\cdot,\Omega)$, according to Definition~\ref{def:quasi-min}. 
    Consequently, Theorem~\ref{thm:density} tells us that $J_u$ is essentially closed 
    (i.e., \eqref{eq:closed-jump-set} holds for $u$). 
    Let $\Omega' := \Omega \setminus \overline{J_u}$. Then, $\Omega'$ is open and not empty 
    (actually, $\mathcal{H}^1\left(\Omega \setminus \Omega'\right) = 0$, by~\eqref{eq:closed-jump-set} and the rectifiability of $J_u$). 
    By~\eqref{eq:sobolev-sbv} we have 
    $u \in W^{1,p(\cdot)}\left( \Omega', \S^{k-1} \right)$,
    and moreover $u$ is a local minimiser of the functional~\eqref{eq:p(x)-energy} in $\Omega'$. 
    By Theorem~\ref{thm:defilippis}, it follows that $u \in C^{1,\beta_0}_{\rm loc}\left(\Omega_0, \S^{k-1}\right)$, for some $\beta_0 \in (0,1)$ depending only on $k$, $p^-$, $p^+$, $[p]_{0,\alpha}$, $\alpha$ and for a  
    set $\Omega_0$, relatively open in $\Omega'$ (hence, in $\Omega$) and satisfying 
    $\mathcal{H}^{2-p^-}\left( \Omega' \setminus \Omega_0 \right) = 0$ (implying also 
    $\mathcal{H}^1\left(\Omega \setminus \Omega_0\right) = 0$). 
    The conclusion follows.
\end{proof}

\begin{remark}\label{rk:regularity-p-const}
    In the particular case $p$ is constant, we can apply standard elliptic regularity results   
    \cite{Hamburger} to obtain 
    $u \in C^{1,\beta_0}_{\rm loc}\left(\Omega \setminus \overline{J_u}, \S^{k-1}\right)$, for some 
    $\beta_0 \in (0,1)$ depending only on $k$ and $p$.
\end{remark}

\paragraph{Data Availability Statement} Data sharing not applicable to this article as no datasets
were generated or analysed during the current study.

\vskip10pt

\paragraph{Declarations} The authors have no competing interests to declare that are relevant to the
content of this article.

\appendix

\section{Closure and compactness theorems for $\SBV^{p(\cdot)}$}\label{app:cpt}

The following results are certainly known to experts but we have 
not found any explicit proof in literature, hence we provide one 
here, for the reader's convenience.

\begin{theorem}[Closure of $\SBV^{p(\cdot)}$]\label{thm:closure-SBVpx}
	Let $\Omega \subset \R^n$ be a bounded open set, 
	let $p : \Omega \to (1,+\infty)$ be a bounded, log-H\"older continuous 
	variable exponent satisfying
	\[
		1 < p^- \leq p(x) \leq p^+ < +\infty.
	\]
	Assume, in addition, $\theta : [0,+\infty) \to [0,+\infty]$ is an 
	increasing lower semicontinuous functions satisfying 
	$\theta(t)/t \to +\infty$ as $t \to 0$, and that 
	$\{u_h\} \subset \SBV^{p(\cdot)}\left(\Omega,\R^k\right)$ is a sequence satisfying 
	\begin{equation}\label{eq:finite-sup}
		\sup_{h\in \N} \left\{ \int_\Omega \abs{\nabla u_h}^{p(x)} \,\d x 
		+ \int_{J_{u_h}} \theta\left(\abs{u^+_h - u^-_h}\right)\,\d \mathcal{H}^{n-1}(J_{u_h}) \right\} < +\infty.
	\end{equation}
	If $\{u_h\}$ weakly$^*$ converges in $\BV\left(\Omega,\R^k\right)$ to some 
	$u$, then $u \in \SBV^{p(\cdot)}\left(\Omega,\R^k\right)$, the approximate gradients 
	$\nabla u_h$ weakly converge to $\nabla u$ in $L^1\left(\Omega,\R^{k \times n}\right)$, 
	${\D^j}u_h$ weakly$^*$ converge to ${\D^j}u$ in $\Omega$, and 
	\begin{equation}\label{eq:liminf-nabla}
		\int_\Omega \abs{\nabla u}^{p(x)} \,{\d}x \leq 
		\liminf_{h \to \infty} \int_\Omega \abs{\nabla u_h}^{p(x)} \,{\d}x,
	\end{equation}
	and
	\begin{equation}\label{eq:liminf-theta}
		\int_{J_u} \theta\left(\abs{u^+ - u^-}\right) \,{\d}\mathcal{H}^{n-1} 
		\leq 
		\int_{J_{u_h}} \theta\left(\abs{u^+_h - u^-_h}\right) \,{\d}\mathcal{H}^{n-1}
	\end{equation}
	if $\theta$ is concave.
	
\end{theorem}

\begin{proof}
	The conclusion follows by combining the classical closure theorem for  
	$\SBV$ (e.g., \cite[Theorem~4.7]{AFP}) with the lower semicontinuity result 
	given by \cite[Theorem~1.1]{DeCLV} in the context of variable exponents.
	Although short, for clarity we divide the easy proof in two steps.
	
	\setcounter{step}{0}
	\begin{step} 
	Since $\Omega$ is bounded, we have, in 
	particular, $\{u_h\} \subset \SBV^{p^-}\left(\Omega,\R^k\right)$. 
	In addition, by assumption, the sequence $\{u_h\}$ weakly$^*$ converges 
	in $\BV\left(\Omega,\R^k\right)$ to $u \in \BV\left(\Omega,\R^k\right)$. By the classical 
	closure theorem of $\SBV\left(\Omega,\R^k\right)$ (\cite[Theorem~4.7]{AFP}, with 
	$\varphi(t) = t^{p^-}$ and $\theta$ as in the statement), it follows that 
	$u \in \SBV\left(\Omega,\R^k\right)$ (actually, \cite[Theorem~4.7]{AFP} 
	also gives $u \in \SBV^{p^-}\left(\Omega,\R^k\right)$, because of the inequality (4.5) there). 
	By \cite[Theorem~4.7]{AFP} and  \cite[Remark~4.9]{AFP}, 
	\eqref{eq:liminf-theta} follows as well.
	\end{step}
	
	\begin{step} 
	From Step~1, we know that $u_h \to u$ in $L^1\left(\Omega,\R^k\right)$ 
	and that $u \in \SBV\left(\Omega,\R^k\right)$. Therefore, we can apply 
	\cite[Theorem~1.1]{DeCLV} (with $f(x,u,z) = \abs{z}^{p(x)}$, 
	$\Psi \equiv 1$, $a \equiv 0$, and any 	$c > 0$) and deduce 
	that~\eqref{eq:liminf-nabla} holds. But then it follows that 
	$u \in \SBV^{p(\cdot)}\left(\Omega,\R^k\right)$, which is the desired conclusion.
	\end{step}
\end{proof}

\begin{remark}\label{rk:closure-M-valued}
	Assume $\mathcal{M}$ is a smooth, compact Riemannian manifold without boundary, 
	isometrically embedded in $\R^k$. 
	If $\{u_h\} \subset \SBV^{p(\cdot)}(\Omega,\mathcal{M})$ and the 
	assumptions of Theorem~\ref{thm:closure-SBVpx} hold, then the limit  
	function $u$ provided by Theorem~\ref{thm:closure-SBVpx} belongs to
	$\SBV^{p(\cdot)}(\Omega,\mathcal{M})$. Indeed, the strong 
	convergence $u_h \to u$ in $L^1\left(\Omega, \R^k\right)$ implies (up to extraction of 
	a --- not relabelled --- subsequence) $u_h(x) \to u(x)$ for a.e. $x \in \Omega$, 
	hence $u(x) \in \mathcal{M}$ for a.e. $x \in \Omega$.
\end{remark}

\begin{corollary}[Compactness theorem for $\SBV^{p(\cdot)}$]\label{cor:compactness-SBVpx}
	Let $\Omega \subset \R^n$ be a bounded open set 
	and let $p : \Omega \to (1,+\infty)$ be a bounded, log-H\"older continuous 
	variable exponent satisfying
	\[
		1 < p^- \leq p(x) \leq p^+ < +\infty.
	\]
	Let $\theta : [0,+\infty) \to [0,+\infty]$ be an increasing lower 
	semicontinuous function satisfying $\theta(t)/t \to +\infty$ as 
	$t \to 0$.  Finally, assume that 
	$\{u_h\} \subset \SBV^{p(\cdot)}\left(\Omega,\R^k\right)$ is a  
	uniformly bounded sequence in $\BV\left(\Omega,\R^k\right)$ and 
	that~\eqref{eq:finite-sup} holds. Then, we may find 
	$u \in \SBV^{p(\cdot)}\left(\Omega,\R^k\right)$ and extract a (not relabelled) 
	subsequence such that $u_h \to u$ weakly$^*$ in $\BV\left(\Omega,\R^k\right)$. 
	Moreover, the Lebesgue part and the jump part of the derivative converge 
	separately, i.e., $\nabla u_h \to \nabla u$ and 
	${\D^j}u_h \to {\D^j}u$ weakly$^*$ in $\Omega$.
\end{corollary}

\begin{proof}
	Again, we divide the simple proof in two steps.
	
	\setcounter{step}{0}
	\begin{step}
	Assumption~\eqref{eq:finite-sup}, along with $p^- > 1$ 
	and the boundedness of $\Omega$, implies that the sequence $\{\nabla u_h\}$ 
	of the approximate differentials is equiintegrabile (see, e.g.,  
	\cite[Proposition~1.27]{AFP}). On the other hand,
	by assumption, we also have that $\{u_h\}$ is uniformly bounded in 
	$\BV\left(\Omega,\R^k\right)$. Since, by~\eqref{eq:finite-sup}, there holds 
	$\sup_h \int_{J_{u_h}} \theta\left( \abs{u^+_h - u^-_h} \right) \,{\d}\mathcal{H}^{n-1} < +\infty$, 
	the hypotheses of the classical compactness theorem for 
	$\SBV\left(\Omega,\R^k\right)$, as stated for instance in \cite{AlbertiMantegazza} 
	(see Theorem~1.4 there), are satisfied. 
	Therefore, there exist $u \in \SBV\left(\Omega,\R^k\right)$ and a (not relabelled) 
	subsequence such that $u_h \to u$ as $h \to +\infty$, with separate 
	convergence for the Lebesgue part and the jump part of the derivative.
	\end{step}
	
	\begin{step}
	Since $u_h \to u$ weakly$^*$ in $\BV\left(\Omega,\R^k\right)$, 
	Theorem~\ref{thm:closure-SBVpx} yields $u \in \SBV^{p(\cdot)}\left(\Omega,\R^k\right)$. 
	This concludes the proof.
	\end{step}
\end{proof}

\section{Poincar\'{e}'s inequality for bounded variation functions in convex sets}
Let $\Omega \subset \R^n$ be an open set and $k \in \N$ be an integer.
We recall that every function in $\BV\left(\Omega,\R^k\right)$ can be approximated by a sequence smooth 
functions. More precisely, the following characterisation holds 
(see, e.~g.,~\cite[Theorem~3.9]{AFP}): a function 
$u \in L^1\left(\Omega,\R^k\right)$ belongs to $\BV\left(\Omega,\R^k\right)$ if and only if there 
exists a sequence $\{u_h\}\subset \left(W^{1,1}\cap C^\infty\right)\left(\Omega,\R^k\right)$ 
converging to $u$ in $L^1\left(\Omega,\R^k\right)$ and such that 
\begin{equation}\label{eq:approx-BV-smooth}
	L := \lim_{h\to\infty} \int_\Omega \abs{\nabla u_h}\,{\d}x \leq 
	\abs{{\D}u}(\Omega).
\end{equation}
By this approximation result and Poincar\'{e}'s inequality for Sobolev 
functions in convex sets 
(e.g.~\cite[Theorem~12.30]{Leoni}), we derive the 
following more precise version of the classical Poincar\'{e}'s inequality in $\BV$ 
(see, e.g., \cite[Theorem~3.44]{AFP}), showing up the explicit 
dependence on the diameter of $\Omega$. Since we have not found an explicit proof of this statement 
in the literature, we provide a detailed one here, for reader's convenience.

\begin{prop}\label{prop:poincare-BV-convex}
	Let $\Omega \subset \R^n$ be a bounded, convex open set. There exists 
	a constant $C=C(n,k)$, depending only on $n$ and $k$, such that for any 
	$u \in \BV(\Omega,\R^k)$ there holds
	\begin{equation}\label{eq:poincare-BV-convex}
		\norm{u-\avg{u}_\Omega}_{L^1(\Omega,\R^k)} \leq C(n,k) (\diam(\Omega)) 
		\abs{{\D}u}(\Omega)
	\end{equation}
\end{prop}

\begin{proof}
	Since $u$ belongs to $\BV\left(\Omega,\R^k\right)$, by \cite[Theorem~3.9]{AFP} there 
	exists $\{u_h\} \subset \left(W^{1,1}\cap C^\infty\right)\left(\Omega,\R^k\right)$, 
    a sequence of smooth functions in $\Omega$, such that $u_h \to u$ in 
    $L^1\left(\Omega,\R^k\right)$ as 
	$h \to +\infty$ and such that~\eqref{eq:approx-BV-smooth} holds. To the 
	functions $u_h$, we may apply {Poincar\'{e}}'s inequality for Sobolev 
	functions in convex sets (e.g., \cite[Theorem~12.30]{Leoni}), from which it 
	follows that
	\begin{equation}\label{eq:poincare-BV-convex-compu1}
		\norm{u_h - \avg{u_h}_\Omega}_{L^1\left(\Omega\right)} \leq C(n,k)
		(\diam\Omega)\norm{\nabla u_h}_{L^1(\Omega)},
	\end{equation}
	where the constant $C(n,k)$ depends only on $n$ and $k$. On the other hand, 
	by~\eqref{eq:approx-BV-smooth}, up to discarding finitely many terms of 
	the sequence $\{u_h\}$, we may assume that
	\begin{equation}\label{eq:poincare-BV-convex-compu2}		
	\sup_{h \in \N} \norm{\nabla u_h}_{L^1(\Omega)} \leq 
		\abs{{\D}u}(\Omega).
	\end{equation} 
	Combining~\eqref{eq:poincare-BV-convex-compu1} and  
	\eqref{eq:poincare-BV-convex-compu2} with the fact that the 
	$L^1$-convergence of $\{u_h\}$ to $u$ obviously implies that 
	\[
		\norm{u_h - \avg{u_h}_\Omega}_{L^1(\Omega)} \to 
		\norm{u - \avg{u}_\Omega}_{L^1(\Omega)}
	\]  
	as $h \to +\infty$, the conclusion follows.
\end{proof}

\section{A counterexample in higher dimension}\label{app:counterex}

In this appendix, expanding on Remark~\ref{rk:obstruction}, we provide an example which shows that the 
construction of the Sobolev approximation in Section~\ref{sec:approx} cannot work in higher dimensions 
under the mere assumption of smallness of the jump set.

Carefully looking at the proof of \cite[Theorem~2.1]{CFI-ARMA}, it is seen that the only point in which 
the assumption on the dimension is used is in the proof of property $({\rm P}_2)$, and more precisely 
in the passage from \cite[(2.8)]{CFI-ARMA} to \cite[(2.9)]{CFI-ARMA}, 
which relies on inequality~\eqref{eq:(2.3)-CFI-ARMA} (i.e., on \cite[(2.3)]{CFI-ARMA}). 
Again by inspection of the proof of $({\rm P}_2)$, it is easily realised that, in order for 
the construction to work in higher dimensions, for any $\eta \in (0,1)$ and any given Borel 
set $J \subset B_{2r}^n$ with $\mathcal{H}^{n-1}(J) < \eta (2r)^{n-1}$ 
we must be able to find a radius $R \in (r,2r)$ satisfying \emph{both}
\begin{equation}\label{eq:(2.2)-CFI-ARMA-higher-dim}
    \mathcal{H}^{n-1}\left( J \cap \partial B_R^n \right) = 0
\end{equation}
and
\begin{equation}\label{eq:(2.3)-CFI-ARMA-higher-dim}
    \mathcal{H}^{n-1}\left( J \cap \left( B_R^n \setminus B_{R -\delta_h}^n\right) \right) < \eta C_n \delta_h^{n-1}, \quad \mbox{for every } h \in \N,
\end{equation}
where $\delta_h := R 2^{-h}$ and $C_n$ is a dimensional constant.

Now, for every $n \in \N$, given any Borel set $J$ in $B_{2r}^n$ with $\mathcal{H}^{n-1}(J) < +\infty$, 
the measure $\mathcal{H}^{n-1} \res J$ is a Radon measure, hence 
equation~\eqref{eq:(2.2)-CFI-ARMA-higher-dim} holds for all $R \in (r,2r)$ except at most countably 
(indeed, Radon measures can charge at most countably many boundaries of encapsulated sets; 
c.f., e.g., the discussion in \cite[Example~1.63]{AFP}). 
However, \eqref{eq:(2.3)-CFI-ARMA-higher-dim} may fail for $n \geq 3$.
We exhibit an example below for $n = 3$ (for ease of notation) that can be easily adapted to any dimension.

\begin{prop}\label{prop:counterex}
    For any $\eps \in (0,1)$, any $C > 0$, and any $r > 0$, there exists a 
    $\mathcal{H}^2$-rectifiable, star-shaped, relatively closed and essentially closed Borel set 
    $J \subset B^3_{2r}$ with the following properties:
    \begin{enumerate}[(i)]
        \item $\mathcal{H}^{n-1}(J) < \eps (2r)^{n-1}$. 
        \item There exists $h_0 \in \N$ so that 
    \begin{equation}\label{eq:no-(2.3)-higher-dim}
        \mathcal{H}^{2}\left( J \cap \left(B_R \setminus B_{R\left(1 - 2^{-h_0}\right)} \right) \right) 
        \geq C \eps \delta_{h_0}^2 
    \end{equation}
    for every $R \in (r,2r)$.
    \end{enumerate}
\end{prop}

\begin{proof}
    Fix arbitrarily $\eps \in (0,1)$ and $C > 0$.
    By scaling, it suffices to consider the case $2r = 1$. Let us set, for brevity, 
    $\mathbb{B} := B_1^3$. Let $\{x_j\}$ be an enumeration of $\mathbb{Q}^3 \cap \partial \mathbb{B}$. 
    For any $j \in \N$, let $\mathcal{C}_j^{(\eps)}$ be the cone with apex the origin, axis the radius 
    $\overrightarrow{O x_j}$, and opening angle $\frac{\eps 2^{-j}}{40\pi}$ (the reason for this 
    choice will be clear after~\eqref{eq:J-compu-2}).
    For any $j \in \N$, let 
    $\partial \mathcal{B}^2_j := \mathcal{C}_j^{(\eps)} \cap \S^2_{1/2}$ denote the boundary 
    of the geodesic ball cut on $\mathbb{S}^2_{1/2} = \partial B^3_{1/2}$ by the solid cone bounded by 
    $\mathcal{C}_j^{(\eps)}$. Correspondingly, the same solid cone cuts geodesic balls 
    $\widehat{\mathcal{B}}_j^2$ on $\S^2$, whose radius is twice that of $\mathcal{B}^2_j$, 
    for each $j \in \N$.
    
    By Vitali's covering lemma, there exists a subsequence 
    $\left\{x_{j_l}\right\}$ such that the balls $\mathcal{B}^2_{j_l}$ have disjoint 
    closures and such that the geodesic balls $5 \mathcal{B}^2_{j_l}$, concentric with the balls 
    $\mathcal{B}^2_{j_l}$ but with radius $5$ times larger, are a covering of 
    $\cup_{h \in \N} \mathcal{B}_j^2$. 
    We notice that, by the choice of the open angle of the cones,
    \begin{equation}\label{eq:J-compu-1}
        \sum_{l \in \N} \mathcal{H}^1\left(\partial \mathcal{B}^2_{j_l}\right) \leq 
        \sum_{l \in \N} 2 \pi \left( \frac{2^{-h_l}}{40\pi} \eps  \right) < \eps/10.
    \end{equation} 
    Upon setting 
    \[
        \overline{J} := \overline{\mathbb{B}} \cap \bigcup_{l \in \N} \mathcal{C}_{j_l}^{(\eps)},
    \]
    we have that $J$ is a countably $\mathcal{H}^2$-rectifiable, star-shaped (because all the cones 
    $\mathcal{C}^{(\eps)}_j$ intersect at the origin), closed set in 
    $\overline{\mathbb{B}}$.  
    Moreover, since for each $l \in \N$ the geodesic ball $\widehat{\mathcal{B}}^2_{j_l}$ has radius 
    twice that of $\mathcal{B}^2_{j_l}$ and since the slant height of each (truncated) cone 
    $\overline{\mathbb{B}} \cap \mathcal{C}^{(\eps)}_j$ is one, by~\eqref{eq:J-compu-1} there holds
    \begin{equation}\label{eq:J-compu-2}
        \mathcal{H}^2\left(\overline{J}\right) = \sum_{l \in \N} \mathcal{H}^1\left(\partial \widehat{\mathcal{B}}^2_{j_l}\right) \leq 10 \sum_{l \in \N} \mathcal{H}^1\left(\partial \mathcal{B}^2_{j_l}\right) < \eps,
    \end{equation}
    which gives both item~{(i)} and the $\mathcal{H}^2$-rectifiability of $\overline{J}$. Then, 
    $J := \mathbb{B} \cap \overline{J}$ is a Borel set in $\mathbb{B}$ (because it is relatively closed in $\mathbb{B}$); moreover, $J$ is $\mathcal{H}^2$-rectifiable, star-shaped, essentially closed 
    (because $J = \overline{J} \setminus \left(\overline{J} \cap \partial \mathbb{B}\right)$, hence 
    $\mathcal{H}^2\left(\overline{J} \setminus J\right) = 0$), 
    and it satisfies $\mathcal{H}^2(J) < \eps$.

    It remains to prove~\eqref{eq:no-(2.3)-higher-dim}. To this purpose, 
    let $\kappa := \sum_{l \in \N} 2^{-j_l}$ and  
    notice that for any choice 
    of $R \in (1/2,1)$ and any $h \in \N$ we have
    \[
    \begin{split}
        \mathcal{H}^2\left( J \cap \left( B_{R} \setminus B_{R - \delta_h} \right) \right) &\geq 
        \delta_h \sum_{l \in \N} 2\pi \left(R-\delta_h\right) \sin\left(\frac{\eps 2^{-j_l}}{40 \pi} \right) \\ 
        &\geq \frac{1}{40} \delta_h \left(R-\delta_h\right) \eps \sum_{l \in \N} 2^{-j_l} 
        = \frac{\kappa}{40} \delta_h \left( R - \delta_h \right) \eps,
    \end{split}
    \]
    where $\delta_h := R 2^{-h}$. From here, \eqref{eq:no-(2.3)-higher-dim} follows as soon as 
    $h_0$ is chosen so that 
    \[
        \frac{\kappa}{40} (R-\delta_{h_0}) \delta_{h_0} \eps 
        = \frac{\kappa}{40} R^2 \left(1-2^{-h_0}\right) 2^{-h_0} \eps
        \geq C R^2 2^{-2h_0} \eps,
    \]
    i.e., for $h_0 \geq \log_2\left(1+\frac{40 C}{\kappa}\right)$.
\end{proof}

\begin{remark}
    Let $\eps \in (0,1)$, $r > 0$, $C > 0$ be arbitrary, let  
    $J \subset B^3_{2r}$ be the corresponding Borel set provided 
    by~Proposition~\ref{prop:counterex}, and let $\Omega$ be the union of the solid cones determined by the sets $\mathcal{C}^{(\eps)}_{j_l}$ defined 
    as in the proof of Proposition~\ref{prop:counterex}, intersected with $\mathbb{B}$. Define 
    $u : B^3_{2r} \to \R^k$ by setting 
    \[
        u := 
        \begin{cases}
            {\bf e}_1 & \mbox{in } \Omega, \\
            0 & \mbox{in } B_1 \setminus \overline{\Omega}.
        \end{cases}
    \]
    Then, $u \in \SBV^p\left(B_{2r}^3, \R^k\right)$ for any $p \geq 1$, $J_u = J$, 
    $\mathcal{H}^2(J \setminus J_u) = 0$, $\mathcal{H}^2\left(J_u\right) < \eps$, 
    $J_u$ is essentially closed and it does not 
    satisfy~\eqref{eq:(2.3)-CFI-ARMA-higher-dim}.  
\end{remark}

\section{A criterion for being out of the jump set}\label{app:not-in-Ju}

In this appendix, we prove a sufficient criterion that allows for excluding that a point belongs 
to the jump set. 
An analogous, but slightly stronger, statement concerning the whole approximate discontinuity set $S_u$ 
is classical (see \cite{DeGiorgiCarrieroLeaci} and \cite[Theorem~7.8]{AFP}) but the classical proof 
uses medians and truncations, which we want to avoid. 
Here, we use only the tools provided by the Sobolev approximation results from Section~\ref{sec:approx}.
\begin{prop}\label{prop:not-in-Ju}
    Let $\Omega \subset \R^2$ be a bounded open set and $p : \Omega \to (1,+\infty)$ be a 
    log-H\"{o}lder continuous, variable exponent satisfying $p^- > 1$ and $p^+ < 2$.
    Let $u \in \left(L^\infty \cap \SBV^{p(\cdot)} \right)\left(\Omega, \R^k\right)$ and 
    $x_0 \in \Omega$. If
    \begin{equation}\label{eq:vanishing-rescaled-energy}
        \lim_{\rho \downarrow 0} \frac{1}{\rho}\left[ \int_{B_\rho(x_0)} \abs{\nabla u}^{p(x)}\,{\d}x + \mathcal{H}^1\left(J_u \cap B_\rho(x)\right) \right] = 0,
    \end{equation}
    then $x_0 \not\in J_u$.
\end{prop}

\begin{proof} We argue by contradiction, exploiting the tools provided by the Sobolev approximation.
    \setcounter{step}{0}
    \begin{step}
    Suppose, for the sake of a contradiction, that $x \in J_u$. Then, by definition 
    (c.f., \cite[Definition~3.67]{AFP}), there exist $a$, $b \in \R^k$ 
    and $\nu_0 \in \S^1$ such that $a \neq b$ and 
    \[
        \lim_{\rho \downarrow 0} \fint_{B_\rho^+(x_0,\nu_0)} \abs{u(x)-a}\,{\d}x = 0, \qquad 
        \lim_{\rho \downarrow 0} \fint_{B_\rho^+(x_0,\nu_0)} \abs{u(x)-b}\,{\d}x = 0.
    \]
    Moreover, the triple $(a,b,\nu_0)$ is uniquely determined, up to a permutation and a change of 
    sign.
    
    Observe that, thanks to the trivial pointwise inequality
    \[
        \forall \xi \in \R^k, \qquad \abs{\xi}^{p^-} \leq 1 + \abs{\xi}^{p(x)}, 
    \]
    holding at each point $x \in \Omega$, condition~\eqref{eq:vanishing-rescaled-energy} implies
    \begin{equation}\label{eq:vanishing-rescaled-energy-bis}
        \lim_{\rho \downarrow 0} \frac{1}{\rho}\left[ \int_{B_\rho(x_0)} \abs{\nabla u}^{p^-}\,{\d}x + \mathcal{H}^1\left(J_u \cap B_\rho(x) \right) \right] = 0,
    \end{equation}
    \end{step}

    \begin{step}\label{step:not-Ju-step-2}
    We are going to show that, given~\eqref{eq:vanishing-rescaled-energy-bis}, we can 
    find a sequence $\rho_j \downarrow 0$ as $j \to +\infty$ and $m \in \R^k$ such that 
    \begin{equation}\label{eq:not-Ju-1}
        \lim_{j \to +\infty}\fint_{B_{\rho_j}(x_0)} \abs{u - m}\,{\d}x = 0.
    \end{equation}
    To this purpose, define the blow-up $u_\rho$ of $u$ at $x_0$ by letting 
    $u_\rho(y) := \rho^{(1-p^-)/p^-} u(x_0 + \rho y)$ for all $y \in B_1$ and each $\rho > 0$.  
    Then, condition~\eqref{eq:vanishing-rescaled-energy-bis} can be recast as
    \begin{equation}\label{eq:not-Ju-2}
        \lim_{\rho \downarrow 0} \left[ \int_{B_1} \abs{\nabla u_\rho}^{p^-}\,{\d}y + \mathcal{H}^1\left(J_{u_\rho} \cap B_1 \right) \right] = 0.
    \end{equation}
    Now, take any sequence $\{\rho_h\}$ such that $\rho_h \downarrow 0$ as $h \to +\infty$.
    Then, for any $h \in \N$,
    \[
        \int_{B_1} \abs{\nabla u_{\rho_h}}^{p^-}\,{\d}y + \mathcal{H}^1\left(J_{u_{\rho_h}} \cap B_1 \right) = \sigma_h,
    \]
    where $\sigma_h \to 0$ as $h \to + \infty$. For each $h \in \N$ so large that 
    $2\sigma_h < \eta$, where $\eta$ is the universal constant in 
    Proposition~\ref{thm:approx-sbv-px-Rk}, 
    let $s_h := 1 - \frac{2\sigma_h}{\eta}$. Then, let $w^{(s_h)}_{\rho_h} \in  W^{1,p^-}\left(B_{s_h}, \R^k\right) \cap \SBV^{p^-}\left(B_1, \R^k\right)$ be the approximation of $u_{\rho_h}$ provided by Theorem~\ref{thm:approx-sbv-px-Rk}. Recall that, since $u$ (and hence each map $u_\rho$) is 
    bounded, each map $w^{(s_h)}_{\rho_h}$ is bounded as well and, indeed, $\norm{w^{(s_h)}_{\rho_h}}_{L^\infty(B_1)} \leq \norm{u}_{L^\infty(B_1)}$ for any $h$. Set
    \[
        m^{(s_h)}_{\rho_h} := \fint_{B_{s_h}} w^{(s_h)}_{\rho_h}\,{\d}y.
    \]
    Since $\norm{w^{(s_h)}_{\rho_h}}_{L^\infty(B_1)}$ is uniformly bounded, 
    the sequence $\left\{ m^{(s_h)}_{\rho_h} \right\}_h$ is 
    bounded in $\R^k$. Hence, we can find a sequence $\{h_j\}_j$ and $m \in \R^k$ so that 
    \[
        m_j := m^{(s_{h_j})}_{\rho_{h_j}} \to m \qquad \mbox{as } j \to +\infty.
    \]
    Let us also set, for brevity's sake,
    \[
       \rho_j := \rho_{h_j}, \qquad s_j := s_{h_j}, \qquad u_j := u_{\rho_{h_j}}, \qquad w_j := w^{(s_{h_j})}_{\rho_{h_j}}.
    \]
    for each $j \in \N$.
    \end{step}

    \begin{step}
        We claim that~\eqref{eq:not-Ju-1} holds for $\{\rho_j\}$ and $m$ as in 
        Step~\ref{step:not-Ju-step-2}. Indeed, \eqref{eq:not-Ju-1} is equivalent to 
        \begin{equation}\label{eq:not-Ju-3}
            \lim_{j \to +\infty}\int_{B_1} \abs{u_j - m}\,{\d}y = 0,
        \end{equation}
        and, for any $j \in \N$, we have
        \[
        \begin{split}
            \int_{B_1} &\abs{u_j - m}\,{\d}y \leq \int_{B_1} \abs{u_j - w_j}\,{\d}y +  \int_{B_1} \abs{w_j - m}\,{\d}y \\
            &\leq \underbrace{\int_{B_1} \abs{u_j - w_j}\,{\d}y}_{:= ({\rm I})_j} +  \underbrace{\int_{B_1} \abs{w_j - m_j}\,{\d}y}_{:= ({\rm II})_j} + 
            \abs{B_1} \abs{m_j - m}
        \end{split}
        \]
        By construction, the last term at right hand side above tends to 0 as $j \to +\infty$. 
        As for $({\rm I})_j$, 
        by items~\ref{Rk-approx-item-i} and~\ref{Rk-approx-item-iii} of 
        Proposition~\ref{thm:approx-sbv-px-Rk} we have
        \[
            \int_{B_1} \abs{u_j - w_j}\,{\d}y \leq 2 \norm{u}_{L^\infty(\Omega)} \abs{ \left\{ u_j \neq w_j \right\} } \lesssim \left( \mathcal{H}^1\left(J_{u_j} \cap B_1\right) \right) = {\rm o}(1)
        \]
        as $j \to +\infty$. We now estimate $({\rm II}_j)$, using that $w_j$ is a Sobolev function in 
        $B_{s_j}$ and item~\ref{Rk-approx-item-v} of Proposition~\ref{thm:approx-sbv-px-Rk}:
        \[
        \begin{split}
            \int_{B_1} \abs{w_j - m_j}\,{\d}y &= \int_{B_{s_j}} \abs{w_j - m_j}\,{\d}y + \int_{B_1 \setminus B_{s_j}} \abs{w_j - m_j}\,{\d}y \\
            &\leq \int_{B_{s_j}} \abs{\nabla w_j}\,{\d}y + 2 \norm{u}_{L^\infty(\Omega)}  \pi \left(1 - s_j^2 \right)\\
            &\lesssim \abs{B_{s_j}}^{1-\frac{1}{p^-}} \left( \int_{B_{s_j}} \abs{\nabla u_j}^{p^-} \,{\d}y \right)^{\frac{1}{p^-}} + {\rm o}(1) \\
            &= {\rm o}(1)
        \end{split}
        \]
        as $j \to +\infty$, because $s_j \to 1$ as $j \to +\infty$ and thanks to item~\ref{Rk-approx-item-v} of Proposition~\ref{thm:approx-sbv-px-Rk} and to~\eqref{eq:not-Ju-2}.
        By the above estimates, ~\eqref{eq:not-Ju-3} follows immediately.
    \end{step}

    \begin{step}[Conclusion]
        By assumption, we should have
         \[
        \lim_{j \to +\infty} \fint_{B_{\rho_j}^+(x_0,\nu_0)} \abs{u(x)-a}\,{\d}x = 0, \qquad 
        \lim_{j \to +\infty} \fint_{B_{\rho_j}^+(x_0,\nu_0)} \abs{u(x)-b}\,{\d}x = 0,
        \]
        with $a \neq b$ and some $\nu_0 \in \S^1$, while~\eqref{eq:not-Ju-3} implies
        \[
        \lim_{j \to +\infty} \fint_{B_{\rho_j}^+(x_0,\nu)} \abs{u(x)-m}\,{\d}x = 0, \qquad 
        \lim_{j \to +\infty} \fint_{B_{\rho_j}^+(x_0,\nu)} \abs{u(x)-m}\,{\d}x = 0,
        \]
        for any $\nu \in \S^1$, 
        and this a contradiction because the triple $(a,b,\nu_0)$ is uniquely determined, up to a 
        permutation and a change of sign.
    \end{step}
\end{proof}

\begin{remark}\label{rk:not-in-Ju}
    In particular, the assumptions of Proposition~\ref{prop:not-in-Ju} are satisfied if 
    $u$ takes values into a compact Riemannian manifold $\mathcal{M}$ and 
    satisfies~\eqref{eq:vanishing-rescaled-energy}. 
\end{remark}

\begin{remark}
    A slightly stronger version of Proposition~\ref{prop:not-in-Ju} is obtained by replacing the 
    assumption $u \in L^\infty\left(\Omega, \R^k\right)$ with the weaker assumption
    \[
        \lim_{\rho \downarrow 0} \fint_{B_\rho(x_0)} \abs{u}^q\,{\d}x < +\infty
    \]
    for some $q > 1$, which is exactly the hypothesis required in the original formulation of the 
    criterion from \cite{DeGiorgiCarrieroLeaci} (in place of the boundedness of $u$).
    Up to the fact that we also need to assume~\eqref{eq:p+p-star} 
    (in order to apply Proposition~\ref{thm:approx-sbv-px-Rk}), 
    the details are exactly as in \cite{DeGiorgiCarrieroLeaci} or, equivalently, \cite[Theorem~7.8]{AFP} and, therefore, omitted.
\end{remark}
\section{Compactness of sequences of $\SBV^{p(\cdot)}$-functions with vanishing jump set}
Here,
arguing along the lines of \cite[Section~3]{CFI-AIHP}, we use 
Theorem~\ref{thm:approx-sbv-px} and compactness results in variable 
exponent Sobolev spaces (see, e.g., \cite[Chapter~8]{DHHR})
to prove a compactness and lower semicontinuity result for sequences 
with vanishing jump set. Such a result is the counterpart of 
\cite[Proposition~3.2]{CFI-AIHP} and it can be applied to regularise sequences of functions in 
$\SBV^{p(x)}(B_\rho,\mathcal{M})$ with vanishing jump set. 
\begin{prop}\label{prop:CFI-AIHP-3.2-type}
	Let $p : \R^2 \to (1,+\infty)$ be a 
	variable exponent satisfying~{(${\rm p}_1$)} and~{(${\rm p}_2$)}. 
	Let $k \in \N$ and let $(\mathcal{M},g)$ be either 
	$\R^k$ (endowed with the canonical metric) or a compact, connected, smooth Riemannian manifold 
	without boundary, isometrically embedded in $\R^k$.  Assume 
	$\{u_h\} \subset \SBV^{p(x)}\left(B_\rho,\mathcal{M}\right)$ is a sequence satisfying 
	\begin{equation}\label{eq:equibdd}
		\sup_{h\in \N} \int_{B_\rho} \abs{\nabla u_h}^{p(x)}\,{\d}x \leq \Lambda_* < +\infty, 
		\qquad 
		\mathcal{H}^1\left(J_{u_h}\right) \to 0 \quad \mbox{as } h \to +\infty.
	\end{equation}
	Then, we can find a (not relabelled subsequence) and a function 
    $u \in W^{1,p(\cdot)}(B_\rho,\mathcal{M})$ such that 
    $u_h \to u$ in $L^{p(\cdot)}(B_\rho)$ as $h \to +\infty$. Moreover,
    \begin{equation}\label{eq:bound-u}
        \int_{B_\rho} \abs{\nabla u}^{p(x)}\,{\d}x \leq \liminf_{h \to +\infty}
        \int_{B_\rho} \abs{\nabla u_h}^{p(x)}\,{\d}x.
    \end{equation}
\end{prop}

\begin{proof}
	The proof is along the lines of \cite[Proposition~3.2]{CFI-AIHP}. 
	We reproduce the main steps, to point out the relevant changes.
	
	For any $s \in [1/2,1)$ and any $h \in \N$ so large that 
    $\mathcal{H}^1\left( u_h \right) < \eta(1-s)$, where $\eta$ is the universal constant 
    provided in Theorem~\ref{thm:approx-sbv-px}, we let 
	$w^{(s)}_h \in \SBV^{p(\cdot)}\left(B_\rho,\mathcal{M}\right) \cap W^{1,p(\cdot)}\left(B_{s \rho},\mathcal{M}\right)$ and $\mathcal{F}^{(s)}_h$ 
	be the function and the family of balls obtained by 
	Theorem~\ref{thm:approx-sbv-px}. By item~\ref{item:approx-iii} in 
	Theorem~\ref{thm:approx-sbv-px}, there holds
	\begin{equation}\label{eq:comp1}
		\int_{B_{s\rho}} \abs{\nabla w^{(s)}_h}^{p(x)}\,{\d}x \lesssim 
			\max\left\{ \norm{\nabla u_{h}}_{L^{p(\cdot)}(B_\rho)}^{p^-}, \norm{\nabla u_{h}}_{L^{p(\cdot)}(B_\rho)}^{p^+}\right\}
	\end{equation}
	where the implicit constant at right hand side depends only on the quantities listed in 
    Remark~\ref{rk:projection}. 
	Moreover, by Poincar\'{e}'s inequality (\cite[Theorem~8.2.4]{DHHR}) and~\eqref{eq:comp1}, 
	there holds 
	\begin{equation}\label{eq:comp2}
		\norm{w^{(s)}_h - m^{(s)}_h}_{L^{p(\cdot)}(B_{s\rho})} \lesssim \left(1+\rho^2\right) 
		\max\left\{ \norm{\nabla u_{h}}_{L^{p(\cdot)}(B_\rho)}, \norm{\nabla u_{h}}_{L^{p(\cdot)}(B_\rho)}^{\frac{p^-}{p^+}}, \norm{\nabla u_{h}}_{L^{p(\cdot)}(B_\rho)}^{\frac{p^+}{p^-}}\right\}
	\end{equation}
	where $m^{(s)}_h := \left\langle w^{(s)}_h \right\rangle_{B_{s \rho}}$ is 
	the average of $w^{(s)}_h$ on $B_{s \rho}$ and the implicit constant at right 
    hand side depends only on the quantities listed in Remark~\ref{rk:projection}.
    For any $s \in [1/2,1)$, by the compactness of $\mathcal{M}$, the sequence 
    $\left\{ m^{(s)}_h\right\}$ is obviously bounded and hence can extract from it 
    a subsequence which converges to some $\overline{m}^{(s)}$ in $\R^k$.
	%
    Also, again by the compactness of $\mathcal{M}$, we get
	\[
		\norm{m^{(s)}_h - m^{(1/2)}_h}_{L^{p(\cdot)}(B_{\rho})} \lesssim 1,
	\]
    up to a constant which depends only on $\mathcal{M}$.
	Therefore, by triangle inequality and~\eqref{eq:comp2}, the sequence 
	$\left\{w^{(s)}_h - m^{(1/2)}_h\right\}$ is bounded 
	in $W^{1,p(\cdot)}\left(B_{s \rho}, \R^k\right)$ for any $s \in [1/2,1)$  
    and we can find a subsequence (depending 
	on $s$ and not relabelled) such that it converges to some $w^{(s)}$ weakly 
	in $W^{1,p(\cdot)}\left(B_{s\rho}, \R^k\right)$, strongly in 
	$L^{p(\cdot)}\left(B_{s\rho},\R^k\right)$ 
	(by \cite[Corollary~8.3.2]{DHHR}), and pointwise a.e. on $B_{s\rho}$ 
	(because, in particular, we have weak convergence in 
	$W^{1,p^-}\left(B_{s \rho}, \R^k\right)$ and hence strong convergence in $
	L^{p^-}\left(B_{s \rho}, \R^k\right)$ by the classical Rellich-Kondrachov theorem).
	
	By item~\ref{item:approx-iv} in Theorem~\ref{thm:approx-sbv-px}, for any $s \in [1/2, 1)$ 
    and, correspondingly, any $h$ sufficiently large, there holds
	\begin{equation}\label{eq:comp3}
		\mathcal{L}^2\left( \cup_{\mathcal{F}^{(s)}_h} B \right) \leq 
		\frac{2\pi\xi}{\eta} \rho \mathcal{H}^1\left(J_{u_h}\right),		
	\end{equation}
	where $\xi$ is the same universal constant as in Theorem~\ref{thm:approx-sbv-px}. 
    By~\eqref{eq:equibdd} and~\eqref{eq:comp3}, it 
	follows that $w^{(s)} = w^{(t)}$ $\mathcal{L}^2$-a.e. on $B_{s\rho}$ if 
	$1/2 \leq s \leq t < 1$. Thus, we can define a limit function $\overline{u}$ on 
	$B_\rho$ by letting $\overline{u} = w^{(s)}$ in $B_{s\rho}$ for 
	all $s \in [1/2,1)$. By definition, it is clear the 
    $\overline{u} \in W^{1,p(\cdot)}_{\rm loc}\left(B_\rho,\R^k\right)$. 
    Moreover, for any $s \in [1/2,1)$,
    \begin{equation}\label{eq:bound-bar-u-1}
    \begin{split}
        \int_{B_{s\rho}} \abs{\nabla \overline{u}}^{p(x)}\,{\d}x  &= \int_{B_{s\rho}} \abs{\nabla w^{(s)}}^{p(x)}\,{\d}x \leq \liminf_{h\to+\infty}\int_{B_{s\rho}} \abs{ \nabla\left(w_h^{(s)} \chi_{{B_\rho} \setminus \cup_{\mathcal{F}_h^{(s)}} B }\right)}^{p(x)}\,{\d}x \\
        &\leq \liminf_{h\to+\infty}\int_{B_{\rho} \setminus \cup_{\mathcal{F}_h^{(s)}} B } \abs{ \nabla w_h^{(s)} }^{p(x)}\,{\d}x  
        \leq \liminf_{h\to+\infty} \int_{B_\rho} \abs{\nabla u_h}^{p(x)}\,{\d}x.
    \end{split}
    \end{equation}
    Indeed, the first inequality follows from the weak convergence $w^{(s)}_h - m^{(1/2)}_h$ to 
    $w^{(s)}$ in $B_{s\rho}$, 
    the lower semicontinuity of the modular with respect to the weak convergence 
    (Proposition~\ref{prop:convergence-modulars}), 
    and~\eqref{eq:comp3}. The second inequality follows because $w^{(s)}_h = u_h$ a.e. outside 
    $\cup_{\mathcal{F}^{(s)}_h} B$. Thus, we can let $s \uparrow 1$ in~\eqref{eq:bound-bar-u-1} and 
    obtain
    \begin{equation}\label{eq:bound-bar-u-2}
        \int_{B_{\rho}} \abs{\nabla \overline{u}}^{p(x)}\,{\d}x \leq 
        \liminf_{h\to+\infty} \int_{B_\rho} \abs{\nabla u_h}^{p(x)}\,{\d}x.
    \end{equation} 
    In particular, $\overline{u} \in W^{1,p(\cdot)}\left(B_\rho, \R^k\right)$.
	
	By~\eqref{eq:equibdd} and~\eqref{eq:comp2}, we can find  
    a (not relabelled) subsequence and a vector $\overline{m} \in \R^k$ such that
    $m_h^{(1/2)} \to \overline{m}$. 
    Furthermore, 
    by {(ii)} in Theorem~\ref{thm:approx-sbv-px} and since 
	$\mathcal{L}^2\left( \cup_{\mathcal{F}^{(s)}_{h}} B\right)$ is 
	infinitesimal, 
    $\left\{u_{h} - m_h^{(1/2)} \right\}$ converges in measure 
    to $\overline{u}$, hence also pointwise a.e., up to extraction of a 
    (not relabelled) subsequence. But then $u := \overline{u} + \overline{m}$ belongs 
    to $W^{1,p(\cdot)}\left(B_\rho,\R^k\right)$, and we have $u_{h} \to u$ 
    almost everywhere in $B_\rho$ as $h \to +\infty$. 
    In turn, this implies $u(x) \in \mathcal{M}$ for a.e. $x \in B_{\rho}$. 
    Moreover, by dominated convergence, $u_{h} \to u$ in $L^{p(\cdot)}(B_\rho)$
    as $h \to +\infty$. Finally, \eqref{eq:bound-u} follows 
    from~\eqref{eq:bound-bar-u-2} because $u$ and $\overline{u}$ differ only 
    by a constant.
\end{proof}

\paragraph{Acknowledgements}
The authors are supported by the project \textsc{Star Plus 2020 - Linea 1
(21-UNINA-EPIG-172)} ``New perspectives in the Variational modelling of Continuum Mechanics''. The authors thank the Hausdorff research Institute for Mathematics (HIM)
for the warm hospitality during the Trimester Program “Mathematics of complex materials”,
funded by the Deutsche Forschungsgemeinschaft (DFG, German Research Foundation) under
Germany’s Excellence Strategy – EXC-2047/1 – 390685813, when part of this work was carried
out. The authors are members of, and have been partially supported by, GNAMPA-INdAM.

\bibliographystyle{plain}
\bibliography{SBVpx}

\Addresses
\end{document}